\documentclass[12pt]{article}

\usepackage{amsmath, amsfonts, amssymb, amsthm, a4wide}

\newcommand{\End}{{\rm{End}\ts}}
\newcommand{\non}{\nonumber}
\newcommand{\wt}{\widetilde}
\newcommand{\wh}{\widehat}
\newcommand{\ot}{\otimes}
\newcommand{\la}{\lambda}
\newcommand{\al}{\alpha}
\newcommand{\be}{\beta}
\newcommand{\ze}{\zeta}
\newcommand{\ga}{\gamma}
\newcommand{\si}{\sigma}
\newcommand{\de}{\delta^{}}
\newcommand{\om}{\omega^{}}
\newcommand{\Om}{\Omega}
\newcommand{\vs}{\varsigma}

\newcommand{\ve}{\varepsilon}
\newcommand{\ts}{\,}
\newcommand{\qin}{q^{-1}}

\newcommand{\zin}{z^{-1}}
\newcommand{\tss}{\hspace{1pt}}
\newcommand{\U}{ {\rm U}}
\newcommand{\Y}{ {\rm Y}}
\newcommand{\CC}{\mathbb{C}\tss}
\newcommand{\Pc}{\mathcal{P}}

\newcommand{\RR}{\mathbb{R}}
\newcommand{\ZZ}{\mathbb{Z}}
\newcommand{\gl}{\mathfrak{gl}}
\newcommand{\oa}{\mathfrak{o}}
\newcommand{\spa}{\mathfrak{sp}}
\newcommand{\g}{\mathfrak{g}}

\newcommand{\sll}{\mathfrak{sl}}
\newcommand{\agot}{\mathfrak{a}}

\newcommand{\Fand}{\qquad\text{and}\qquad}

\newcommand{\beq}{\begin{equation}}
\newcommand{\eeq}{\end{equation}}
\newcommand{\ben}{\begin{equation*}}
\newcommand{\een}{\end{equation*}}
\newcommand{\bal}{\begin{aligned}}
\newcommand{\eal}{\end{aligned}}

\newcommand{\bpf}{\begin{proof}}
\newcommand{\epf}{\end{proof}}

\def\beql#1{\begin{equation}\label{#1}}

\newtheorem{thm}{Theorem}[section]
\newtheorem{lemma}[thm]{Lemma}
\newtheorem{prop}[thm]{Proposition}
\newtheorem{cor}[thm]{Corollary}
\newtheorem{conj}[thm]{Conjecture}

\theoremstyle{definition}
\newtheorem{definition}[thm]{Definition}
\newtheorem{example}[thm]{Example}

\theoremstyle{remark}
\newtheorem{remark}[thm]{Remark}

\newcommand{\bth}{\begin{thm}}
\renewcommand{\eth}{\end{thm}}
\newcommand{\bpr}{\begin{prop}}
\newcommand{\epr}{\end{prop}}
\newcommand{\ble}{\begin{lemma}}
\newcommand{\ele}{\end{lemma}}
\newcommand{\bco}{\begin{cor}}
\newcommand{\eco}{\end{cor}}
\newcommand{\bconj}{\begin{conj}}
\newcommand{\econj}{\end{conj}}
\newcommand{\bex}{\begin{example}}
\newcommand{\eex}{\end{example}}
\newcommand{\bde}{\begin{definition}}
\newcommand{\ede}{\end{definition}}
\newcommand{\bre}{\begin{remark}}
\newcommand{\ere}{\end{remark}}

\begin{document}
\title{Representations of twisted $q$-Yangians}
\author{by Lucy Gow and Alexander Molev}

\date{}
\maketitle

\begin{abstract}
The twisted $q$-Yangians are coideal subalgebras
of the quantum affine algebra associated with
$\gl_N$. We prove a classification theorem
for finite-dimensional
irreducible representations of the twisted $q$-Yangians
associated with the symplectic Lie algebras $\spa_{2n}$.
The representations
are parameterized by their highest weights or by
their Drinfeld polynomials. In the simplest case of $\spa_2$
we give an explicit description of all the representations
as tensor products of evaluation modules.
We prove analogues of the Poincar\'e--Birkhoff--Witt theorem
for the quantum affine algebra and for the
twisted $q$-Yangians. We also reproduce a proof of the
classification theorem for finite-dimensional
irreducible representations of the quantum affine algebra
by relying on its $R$-matrix presentation.
\end{abstract}

\tableofcontents

\newpage

\section{Introduction}

The Yangian $\Y(\agot)$
and quantum affine algebra $\U_q(\wh\agot)$
associated with a simple
Lie algebra $\agot$ are known as `infinite-dimensional
quantum groups'. They are deformations of the universal
enveloping algebras $\U(\agot[z])$ and $\U(\wh\agot)$,
respectively, in the class of
Hopf algebras, and were introduced by Drinfeld~\cite{d:ha}
and Jimbo~\cite{j:qd}.
Here $\agot[z]$ denotes the Lie algebra of polynomials
in a variable $z$ with coefficients in $\agot$, while
$\wh\agot$ is the affine Kac--Moody algebra,
i.e., a central extension of the Lie algebra
$\agot[z,z^{-1}]$ of Laurent polynomials in $z$.

The case of $\agot=\sll_N$ (the $A$ type) is exceptional in the sense
that only in this case
do there exist epimorphisms $\Y(\agot)\to\U(\agot)$
and $\U_q(\wh\agot)\to\U_q(\agot)$,
called the evaluation homomorphisms,
where $\U_q(\agot)$
is the corresponding quantized enveloping algebra.
These epimorphisms have important applications in the
representation theory of both the finite- and infinite-dimensional
quantum groups.
For the classical Lie algebra $\agot$ (of type $B$, $C$ or $D$)
there are `twisted' analogues of the
Yangian and quantum affine algebra for which the corresponding
epimorphisms do exist. Namely, the twisted Yangians $\Y'(\oa_N)$
and $\Y'(\spa_{2n})$ associated with the orthogonal
and symplectic Lie algebras were introduced
by Olshanski~\cite{o:ty},
while their $q$-analogues $\Y'_q(\oa_N)$ and $\Y'_q(\spa_{2n})$,
called the twisted $q$-Yangians,
appeared in Molev, Ragoucy and Sorba~\cite{mrs:cs}.
These algebras do not possess natural Hopf algebra structures,
but they are coideal subalgebras of the $A$ type Yangian
and quantum affine algebra, respectively. The evaluation
homomorphisms have the form
\ben
\Y'(\g_N)\to\U(\g_N),\qquad \Y'_q(\g_N)\to\U'_q(\g_N),
\een
where $\g_N$ denotes either the orthogonal Lie algebra $\oa_N$
or the symplectic Lie algebra $\spa_N$ (the latter with $N=2n$)
and $\U'_q(\g_N)$ is the twisted (or nonstandard) quantized
enveloping algebra associated with $\g_N$ which was
defined in \cite{gk:qd},
\cite{nr:qg} and \cite{n:ms}.

Finite-dimensional irreducible representations
of the Yangians $\Y(\agot)$ were classified by Drinfeld~\cite{d:nr}.
The particular case $\agot=\sll_2$ plays a key role
in the arguments and it was done earlier by Tarasov~\cite{t:sq, t:im};
see \cite[Ch.~3]{m:yc} for a detailed exposition of these results.
The classification theorem for the representations of
the quantum affine algebras was proved by
Chari and Pressley~\cite{cp:qaa}, \cite[Ch.~12]{cp:gq}.
Again, the case of $\U_q(\wh\sll_2)$ is crucial, and it is possible
to prove the theorem here
following Tarasov's arguments~\cite{t:sq, t:im}.
The corresponding proof was also outlined in
\cite[Sec. 3.5]{m:yc} and we give more details
below (Section~\ref{sec:qaff}),
as the same approach will be used for the twisted $q$-Yangians.

A classification of
finite-dimensional irreducible representations
of the twisted Yangians $\Y(\oa_N)$ and $\Y(\spa_{2n})$
was obtained in \cite{m:fd}; see also \cite{m:yc}
for a detailed exposition, more references and
applications to representation theory of the classical Lie algebras.
Recent renewed interest in the representation theory of Yangians and
twisted Yangians was caused by its surprising connection
with the theory of finite $W$-algebras (see
\cite{b:ty}, \cite{bk:sy}, \cite{r:ty}, \cite{rs:yr})
and by a generalized Howe duality
(see \cite{kn:yam1}, \cite{kn:ty}).
Note also the applications
of the twisted Yangians and their $q$-analogues to the
soliton spin chain models with special boundary conditions
(see \cite{aacdfr:gb}, \cite{acdfr:sb}) and to
the Teichm\"{u}ller theory of hyperbolic Riemann surfaces
\cite{cm:id}.

In this paper we prove a classification theorem for
finite-dimensional irreducible representations
of the twisted $q$-Yangians associated with the
symplectic Lie algebras $\spa_{2n}$.
The results and the arguments turn out to be parallel to both
the twisted Yangians and quantum affine algebras; cf.
\cite[Ch.~12]{cp:gq} and
\cite[Sec. 3.5 and 4.3]{m:yc}.
First we prove that every finite-dimensional irreducible representation
of $\Y'_q(\spa_{2n})$ is a highest weight representation.
Then we give necessary and sufficient conditions
on the highest weight representations to be
finite-dimensional.
These conditions involve
a family of polynomials $P_1(u),\dots,P_n(u)$ in $u$
(analogues of the Drinfeld polynomials)
so that the finite-dimensional
irreducible representations are essentially parameterized
by $n$-tuples $\big(P_1(u),\dots,P_n(u)\big)$.
In the case of $\Y'_q(\spa_{2})$ we
give an explicit construction
of all finite-dimensional irreducible representations
as tensor products of the evaluation modules over $\U_q(\wh\gl_2)$.

An important ingredient in our arguments is
the Poincar\'e--Birkhoff--Witt theorem for the
quantum affine algebra $\U_q(\wh\gl_N)$ in its
$RTT$-presentation, where $q$
is a fixed nonzero complex number;
see Corollaries~\ref{cor:pbwaffspe}
and \ref{cor:pbwaffrep} below.
This allows us to derive
a new proof of this theorem for the twisted $q$-Yangians
$\Y'_q(\oa_{N})$ and $\Y'_q(\spa_{2n})$; cf. \cite{mrs:cs}.
A version of the PBW theorem in terms of
the `new realization' of
the quantum affine algebra $\U_q(\wh\agot)$
over the field of rational functions in $q$
was given by Beck~\cite{b:bg}, with the case
of $\wh\sll_2$ previously done by Damiani~\cite{d:bt}.
An integral PBW basis of $\U_q(\wh\agot)$ was
constructed by Beck, Chari and Pressley~\cite{bcp:ac},
thus providing a PBW basis of $\U_q(\wh\agot)$, where
$q$ is specialized to any nonzero
complex number;
see also Hernandez~\cite{h:rq} for a weak version
of the PBW theorem for the quantum affinizations of
symmetrized quantum Kac--Moody algebras.
The existence of PBW type bases
for the quantum affine algebras follows also from the general
results of Kharchenko~\cite{h:pbw}.
Although the $RTT$-presentation and the
new realization of the algebra $\U_q(\wh\gl_N)$
are related by the Ding--Frenkel isomorphism~\cite{df:it},
the results of \cite{bcp:ac} do not
immediately imply
a PBW-type theorem for the $RTT$-presentation.

The general approach of this paper
developed for the $C$ type twisted $q$-Yangians
should be applicable
to the $B$ and $D$ types as well,
although some additional arguments
will be needed in order to obtain analogous classification
theorems for representations of the algebras
$\Y'_q(\oa_{N})$; cf. \cite{ik:ct}
and \cite{m:rtq}.


We are grateful to David Hernandez and Mikhail Kotchetov
for discussions of the Poincar\'e--Birkhoff--Witt theorem
for the quantum affine algebras and to Vyjayanthi Chari
for a valuable comment.

\section{Poincar\'e--Birkhoff--Witt theorem}
\setcounter{equation}{0}

We start by reviewing and proving analogues of the PBW theorem
for some quantum algebras. In particular, we prove
it for the $RTT$ presentation of the
quantum affine algebra $\U_q(\wh\gl_N)$ and then use it to get
a new proof of the theorem for the twisted $q$-Yangians.

\subsection{Quantized enveloping algebra $\U_q(\gl_N)$
and its representations}
\label{subsec:que}

Fix a nonzero complex number $q$.
Following \cite{j:qu} and \cite{rtf:ql}, consider
the $R$-matrix presentation of the
quantized enveloping algebra $\U_q(\gl_N)$.
The $R$-matrix is given by
\beql{rmatrixc}
R=q\ts\sum_i E_{ii}\ot E_{ii}+\sum_{i\ne j} E_{ii}\ot E_{jj}+
(q-\qin)\sum_{i<j}E_{ij}\ot E_{ji}
\eeq
which is an element of $\End\CC^N\ot \End\CC^N$, where
the $E_{ij}$ denote the standard matrix units and the indices run over
the set $\{1,\dots,N\}$. The $R$-matrix satisfies
the Yang--Baxter equation
\beql{YBEconst}
R_{12}\ts R_{13}\ts  R_{23} =  R_{23}\ts  R_{13}\ts  R_{12},
\eeq
where both sides take values in
$\End\CC^N\ot \End\CC^N\ot \End\CC^N$ and
the subscripts indicate the copies of $\End\CC^N$, e.g.,
$R_{12}=R\ot 1$ etc.

The algebra $\U_q(\gl_N)$ is generated
by elements $t_{ij}$ and $\bar t_{ij}$ with $1\leqslant i,j\leqslant N$
subject to the relations
\beql{defrel}
\bal
t_{ij}&=\bar t_{ji}=0, \qquad 1 \leqslant i<j\leqslant N,\\
t_{ii}\ts \bar t_{ii}&=\bar t_{ii}\ts t_{ii}=1,
\qquad 1\leqslant i\leqslant N,\\
R\ts T_1T_2&=T_2T_1R,\qquad R\ts \overline T_1\overline T_2=
\overline T_2\overline T_1R,\qquad
R\ts \overline T_1T_2=T_2\overline T_1R.
\eal
\eeq
Here $T$ and $\overline T$ are the matrices
\beql{matrt}
T=\sum_{i,j}t_{ij}\ot E_{ij},\qquad \overline T=\sum_{i,j}
\overline t_{ij}\ot E_{ij},
\eeq
which are regarded as elements of the algebra $\U_q(\gl_N)\ot \End\CC^N$.
Both sides of each of the $R$-matrix relations in \eqref{defrel}
are elements of $\U_q(\gl_N)\ot \End\CC^N\ot \End\CC^N$ and the
subscripts
of $T$ and $\overline T$ indicate the copies of $\End\CC^N$ where
$T$ or $\overline T$ acts; e.g. $T_1=T\ot 1$. In terms of the
generators the defining relations between the $t_{ij}$
can be written as
\beql{defrelg}
q^{\delta_{ij}}\ts t_{ia}\ts t_{jb}-
q^{\delta_{ab}}\ts t_{jb}\ts t_{ia}
=(q-\qin)\ts (\de_{b<a} -\de_{i<j})
\ts t_{ja}\ts t_{ib}
\eeq
where $\de_{i<j}$ equals $1$ if $i<j$ and $0$ otherwise.
The relations between the $\bar t_{ij}$
are obtained by replacing $t_{ij}$ by $\bar t_{ij}$ everywhere in
\eqref{defrelg}:
\beql{defrelgbar}
q^{\delta_{ij}}\ts \bar t_{ia}\ts \bar t_{jb}-
q^{\delta_{ab}}\ts \bar t_{jb}\ts \bar t_{ia}
=(q-\qin)\ts (\de_{b<a} -\de_{i<j})
\ts \bar t_{ja}\ts \bar t_{ib},
\eeq
while the relations involving both
$t_{ij}$ and $\bar t_{ij}$ have the form
\beql{defrelg2}
q^{\delta_{ij}}\ts \bar t_{ia}\ts t_{jb}-
q^{\delta_{ab}}\ts t_{jb}\ts \bar t_{ia}
=(q-\qin)\ts (\de_{b<a}\ts t_{ja}\ts \bar t_{ib} -\de_{i<j}\ts
\ts \bar t_{ja}\ts t_{ib}).
\eeq
Note that for any nonzero complex number $d$ the mapping
\beql{automd}
t_{ij}\mapsto d\ts t_{ij},\qquad
\bar t_{ij}\mapsto d^{-1}\ts\bar t_{ij}
\eeq
defines an automorphism of the algebra $\U_q(\gl_N)$.

Let $z$ denote an indeterminate.
Introduce the algebra $\U_z(\gl_N)$
over $\CC(z)$ with the generators
$t_{ij}$ and $\bar t_{ij}$ with $1\leqslant i,j\leqslant N$
subject to the relations
given in \eqref{defrel}
with $q$ replaced by $z$. Furthermore, we denote by
$\U^{\circ}_z(\gl_N)$
the algebra defined over the ring of
Laurent polynomials $\CC[z,z^{-1}]$ with the same set
of generators and relations. Then we have the isomorphism
\beql{isom}
\U^{\circ}_z(\gl_N)\ot_{\CC[z,z^{-1}]}\CC\cong \U_q(\gl_N),
\eeq
where the $\CC[z,z^{-1}]$-module $\CC$ is defined
via the evaluation of the Laurent polynomials at $z=q$.

The quantized enveloping algebras admit families
of PBW bases depending on choices of
reduced decompositions
of the longest element of the Weyl group; see
Lusztig~\cite{l:fd}. In the $A$ type such bases
were previously constructed by Rosso~\cite{r:ap}
and Yamane~\cite{y:pbw}.
These constructions use the Drinfeld--Jimbo presentation
of the quantized enveloping algebras. In this presentation,
the algebra $\U_z(\gl_N)$ over $\CC(z)$ is generated by
the elements $t_1,\dots,t_N,t_1^{-1},\dots,t_N^{-1}$,
$e_1,\dots,e_{N-1}$ and
$f_1,\dots,f_{N-1}$ subject to the defining relations
\ben
\bal
t_i\tss t_j=t_j\tss t_i, &\qquad t^{}_i\tss
t_i^{-1}=t_i^{-1}\tss t^{}_i=1, \\
t^{}_i\tss e^{}_j\tss t_i^{-1}=e^{}_j \ts
z^{\tss\delta_{ij}-\delta_{i,j+1}},&\qquad
t^{}_i\tss f^{}_j\tss t_i^{-1}=f^{}_j \ts
z^{-\delta_{ij}+\delta_{i,j+1}},\\
[e_i,f_j]=\delta_{ij}\ts\frac{k^{}_i-k_i^{-1}}
{z-\zin}&\qquad \text{with}\ \
k_i=t^{}_i\tss t_{i+1}^{-1},\\
[e_i,e_j]=[f_i,f_j]=0&\qquad\text{if}\ \ |i-j|>1,\\
e^2_i\tss e^{}_{j}-(z+\zin)\tss e^{}_i\tss e^{}_{j}\tss e^{}_i
&+e^{}_{j}\tss e^2_i=0\qquad\text{if}\ \ |i-j|=1,\\
f^2_i\tss f^{}_{j}-(z+\zin)\tss f^{}_i\tss f^{}_{j}\tss f^{}_i
&+f^{}_{j}\tss f^2_i=0\qquad\text{if}\ \ |i-j|=1.
\eal
\een
The root vectors can be defined inductively by
\begin{alignat}{2}
{}&e_{i,i+1}=e_i,\qquad {}&&e_{i+1,i}=f_i,
\non\\
\label{rootv}
{}&e_{ij}=e_{ip}\tss e_{pj}-z\ts e_{pj}\tss e_{ip}
\qquad &&\text{for}\ \ i<p<j,\\
{}&e_{ij}=e_{ip}\tss e_{pj}-\zin\tss e_{pj}\tss
e_{ip}\qquad &&\text{for}\ \ i>p>j,
\non
\end{alignat}
and the elements $e_{ij}$ are independent of the choice of values
of the index $p$.

An isomorphism between the two presentations
of $\U_z(\gl_N)$ is given by the formulas
\beql{isomgln}
t_{ii}\mapsto t_i,\qquad \bar t_{ii}\mapsto t^{-1}_i,\qquad
\bar t_{ij}\mapsto -(z-\zin)\tss e_{ij}\ts t_i^{-1},\qquad
t_{ji}\mapsto (z-\zin)\tss t_i\ts e_{ji}
\end{equation}
for $i<j$; see \cite{df:it}, \cite{rtf:ql}.
We shall identify the corresponding elements of $\U_z(\gl_N)$
via this isomorphism.

The quantized enveloping algebra $\U_z(\sll_N)$ can be
defined as the $\CC(z)$-subalgebra of $\U_z(\gl_N)$ generated
by the elements $k_i,k_i^{-1},e_i,f_i$
for $i=1,\dots,N-1$. Similarly, if $q$ is a nonzero
complex number such that $q^2\ne 1$, then $\U_q(\sll_N)$
can be defined as the subalgebra of $\U_q(\gl_N)$
generated by the same elements.

We will be using the following form of the PBW theorem
for the quantized enveloping algebra
associated with $\gl_N$.

\bpr\label{prop:pbwfin}
The monomials
\begin{multline}\label{monom}
t^{k_{N,N-1}}_{N,N-1}\ts t^{k_{N,N-2}}_{N,N-2}\ts
t^{k_{N-1,N-2}}_{N-1,N-2}
\dots\ts t^{k_{N2}}_{N2}
\dots\ts t^{k_{32}}_{32}\ts t^{k_{N1}}_{N1}\ts
\dots\ts t^{k_{21}}_{21}\\
{}\times t_{11}^{l_1}\dots t_{NN}^{l_N}\ts
\bar t^{\ts k_{12}}_{12}\dots\ts \bar t^{\ts k_{1N}}_{1N}\ts
\bar t^{\ts k_{23}}_{23}\dots\ts \bar t^{\ts k_{2N}}_{2N}\dots\ts
\bar t^{\ts k_{N-1,N}}_{N-1,N},
\end{multline}
where the $k_{ij}$ run over non-negative integers
and the $l_i$ run over all integers,
form a basis of the $\CC[z,\zin]$-algebra $\U^{\circ}_z(\gl_N)$.
\epr

\bpf
It follows easily from the defining relations of
$\U^{\circ}_z(\gl_N)$ that the monomials span the algebra
over $\CC[z,\zin]$. Suppose now that there is a nontrivial
linear combination of the monomials \eqref{monom}
with coefficients in $\CC[z,\zin]$ equal to zero.
Applying the isomorphism \eqref{isomgln} and the relations
$t_{i}\ts e_{jb}=z^{\de_{ij}-\de_{ib}}\ts e_{jb}\ts t_{i}$
we then obtain a nontrivial linear combination
over $\CC[z,\zin]$ of the monomials
\begin{multline}\label{monomdj}
e^{k_{N,N-1}}_{N,N-1}\ts e^{k_{N,N-2}}_{N,N-2}\ts
e^{k_{N-1,N-2}}_{N-1,N-2}
\dots\ts e^{k_{N2}}_{N2}
\dots\ts e^{k_{32}}_{32}\ts e^{k_{N1}}_{N1}\ts
\dots\ts e^{k_{21}}_{21}\\
{}\times t_1^{\ts l_1}\dots\ts t_N^{\ts l_N}\ts
e^{k_{12}}_{12}\dots\ts e^{k_{1N}}_{1N}\ts
e^{k_{23}}_{23}\dots\ts e^{k_{2N}}_{2N}\dots\ts
e^{k_{N-1,N}}_{N-1,N}
\end{multline}
equal to zero. Here the $k_{ij}$ run over non-negative integers
and the $l_i$ run over all integers. However, by
the PBW theorem for the Drinfeld--Jimbo
presentation of the algebra $\U_z(\gl_N)$
(see \cite{l:fd}, \cite{r:ap}, \cite{y:pbw}),
the monomials \eqref{monomdj} form a basis of
$\U_z(\gl_N)$ over $\CC(z)$. This makes a contradiction.
\epf

The following corollary is immediate from the
isomorphism \eqref{isom}.

\bco\label{cor:pbwfin}
Let $q$ be a nonzero complex number. Then the monomials
\eqref{monom} form a basis of $\U_q(\gl_N)$ over $\CC$.
\qed
\eco

Note that in the particular case $q=1$ the algebra
$\U_1(\gl_N)$ is commutative.
Using Corollary~\ref{cor:pbwfin}
we will identify it
with the algebra of polynomials $\Pc_N$
in the variables $\bar x_{ij},x_{ji}$ with
$1\leqslant i< j\leqslant N$ and $x_{ii},\bar x_{ii}$
with $i=1,\dots,N$ subject to the relations $x_{ii}\bar x_{ii}=1$
for all $i$. Thus, due to \eqref{isom} we have the isomorphism
\beql{isomqfone}
\U^{\circ}_z(\gl_N)\ot_{\CC[z,z^{-1}]}\CC\cong \Pc_N,
\eeq
where the $\CC[z,z^{-1}]$-module $\CC$ is defined
via the evaluation of the Laurent polynomials at $z=1$.

In the other degenerate case
$q=-1$ the algebra $\U_{-1}(\gl_N)$
is essentially a `quasi-polynomial' algebra;
see e.g. \cite[Sec.~1.8]{dck:rq}.
It is well known that quasi-polynomial algebras admit
PBW bases.

We will also use an extended version of the
quantized enveloping algebra considered in \cite{m:rtq}.
Denote by $\U^{\tss\text{\rm ext}}_z(\gl_N)$
the algebra over $\CC[z,z^{-1}]$ generated
by elements $t_{ij}$ and $\bar t_{ij}$
with $1\leqslant i,j\leqslant N$
and elements $t^{-1}_{ii}$ and $\bar t_{ii}^{\ts\ts-1}$
with $1\leqslant i\leqslant N$
subject to the relations
\beql{defrelext}
\bal
t_{ij}&=\bar t_{ji}=0, \qquad 1 \leqslant i<j\leqslant N,\\
t_{ii}\ts \bar t_{ii}&=\bar t_{ii}\ts t_{ii},\qquad
t^{}_{ii}\ts t^{-1}_{ii}=t^{-1}_{ii}\ts t^{}_{ii}=1,\qquad
\bar t^{}_{ii}\ts \bar t_{ii}^{\ts\ts-1}=
\bar t^{\ts\ts-1}_{ii}\ts \bar t_{ii}^{}=1,
\qquad 1\leqslant i\leqslant N,\\
R\ts T_1T_2&=T_2T_1R,\qquad R\ts \overline T_1\overline T_2=
\overline T_2\overline T_1R,\qquad
R\ts \overline T_1T_2=T_2\overline T_1R,
\eal
\eeq
where we use the notation of \eqref{defrel} with $q$ replaced by $z$
in the definition of $R$. Although we use
the same notation for the generators
of the algebras $\U^{\tss\text{\rm ext}}_z(\gl_N)$
and $\U^{\circ}_z(\gl_N)$, it should always be clear
from the context which algebra is considered at
any time. There is a natural epimorphism
$\U^{\tss\text{\rm ext}}_z(\gl_N)\to\U^{\circ}_z(\gl_N)$
which takes the generators $t_{ij}$ and $\bar t_{ij}$
to the elements
with the same name.
The kernel $K$ of this epimorphism
is the two-sided ideal of
the algebra $\U^{\tss\text{\rm ext}}_z(\gl_N)$
generated by the elements $t_{ii}\tss \bar t_{ii}-1$
for $i=1,\dots,N$. All these
elements are central in this algebra and
we have the isomorphism
$\U^{\tss\text{\rm ext}}_z(\gl_N)/K\cong \U^{\circ}_z(\gl_N)$.

The following analogue of the PBW theorem is implied
by Proposition~\ref{prop:pbwfin}; see also \cite{m:rtq}.

\bpr\label{prop:pbwfinext}
The monomials
\begin{multline}\label{monomext}
t^{k_{N,N-1}}_{N,N-1}\ts t^{k_{N,N-2}}_{N,N-2}\ts
t^{k_{N-1,N-2}}_{N-1,N-2}
\dots\ts t^{k_{N2}}_{N2}
\dots\ts t^{k_{32}}_{32}\ts t^{k_{N1}}_{N1}\ts
\dots\ts t^{k_{21}}_{21}\ts
t_{11}^{l_1}\dots t_{NN}^{l_N}\\
{}\times
\bar t_{11}^{\ts m_1}\dots \bar t_{NN}^{\ts m_N}\ts
\bar t^{\ts k_{12}}_{12}\dots\ts \bar t^{\ts k_{1N}}_{1N}\ts
\bar t^{\ts k_{23}}_{23}\dots\ts \bar t^{\ts k_{2N}}_{2N}\dots\ts
\bar t^{\ts k_{N-1,N}}_{N-1,N},
\end{multline}
where the $k_{ij}$ run over non-negative integers
and the $l_i$ and $m_i$ run over all integers,
form a basis of the $\CC[z,\zin]$-algebra
$\U^{\tss\text{\rm ext}}_z(\gl_N)$.
\qed
\epr

By specializing $z$ to a nonzero complex number $q$
in the definition of $\U^{\tss\text{\rm ext}}_z(\gl_N)$
we obtain an algebra $\U^{\tss\text{\rm ext}}_q(\gl_N)$
over $\CC$ defined by the same set of relations \eqref{defrelext}.
So we have the isomorphism
\beql{isomext}
\U^{\tss\text{\rm ext}}_z(\gl_N)
\ot_{\CC[z,z^{-1}]}\CC\cong \U^{\tss\text{\rm ext}}_q(\gl_N),
\eeq
where the $\CC[z,z^{-1}]$-module $\CC$ is defined
via the evaluation of the Laurent polynomials at $z=q$.
The corresponding monomials \eqref{monomext} form a basis
of $\U^{\tss\text{\rm ext}}_q(\gl_N)$.
In the particular case $q=1$ the algebra
$\U^{\tss\text{\rm ext}}_1(\gl_N)$ can be identified
with the algebra of polynomials $\Pc^{\tss\text{\rm ext}}_N$
in the variables $\bar x_{ij},x_{ji}$ with
$1\leqslant i< j\leqslant N$ and
$x^{}_{ii}, x^{-1}_{ii}, \bar x^{}_{ii}, \bar x_{ii}^{-1}$
with $i=1,\dots,N$.
Thus we have the isomorphism
\beql{isomqfoneext}
\U^{\tss\text{\rm ext}}_z(\gl_N)
\ot_{\CC[z,z^{-1}]}\CC\cong \Pc^{\tss\text{\rm ext}}_N,
\eeq
where the $\CC[z,z^{-1}]$-module $\CC$ is defined
via the evaluation of the Laurent polynomials at $z=1$.

Suppose now that $q$ is a nonzero complex number which
is not a root of unity.
A description of finite-dimensional irreducible
representations of the algebra $\U^{\tss\text{\rm ext}}_q(\gl_N)$
can be easily obtained from the corresponding results
for the algebras $\U_q(\gl_N)$ and $\U_q(\sll_N)$; see e.g.
\cite[Ch.~10]{cp:gq}.
A representation $L$ of $\U^{\tss\text{\rm ext}}_q(\gl_N)$ is called
a {\it highest weight representation\/} if $L$ is generated
by a nonzero vector $\ze$ (the {\it highest vector\/})
such that
\begin{alignat}{2}
\bar t_{ij}\ts\ze&=0 \qquad &&\text{for} \qquad
1\leqslant i<j\leqslant N,\qquad \text{and}
\non\\
t_{ii}\ts\ze&=\la_i\tss\ze, \qquad
&&\bar t_{ii}\ts\ze=\bar\la_i\tss\ze,\qquad
\text{for} \quad 1\leqslant i\leqslant N,
\non
\end{alignat}
for some nonzero complex numbers $\la_i$ and $\bar\la_i$.
The tuple $(\la_1,\dots,\la_{N};\bar\la_1,\dots,\bar\la_N)$ is called
the {\it highest weight\/} of $L$.
Due to Proposition~\ref{prop:pbwfinext}, for
any $N$-tuples of nonzero complex numbers
$\la=(\la_1,\dots,\la_{N})$ and $\bar\la=(\bar\la_1,\dots,\bar\la_N)$,
there exists an irreducible highest representation
$L(\la;\bar\la)$ with the highest weight $(\la;\bar\la)$.
This representation can be defined as a quotient
of the corresponding Verma module in a standard way.

The irreducible highest weight representations $L(\mu)$,
$\mu=(\mu_1,\dots,\mu_N)$,
over the algebra $\U_q(\gl_N)$ are defined in a similar way
with the above conditions on the highest vector
replaced by
\begin{alignat}{2}
\bar t_{ij}\ts\ze&=0 \qquad &&\text{for} \qquad
1\leqslant i<j\leqslant N,\qquad \text{and}
\non\\
t_{ii}\ts\ze&=\mu_i\tss\ze, \qquad
&&
\text{for} \qquad 1\leqslant i\leqslant N.
\non
\end{alignat}
The representation $L(\mu)$ is finite-dimensional
if and only if there exist nonnegative integers $m_i$ satisfying
$m_1\geqslant m_2\geqslant\dots\geqslant m_N$,
elements $\ve_i\in\{-1,1\}$ for $i=1,\dots,N$,
and a nonzero complex number $d$ such that
\ben
\mu_i=d\ts \ve_i\ts q^{m_i},\qquad i=1,\dots,N.
\een

\bpr\label{prop:fdext}
Every finite-dimensional irreducible representation
of $\U^{\tss\text{\rm ext}}_q(\gl_N)$ is isomorphic
to a highest weight representation $L(\la;\bar\la)$ such
that
\ben
\la_i-a\ts q^{2m_i}\bar\la^{}_i=0,\qquad i=1,\dots,N,
\een
for
some nonnegative integers $m_i$ satisfying
$m_1\geqslant m_2\geqslant\dots\geqslant m_N$
and a nonzero complex number $a$.
\epr

\bpf
A standard argument shows that
every finite-dimensional irreducible representation
of $\U^{\tss\text{\rm ext}}_q(\gl_N)$ is isomorphic
to a highest weight representation; cf. \cite[Ch.~10]{cp:gq}.
Hence we only need to determine when the representation
$L(\la;\bar\la)$ is finite-dimensional.
Each central element
$t_{ii}\tss \bar t_{ii}$ of $\U^{\tss\text{\rm ext}}_q(\gl_N)$ acts
on $L(\la;\bar\la)$ as multiplication
by the scalar $\la_i\bar\la_i$. Fix constants
$c_1,\dots,c_N$ such that
\ben
c^2_{i}=\la_i\bar\la_i,\qquad i=1,\dots,N.
\een
Then the mapping
\ben
t_{ij}\mapsto c_i\tss t_{ij},\qquad
\bar t_{ij}\mapsto c_i\tss \bar t_{ij}
\een
defines an epimorphism
$\U^{\tss\text{\rm ext}}_q(\gl_N)\to\U_q(\gl_N)$
whose kernel is generated by the elements
\beql{centiext}
t_{ii}\tss \bar t_{ii}-\la_i\bar\la_i,
\qquad i=1,\dots,N.
\eeq
Hence, identifying $\U_q(\gl_N)$ with the quotient
of $\U^{\tss\text{\rm ext}}_q(\gl_N)$ by this kernel,
we can equip
$L(\la;\bar\la)$ with the structure of an irreducible
highest weight representation of
$\U_q(\gl_N)$. Its
highest weight $(\mu_1,\dots,\mu_{N})$ is given by
\ben
\mu_{i}=c^{-1}_{i}\ts\la_i,\qquad i=1,\dots,N.
\een
This representation is finite-dimensional
if and only if
\ben
c^{-1}_{i}\ts\la_i=d\ts \ve_i\ts q^{m_i}
\een
for some nonnegative integers $m_i$ satisfying
$m_1\geqslant m_2\geqslant\dots\geqslant m_N$,
a nonzero complex number $d$,
and some elements $\ve_i\in\{-1,1\}$ for $i=1,\dots,N$.
By our choice of the constants $c_i$, this
is equivalent to the relations
$\la_i\ts\bar\la_i^{\ts-1}=a\ts q^{2m_i}$ with $a=d^{\tss2}$,
as required.
\epf

\subsection{Twisted quantized enveloping algebras
$\U'_q(\oa_N)$ and $\U'_q(\spa_{2n})$}
\label{subsec:twrep}

The twisted
quantized enveloping algebra $\U'_q(\oa_N)$
associated with the orthogonal Lie algebra $\oa_N$
was introduced independently in \cite{gk:qd} and
\cite{nr:qg}. Its $R$-matrix presentation was given
in \cite{n:ms}. We follow the notation of \cite{mrs:cs}
and define $\U'_q(\oa_N)$
as the subalgebra of $\U_q(\gl_N)$ generated by the matrix
elements $s_{ij}$ of the matrix $S=T\ts \overline T^{\ts t}$,
where $t$ denotes the usual matrix transposition.
More explicitly, the elements $s_{ij}$ are given by
\beql{sijdeo}
s_{ij}=\sum_{a=1}^N t_{ia}\bar t_{ja}.
\eeq
Hence, \eqref{defrel} implies
\begin{align}\label{sijo}
s_{ij}&=0, \qquad 1 \leqslant i<j\leqslant N,\\
\label{sii1}
s_{ii}&=1,\qquad 1\leqslant i\leqslant N.
\end{align}
Furthermore, $\U'_q(\oa_N)$ is isomorphic to the
algebra with (abstract) generators $s_{ij}$
with the condition $i,j\in\{1,\dots,N\}$
subject to the defining relations \eqref{sijo}, \eqref{sii1} and
\beql{rsrs}
R\ts S_1 R^{\ts t_1} S_2=S_2R^{\ts t_1} S_1R,
\eeq
where
\beql{rt}
R^{\ts t_1}=q\ts\sum_i E_{ii}\ot E_{ii}+\sum_{i\ne j} E_{ii}\ot E_{jj}+
(q-\qin)\sum_{i<j}E_{ji}\ot E_{ji}.
\eeq
In terms of the generators,
relation \eqref{rsrs} takes the form
\beql{drabs}
\bal
q^{\delta_{aj}+\delta_{ij}}\ts s_{ia}\ts s_{jb}-
q^{\delta_{ab}+\delta_{ib}}\ts s_{jb}\ts s_{ia}
{}&=(q-\qin)\ts q^{\delta_{ai}}\ts (\de_{b<a} -\de_{i<j})
\ts s_{ja}\ts s_{ib}\\
{}&+(q-\qin)\ts \big(q^{\delta_{ab}}\ts \de_{b<i}\ts s_{ji}\ts s_{ba}
- q^{\delta_{ij}}\ts \de_{a<j}\ts s_{ij}\ts s_{ab}\big)\\
{}&+ (q-\qin)^2\ts  (\de_{b<a<i} -\de_{a<i<j})\ts s_{ji}\ts s_{ab},
\eal
\end{equation}
where $\de_{i<j}$ or $\de_{i<j<k}$ equals $1$ if
the subscript inequality
is satisfied and $0$ otherwise.

An analogue of the PBW theorem for
the algebra $\U'_q(\oa_N)$ was proved in \cite{ik:nd};
see also \cite{m:rtq}, \cite{mrs:cs}.
Yet another proof is obtained
from Proposition~\ref{prop:pbwfin}.
We regard $q$ as a nonzero complex number.

\bpr\label{prop:pbwo}
The monomials
\beql{monomo}
s^{\ts k_{21}}_{21}\ts
s^{\ts k_{32}}_{32}\ts s^{\ts k_{31}}_{31}\dots\ts
s^{\ts k_{N1}}_{N1}\ts
s^{\ts k_{N2}}_{N2}\dots\ts s^{\ts k_{N,N-1}}_{N,N-1},
\eeq
where the $k_{ij}$ run over non-negative integers
form a basis of the algebra $\U'_q(\oa_N)$.
\epr

\bpf
Let us consider
the $\CC[z,z^{-1}]$-subalgebra $\U^{\circ}_z(\oa_N)$
of $\U^{\circ}_z(\gl_N)$
generated by the elements $s_{ij}$ defined
by \eqref{sijdeo} and show that
the monomials \eqref{monomo} form its basis.
It follows easily from the defining relations that
the monomials span the algebra;
see \cite[Lemma~2.1]{mrs:cs}.
Suppose now that
a nontrivial $\CC[z,z^{-1}]$-linear combination
of the monomials \eqref{monomo}
is zero.
By Proposition~\ref{prop:pbwfin}
we may suppose that at least one coefficient
of the combination does not vanish at $z=1$.
Using the isomorphism \eqref{isomqfone} we then
get a nontrivial $\CC$-linear combination
of the corresponding monomials
in the polynomial algebra $\Pc_N$. We will come
to a contradiction if we show that the images $\si_{ij}$
of the generators of $s_{ij}$, $i>j$,
in $\Pc_N$ are algebraically independent.

We have
\ben
\si_{ij}=\sum_{a=1}^N x_{ia}\ts \bar x_{ja}.
\een
It suffices to verify that the differentials
$d\si_{ij}$ are linearly independent.
Calculate the differentials in terms of
$dx_{ia}$ and
$d\bar x_{ia}$ and specialize the coefficient matrix
by setting
$x_{ij}=\bar x_{ij}=\de_{ij}$.
Then
$
d\si_{ij}=d\bar x_{ji}+dx_{ij}
$
which implies that the differentials $d\si_{ij}$
are linearly independent even under the specialization.

This proves that the monomials \eqref{monomo}
form a basis over $\CC[z,\zin]$
in the subalgebra $\U^{\circ}_z(\oa_N)$.
The application of the isomorphism \eqref{isom}
shows that the monomials \eqref{monomo}
form a basis over $\CC$ in $\U'_q(\oa_N)$.
\epf

Finite-dimensional irreducible representations of
the algebra $\U'_q(\oa_N)$ were classified in \cite{ik:ct}.
Moreover, that paper also contains explicit realization
of the representations of `classical type' via
Gelfand--Tsetlin bases.

\medskip

The twisted
quantized enveloping algebra $\U'_q(\spa_{2n})$
associated with the symplectic Lie algebra $\spa_{2n}$
was first introduced in \cite{n:ms}.
In order to define it, consider the $2n\times 2n$ matrix $G$
given by
\beql{g}
G=q\ts \sum_{k=1}^n E_{2k-1,2k}-\sum_{k=1}^n E_{2k,2k-1},
\eeq
that is,
\ben
G=\begin{bmatrix}
                     \phantom{-}0&q&\cdots&\phantom{-}0&0\ \\
                     -1&0&\cdots&\phantom{-}0&0\ \\
                     \vdots&\vdots&\ddots&\vdots&\vdots\\
                     \phantom{-}0&0&\cdots&\phantom{-}0&q\ \\
                     \phantom{-}0&0&\cdots&-1&0\
\end{bmatrix}.
\een
We define $\U'_q(\spa_{2n})$ as the subalgebra
of $\U^{\tss\text{\rm ext}}_q(\gl_{2n})$
generated by the matrix elements
$s_{ij}$ of the matrix
$S= T\ts G\ts \overline T^{\ts t}$
together with the elements
\beql{sinv}
s^{-1}_{i,i+1}= q^{-1}\ts t^{-1}_{ii}\ts \bar t_{i+1,i+1}^{\ts\ts-1}
\eeq
for $i=1,3,\dots,2n-1$. More explicitly,
\beql{sijtij}
s_{ij}= q\ts \sum_{a=1}^n
t_{i,2a-1}\bar t_{j,2a}
-\sum_{a=1}^n t_{i,2a}\bar t_{j,2a-1}.
\eeq
By \eqref{defrel} we have
\beql{symexcf}
s_{ij}=0\quad\text{for}\quad i<j \quad\text{unless}\ \
j= i+1\quad\text{with}\quad i\ \  \text{odd}.
\eeq
All matrix elements $\bar s_{ij}$ of the matrix
$\overline S=\overline T\ts G\ts T^{\ts t}$ also
belong to the subalgebra.
It was proved in \cite{mrs:cs} that
$\U'_q(\spa_{2n})$ is isomorphic to the
algebra with (abstract) generators $s_{ij}$ with
$i,j\in\{1,\dots,2n\}$ and
$s^{-1}_{i,i+1}$ with $i=1,3,\dots,2n-1$,
subject to the defining relations \eqref{rsrs} (with $N=2n$),
\eqref{symexcf} and
\beql{invrel}
s_{i,i+1}^{}\ts s_{i,i+1}^{-1}=s_{i,i+1}^{-1}\ts s^{}_{i,i+1}=1,\qquad
i=1,3,\dots,2n-1.
\eeq

Our definition of $\U'_q(\spa_{2n})$ follows closer
the original paper \cite{n:ms}, while
a slightly different algebra $\U^{\tss\text{\rm tw}}_q(\spa_{2n})$
was studied in \cite{mrs:cs}. The latter was defined
as a subalgebra of $\U_q(\gl_{2n})$ by the same formulas
\eqref{sinv} and \eqref{sijtij} which lead to extra
relations for the generators: for any odd $i$
\beql{sqcube}
s_{i+1,i+1}\ts s_{ii}-q^2\ts s_{i+1,i}\ts
s_{i,i+1}=q^3.
\eeq
They are implied by the relations $t_{ii}\ts\bar t_{ii}=1$
which hold in the algebra $\U_q(\gl_{2n})$ but not
in $\U^{\tss\text{\rm ext}}_q(\gl_{2n})$.
Moreover,
the elements $s_{i+1,i+1}\ts s_{ii}-q^2\ts s_{i+1,i}\ts
s_{i,i+1}$ are central in the algebra $\U'_q(\spa_{2n})$
and its quotient by the relations \eqref{sqcube}
is isomorphic to $\U^{\tss\text{\rm tw}}_q(\spa_{2n})$;
see also \cite{m:rtq}, where a slightly
different notation was used.
The latter algebra is a deformation
of the universal enveloping algebra $\U(\spa_{2n})$;
see \cite{mrs:cs}.

An analogue of the PBW theorem for
the algebra $\U'_q(\spa_{2n})$ was proved in
\cite{m:rtq}. That paper and \cite{mrs:cs}
also contain proofs of this theorem for the
quotient algebra $\U^{\tss\text{\rm tw}}_q(\spa_{2n})$.
Here we give a
more direct proof based on Proposition~\ref{prop:pbwfinext}
for a slightly different order on the set of generators.

We define a function
$\vs: \{ 1, 2, \ldots, 2n\} \to \{\pm 1, \pm 3,\ldots, \pm(2n-1) \}$
by
\beql{nu}
\vs (i) = \left\{ \begin{array}{ll}
\phantom{ -}i & \quad \text{ if\ \  }i \text{\ \  is odd},\\
-i+1 & \quad \text{ if\ \  }i\text{\ \  is even}.
\end{array} \right.
\eeq
We say $s_{ia} < s_{jb}$ if
$\big(\vs(i)+\vs(a), \vs(i)\big) < \big(\vs(j) + \vs(b), \vs(j)\big)$
when ordered
lexicographically.

In the next proposition we consider
the corresponding ordered monomials in the
generators $s_{ij}$ with $i\geqslant j$ together with
$s_{i,i+1}$ and $s_{i,i+1}^{-1}$ with odd $i$.

\bpr\label{prop:modPBWforU}
The ordered monomials
\beql{monomsp}
s^{\ts k_{2n,2n}}_{2n,2n}\ts
s^{\ts k_{2n,2n-2}}_{2n,2n-2}\ts \dots\ts
s^{\ts k_{2n,2n-1}}_{2n,2n-1}\dots\ts
s^{\ts k_{21}}_{21}\ts
s^{\ts k_{12}}_{12}\ts\dots\ts
s^{\ts k_{2n-1,2n}}_{2n-1,2n}\dots\ts
s^{\ts k_{2n-1,2n-3}}_{2n-1,2n-3}\ts s^{\ts k_{2n_1,2n-1}}_{2n-1,2n-1},
\eeq
where $k_{12},k_{34},\dots, k_{2n-1,2n}$ run over all integers and
the remaining powers $k_{ij}$
run over non-negative integers,
form a basis of the algebra $\U'_q(\spa_{2n})$.
\epr

\bpf
Let us consider
the $\CC[z,z^{-1}]$-subalgebra $\U^{\circ}_z(\spa_{2n})$
of $\U^{\tss\text{\rm ext}}_z(\gl_{2n})$
generated by the elements \eqref{sinv}
and \eqref{sijtij} with $q$ replaced by $z$ and show that
the monomials \eqref{monomsp} form its basis.
The application of the isomorphism \eqref{isomext}
will then imply that the monomials
form a basis over $\CC$ in $\U'_q(\spa_{2n})$.

First we prove that an arbitrary monomial
in the generators can be written as
a linear combination of the ordered monomials;
cf. \cite[Lemma~2.1]{mrs:cs}. Due to the relations
\ben
s^{}_{i,i+1}\ts s^{}_{kl}=
z^{\de_{ik}+\de_{il}-\de_{i+1,k}-\de_{i+1,l}}\ts
s^{}_{kl}\ts s^{}_{i,i+1},\qquad i=1,3,\dots,2n-1,
\een
we can restrict our attention to those
monomials where all generators occur in non-negative powers.
We define the \emph{degree} of a monomial
$s_{i_{1}a_{1}} \dots s_{i_{p}a_{p}}$,
to be $d =i_{1} + \dots + i_{k}$ and we argue
by induction on the degree $d$.
Modulo products of degree less than $i+j$, the
relations \eqref{drabs} (with $q$ replaced by $z$) imply:
\begin{eqnarray}
\nonumber
\lefteqn{z^{\delta_{aj} + \delta_{ij}}
s_{ia} s_{jb} - z^{\delta_{ab} + \delta_{ib}} s_{jb} s_{ia}}\\
\label{Ugradeddr}
&\equiv & (z-z^{-1})\tss z^{\delta_{ai}}
(\delta_{b<a} - \delta_{i<j})\tss s_{ja} s_{ib}.
\end{eqnarray}
Swapping here $i$ with $j$ and $a$ with $b$ we can also
write this in the form
\begin{eqnarray}
\nonumber
\lefteqn{z^{\delta_{aj} + \delta_{ab}}
s_{ia} s_{jb} - z^{\delta_{ij} + \delta_{ib}} s_{jb} s_{ia}}\\
\label{Ugradeddrop}
&\equiv & (z^{-1}-z)\tss z^{\delta_{bj}}
(\delta_{a<b} - \delta_{j<i})\tss s_{ib} s_{ja}.
\end{eqnarray}
Suppose $\vs(i) + \vs(a) > \vs(j) + \vs(b)$.
Then if $\vs(i) + \vs(b) > \vs(j) + \vs(a)$, the equation
\eqref{Ugradeddr} allows us to write $s_{ia} s_{jb}$
as a linear combination of ordered monomials
and monomials of lower
degree.
On the other hand, if $\vs(i) + \vs(b) < \vs(j) + \vs(a)$
then the same outcome is achieved by using \eqref{Ugradeddrop}.
In the case $\vs(i) + \vs(b) = \vs(j) + \vs(a)$
we have either $\vs(i)>\vs(j)$ or $\vs(i)<\vs(j)$
and we use \eqref{Ugradeddr} or \eqref{Ugradeddrop},
respectively; the equality $\vs(i)=\vs(j)$ is impossible as it would
imply $\vs(i) + \vs(a) = \vs(j) + \vs(b)$.

Now suppose that we have a pair of generators
$s_{ia}, s_{jb}$ such that
$\vs(i) + \vs(a) = \vs(j) + \vs (b)$, and that $\vs(i) > \vs(j)$.
Then $\vs(a) < \vs(b)$, and so
\ben
\vs(i) + \vs(b)>\vs(j) + \vs(a).
\een
This means that by applying \eqref{Ugradeddr}, we can
write $s_{ia} s_{jb}$ as a linear combination of
ordered monomials.  Thus, given an arbitrary monomial,
we may rearrange each pair of generators
in turn to write the monomial as a linear combination of
ordered monomials and monomials of lower degree.

Suppose now that
a nontrivial $\CC[z,z^{-1}]$-linear combination
of the monomials \eqref{monomsp}
is zero. By Proposition~\ref{prop:pbwfinext}
we may suppose that at least one coefficient
of the combination does not vanish at $z=1$.
Using the isomorphism \eqref{isomqfoneext} we then
get a nontrivial $\CC$-linear combination
of the corresponding monomials
in the polynomial algebra $\Pc^{\tss\text{\rm ext}}_{2n}$.

Let $\si_{ij}$ denote the image of $s_{ij}$
in $\Pc^{\tss\text{\rm ext}}_{2n}$. Hence
\ben
\si_{ij}=\sum_{a=1}^n
\big(x_{i,2a-1}\ts \bar x_{j,2a}
-x_{i,2a}\ts \bar x_{j,2a-1}\big).
\een
It suffices to verify that
the polynomials $\si_{ij}$ with $i\geqslant j$
and $\si_{i,i+1}$ with odd $i$
are algebraically independent in $\Pc^{\tss\text{\rm ext}}_{2n}$.
Calculate their differentials in terms of
$dx_{ia}$ and
$d\bar x_{ia}$ and specialize the coefficient matrix
by setting
$x_{ij}=\bar x_{ij}=\de_{ij}$.
Then
\ben
d\si_{ij}=\begin{cases}
       d\bar x_{j,i+1}+dx_{i,j-1}
       \qquad&\text{if}\quad i \ \ \text{is odd,}
       \quad j\ \ \text{is even,}\\
       d\bar x_{j,i+1}-dx_{i,j+1}
       \qquad&\text{if}\quad i \ \ \text{is odd,}
       \quad j\ \ \text{is odd,}\\
       -d\bar x_{j,i-1}+dx_{i,j-1}
       \qquad&\text{if}\quad i \ \ \text{is even,}
       \quad j\ \ \text{is even,}\\
       -d\bar x_{j,i-1}-dx_{i,j+1}
       \qquad&\text{if}\quad i \ \ \text{is even,}
       \quad j\ \ \text{is odd,}
                      \end{cases}
\een
so that the differentials are linearly independent.
\epf

Finally, for the use in the next sections we reproduce
the classification theorem for finite-dimensional
irreducible representations of the algebra $\U'_q(\spa_{2n})$.
This theorem was proved in \cite{m:rtq} for the quotient
$\U^{\tss\text{\rm tw}}_q(\spa_{2n})$
of this algebra by the relations \eqref{sqcube},
and it is not difficult
to get the corresponding results for the algebra
$\U'_q(\spa_{2n})$. For the rest of this section
we suppose that $q$ is a nonzero complex number which is not
a root of unity.

A representation $V$ of $\U'_q(\spa_{2n})$ is called
a {\it highest weight representation\/} if $V$ is generated
by a nonzero vector $\xi$ (the {\it highest vector\/}) such that
\begin{alignat}{2}
s_{ij}\ts \xi&=0 \qquad &&\text{for} \quad i=1,3,\dots,2n-1,\quad
j=1,2,\dots,i, \qquad \text{and}
\non\\
s_{2i,2i-1}\ts \xi&=\mu^{}_i\ts \xi,\qquad
&&s_{2i-1,2i}\ts \xi=\mu'_i\ts \xi,
\qquad\qquad\text{for} \quad i=1,2,\dots,n,
\non
\end{alignat}
for some complex numbers $\mu_i$ and $\mu'_i$.
The numbers $\mu'_i$ have to be nonzero due to
the relation \eqref{invrel}.
The tuple $(\mu_1,\dots,\mu_{n};\mu'_1,\dots,\mu'_n)$ is called
the {\it highest weight\/} of $V$.

Due to the PBW theorem for the algebra $\U'_q(\spa_{2n})$
(Proposition~\ref{prop:modPBWforU}), given
any two $n$-tuples of complex numbers
$\mu=(\mu_1,\dots,\mu_{n})$ and $\mu'=(\mu'_1,\dots,\mu'_n)$,
where all $\mu'_i$ are nonzero, there exists an
irreducible highest weight representation $V(\mu;\mu')$
with the highest weight $(\mu;\mu')$.
It is defined as the unique irreducible quotient of the
corresponding Verma module $M(\mu;\mu')$; cf. \cite{m:rtq}.
By definition, $M(\mu;\mu')$ is the quotient
of $\U'_q(\spa_{2n})$ by the left ideal generated by
the elements $s_{ij}$ with $i=1,3,\dots,2n-1$,
$j=1,2,\dots,i$, and by $s_{2i,2i-1}-\mu^{}_i$,
$s_{2i-1,2i}-\mu'_i$ with $i=1,\dots,n$.

\bpr\label{prop:fdco}
Every finite-dimensional irreducible representation
of $\U'_q(\spa_{2n})$ is isomorphic
to a highest weight representation $V(\mu;\mu')$ such
that
\ben
\mu'_i+ q^{2p_i+1}\mu^{}_i=0,\qquad i=1,\dots,n,
\een
for some nonnegative integers $p_i$ satisfying
$p_1\leqslant p_2\leqslant\dots\leqslant p_n$.
\epr

\bpf
A standard argument as in \cite{m:rtq}
shows that
every finite-dimensional irreducible representation
of $\U'_q(\spa_{2n})$ is isomorphic to $V(\mu;\mu')$
for certain $\mu$ and $\mu'$.
In order to find out when an irreducible
highest weight representation $V(\mu;\mu')$
is finite-dimensional, consider
first the case $n=1$. Let $M(\mu_1;\mu'_1)$ be the Verma module
over $\U'_q(\spa_{2})$ with the highest vector $\xi$.
The vectors
$s_{22}^k\tss\xi$ with $k\geqslant 0$ form
a basis of $M(\mu_1;\mu'_1)$.
The central element $s_{22}s_{11}-q^2s_{21}s_{12}$
acts on $M(\mu_1;\mu'_1)$ as multiplication by the scalar
$-q^2\mu_1\mu'_1$.
Hence using the defining relations \eqref{drabs} we derive
by induction on $k$ that
\ben
s^{}_{11}s^k_{22}\tss\xi=
\big(q^{-2k}-1\big)\big(q^2\mu_1\mu'_1+(\mu'_1)^2 q^{3-2k}\big)
\ts s^{k-1}_{22}\tss\xi.
\een
Since $\mu'_1\ne0$, this implies that if $\mu_1=0$ then
$M(\mu_1;\mu'_1)$ is irreducible and so the
representation $V(\mu_1;\mu'_1)$ is infinite-dimensional.

By embedding $\U'_q(\spa_{2})$ into $\U'_q(\spa_{2n})$
as the subalgebra generated by the elements
$s^{}_{2i-1,2i}$, $s^{-1}_{2i-1,2i}$, $s^{}_{2i,2i-1}$,
$s^{}_{2i-1,2i-1}$ and $s^{}_{2i,2i}$
for $i\in\{1,\dots,n\}$,
we can conclude that if the representation $V(\mu;\mu')$
of $\U'_q(\spa_{2n})$ is finite-dimensional,
then all components $\mu_i$ must be nonzero.
Furthermore, each central element
$s_{2i,2i}\tss s_{2i-1,2i-1}-q^2\ts s_{2i,2i-1}\tss
s_{2i-1,2i}$ acts
in $V(\mu;\mu')$ as multiplication
by the nonzero scalar $-q^2\mu_i\mu'_i$.

On the other hand,
the quotient of
$\U'_q(\spa_{2n})$ by the ideal generated by the elements
\beql{centi}
s_{2i,2i}\tss s_{2i-1,2i-1}-q^2\ts s_{2i,2i-1}\tss s_{2i-1,2i}
+q^2\mu_i\mu'_i,\qquad i=1,\dots,n,
\eeq
is isomorphic to the algebra $\U^{\tss\text{\rm tw}}_q(\spa_{2n})$.
Indeed, the mapping
$
s_{ij}\mapsto c_i\tss c_j\tss s_{ij}
$
for nonzero scalars $c_1,\dots,c_{2n}$ such that
\ben
q^2\mu_i\mu'_i=-q^3\ts c^2_{2i-1}c^2_{2i},\qquad i=1,\dots,n,
\een
defines an epimorphism
$\U'_q(\spa_{2n})\to \U^{\tss\text{\rm tw}}_q(\spa_{2n})$
whose kernel is generated by the elements \eqref{centi}.
Thus, $V(\mu;\mu')$ becomes an irreducible
highest weight representation of the algebra
$\U^{\tss\text{\rm tw}}_q(\spa_{2n})$
whose highest weight $\la=(\la_1,\la_3,\dots,\la_{2n-1})$
in the notation of \cite[Sec.~4]{m:rtq} is found by
\ben
\la_{2i-1}=c^{-1}_{2i-1}c^{-1}_{2i}\ts\mu'_i,\qquad i=1,\dots,n.
\een
This implies $\la^2_{2i-1}=-q\mu'_i\ts\mu_i^{-1}$. By
\cite[Theorem~6.3]{m:rtq} we must have
\ben
\la^2_{2i-1}=q^{2m_i},\qquad i=1,\dots,n,
\een
for some positive integers $m_i$ satisfying
$m_1\leqslant m_2\leqslant\dots\leqslant m_n$.
This gives the desired conditions on the highest weight
$(\mu;\mu')$.
\epf

\bre\label{rem:isomsltwo}
If $q^2\ne 1$ then the algebra
$\U^{\tss\text{\rm tw}}_q(\spa_{2})$ is isomorphic to
$\U_q(\sll_2)$. An isomorphism can be given by
\ben
k\mapsto\qin s_{12},\qquad e\mapsto\frac{s_{11}}{q^3-q},\qquad
f\mapsto\frac{s^{-1}_{12}\tss s^{}_{22}}{1-q^2},
\een
where $e,f,k,k^{-1}$ are the standard generators
of $\U_q(\sll_2)$ satisfying the relations
\ben
k\tss e=q^2\tss e\tss k, \qquad k\tss f=q^{-2}\tss f\tss k, \qquad
e\tss f-f\tss e=\frac{k-k^{-1}}{q-\qin}.
\een
This isomorphism can be used to get
a description of finite-dimensional irreducible
representations of $\U^{\tss\text{\rm tw}}_q(\spa_{2})$;
cf. \cite{m:rtq}.
\qed
\ere

\subsection{Quantum affine algebra $\U_q(\wh\gl_N)$}
\label{subsec:quaff}

We start by recalling some well-known facts about the
quantum affine algebra (or quantum loop algebra)
associated with $\gl_N$.
We will keep the notation $q$ for a fixed
nonzero complex number.
Consider the Lie algebra of Laurent polynomials $\gl_N[\la,\la^{-1}]$
in an indeterminate $\la$. We denote it by
$\wh\gl_N$ for brevity. The {\it quantum affine algebra\/}
$\U_q(\wh\gl_N)$ (with the trivial central charge)
has countably many
generators $t_{ij}^{(r)}$ and $\bar t_{ij}^{\ts(r)}$ where
$1\leqslant i,j\leqslant N$ and $r$ runs over nonnegative integers.
They are combined into the matrices
\beql{taff}
T(u)=\sum_{i,j=1}^N t_{ij}(u)\ot E_{ij},\qquad
\overline T(u)=\sum_{i,j=1}^N \bar t_{ij}(u)\ot E_{ij},
\end{equation}
where $t_{ij}(u)$ and $\bar t_{ij}(u)$ are formal series
in $u^{-1}$ and $u$, respectively:
\beql{expa}
t_{ij}(u)=\sum_{r=0}^{\infty}t_{ij}^{(r)}\ts u^{-r},\qquad
\bar t_{ij}(u)=\sum_{r=0}^{\infty}\bar t_{ij}^{\ts(r)}\ts u^{r}.
\end{equation}
The defining relations are
\beql{defrelaff}
\bal
t_{ij}^{(0)}&=\bar t_{ji}^{\ts(0)}=0, \qquad 1 \leqslant i<j\leqslant N,\\
t_{ii}^{(0)}\ts \bar t_{ii}^{\ts(0)}&=\bar t_{ii}^{\ts(0)}
\ts t_{ii}^{(0)}=1,\qquad 1\leqslant i\leqslant N,\\
R(u,v)\ts T_1(u)T_2(v)&=T_2(v)T_1(u)R(u,v),\\
R(u,v)\ts \overline T_1(u)\overline T_2(v)&=
\overline T_2(v)\overline T_1(u)R(u,v),\\
R(u,v)\ts \overline T_1(u)T_2(v)&=T_2(v)\overline T_1(u)R(u,v),
\eal
\end{equation}
where $R(u,v)$ is
the trigonometric $R$-matrix given by
\beql{trRm}
\bal
R(u,v)={}&(u-v)\sum_{i\ne j}E_{ii}\ot E_{jj}+(\qin u-q\tss v)
\sum_{i}E_{ii}\ot E_{ii} \\
{}+ {}&(\qin-q)\tss u\tss\sum_{i> j}E_{ij}\ot
E_{ji}+ (\qin-q)\tss v\tss\sum_{i< j}E_{ij}\ot E_{ji}.
\eal
\end{equation}
Both sides of each of the $R$-matrix relations
are series with coefficients in the algebra
$\U_q(\wh\gl_N)\ot \End\CC^N\ot \End\CC^N$ and the subscripts
of $T(u)$ and $\overline T(u)$ indicate the copies of $\End\CC^N$;
e.g. $T_1(u)=T(u)\ot 1$.
In terms of the
generators these relations
can be written more explicitly as
\beql{dfaff1}
\bal
&(q^{-\delta_{ij}} u - q^{\delta_{ij}}v)\ts
t_{ia}(u)\ts t_{jb}(v)
+ (q^{-1} - q)\ts (u\ts \delta_{i>j} +v\ts \delta_{i<j})\ts
t_{ja}(u)\ts t_{ib}(v) \\
&{}=(q^{-\delta_{ab}}u - q^{\delta_{ab}} v)\ts t_{jb}(v)\ts t_{ia}(u)
+ (q^{-1} - q)\ts (u\ts \delta_{a<b}
+ v\ts \delta_{a>b})\ts t_{ja}(v)\ts t_{ib}(u)
\eal
\eeq
for the relations involving the $t_{ij}^{(r)}$,
\beql{dfaff2}
\bal
&(q^{-\delta_{ij}} u - q^{\delta_{ij}}v)\ts
\bar{t}_{ia}(u)\ts \bar{t}_{jb}(v)
+ (q^{-1} - q)\ts (u\ts \delta_{i>j} +v\ts \delta_{i<j})\ts
\bar{t}_{ja}(u)\ts \bar{t}_{ib}(v) \\
&{}=(q^{-\delta_{ab}}u - q^{\delta_{ab}} v)\ts
\bar{t}_{jb}(v)\ts \bar{t}_{ia}(u)
+ (q^{-1} - q)\ts (u\ts \delta_{a<b}
+ v\ts \delta_{a>b})\ts \bar{t}_{ja}(v)\ts \bar{t}_{ib}(u)
\eal
\eeq
for the relations involving the $\bar t_{ij}^{\ts(r)}$ and
\beql{dfaff3}
\bal
&(q^{-\delta_{ij}} u - q^{\delta_{ij}}v)\ts
\bar{t}_{ia}(u)\ts t_{jb}(v)
+ (q^{-1} - q)\ts (u\ts \delta_{i>j} +v\ts \delta_{i<j})\ts
\bar{t}_{ja}(u)\ts t_{ib}(v) \\
&{}=(q^{-\delta_{ab}}u - q^{\delta_{ab}} v)
\ts t_{jb}(v)\ts \bar{t}_{ia}(u)
+ (q^{-1} - q)\ts (u\ts \delta_{a<b}
+ v\ts \delta_{a>b})\ts t_{ja}(v)\ts \bar{t}_{ib}(u)
\eal
\eeq
for the relations involving both $t_{ij}^{(r)}$ and
$\bar t_{ij}^{\ts(r)}$.

Note that the last relation in \eqref{defrelaff}
can be equivalently written in the form
\beql{altrel}
R(u,v)\ts T_1(u)\ts \overline T_2(v)
=\overline T_2(v)\ts T_1(u)\ts R(u,v).
\eeq
Indeed, we have the identity
\ben
R(u,v)\ts R_{\qin}(u,v)=(q\tss u-\qin v)(\qin u-q\tss v)\ts 1\ot 1,
\een
where $R_{\qin}(u,v)$ is obtained from $R(u,v)$ by replacing
$q$ with $\qin$. Therefore, the last relation
in \eqref{defrelaff}
can be written as
\ben
R_{\qin}(u,v)\ts T_2(v)\ts \overline T_1(u)=\overline T_1(u)\ts
\ts T_2(v)\ts R_{\qin}(u,v).
\een
Now conjugate both sides by the permutation operator
\beql{permu}
P=\sum_{i,j=1}^N E_{ij}\ot E_{ji},
\eeq
then swap $u$ and $v$ to get \eqref{altrel}, as
\ben
R(u,v)=-P\ts R_{\qin}(v,u)\ts P.
\een
In terms of the generators the relation \eqref{altrel} takes the form
\beql{dfaff4}
\bal
&(q^{-\delta_{ij}} u - q^{\delta_{ij}}v)\ts
{t}_{ia}(u)\ts \bar t_{jb}(v)
+ (q^{-1} - q)\ts (u\ts \delta_{i>j} +v\ts \delta_{i<j})\ts
{t}_{ja}(u)\ts\bar t_{ib}(v) \\
&{}=(q^{-\delta_{ab}}u - q^{\delta_{ab}} v)
\ts\bar t_{jb}(v)\ts {t}_{ia}(u)
+ (q^{-1} - q)\ts (u\ts \delta_{a<b}
+ v\ts \delta_{a>b})\ts\bar t_{ja}(v)\ts {t}_{ib}(u).
\eal
\eeq

Let $f(u)$ and $\bar f(u)$ be formal power series in $u^{-1}$
and $u$, respectively,
\ben
\bal
f(u)&=f_0+f_1\ts u^{-1}+f_2\ts  u^{-2}+\dots,\\
\bar f(u)&=
\bar f_0+\bar f_1\ts u+\bar f_2\ts  u^2+\dots,
\eal
\een
such that $f_0\ts \bar f_0=1$.
Then it is immediate from the defining relations that the mapping
\beql{automsln}
T(u)\mapsto f(u)\ts T(u),\qquad
\overline T(u)\mapsto \bar f(u)\ts \overline T(u)
\eeq
defines an automorphism of the algebra $\U_q(\wh\gl_N)$.

We will also use an involutive automorphism of the
algebra $\U_q(\wh\gl_N)$ given by
\beql{autotransp}
T(u)\mapsto \overline T(u^{-1})^t,\qquad
\overline T(u)\mapsto T(u^{-1})^t,
\eeq
where $t$ denotes the matrix transposition.
The first two sets of relations in \eqref{defrelaff}
are obviously preserved by the map \eqref{autotransp}.
In order to verify that the $R$-matrix relations
are preserved as well, apply the transposition
$t_1$ in the first copy of $\End\CC^N$
to each of them, followed by the transposition $t_2$
in the second copy of $\End\CC^N$.
Then conjugate both sides by the permutation
operator \eqref{permu}, replace $u$ and $v$ by
$v^{-1}$ and $u^{-1}$ respectively, and observe that
\ben
u\tss v\ts PR^{t_1t_2}(v^{-1},u^{-1})\ts P=R(u,v).
\een

Another involutive automorphism is defined by
the mapping
\beql{aautoaffsign}
t_{ij}(u)\mapsto \ve_i\ts t_{ij}(u),\qquad
\bar t_{ij}(u)\mapsto \ve_i\ts \bar t_{ij}(u),
\eeq
where each $\ve_i$ equals $1$ or $-1$.

It follows easily from the defining relations \eqref{defrelaff}
that the mapping
\beql{antiautopr}
t_{ij}(u)\mapsto t_{N-j+1,N-i+1}(u),
\qquad \bar t_{ij}(u)\mapsto \bar t_{N-j+1,N-i+1}(u)
\eeq
defines an involutive
anti-automorphism of the algebra $\U_q(\wh\gl_N)$.

Ding and Frenkel \cite{df:it} used the Gauss decompositions
of the matrices $T(u)$ and $\overline T(u)$
to construct an isomorphism  between
the $RTT$-presentation \eqref{defrelaff}
and Drinfeld's `new realization'
of $\U_q(\wh\gl_N)$; see also \cite{fm:ha}.
However, the version of the PBW theorem given in
\cite{bcp:ac} for the new realization of
the quantum affine algebras $\U_q(\wh\agot)$
does not immediately imply
a PBW-type theorem for the $RTT$-presentation.
Our next goal is to prove this theorem, where
$q$ is considered to be an arbitrary fixed
nonzero complex number.

As before, we let $z$ denote an indeterminate.
Introduce the algebra $\U^{\circ}_z(\wh\gl_N)$
over $\CC[z,z^{-1}]$
by the respective generators and relations
given in \eqref{defrelaff}
with $q$ replaced by $z$. Then we have the isomorphism
\beql{isomaff}
\U^{\circ}_z(\wh\gl_N)\ot_{\CC[z,z^{-1}]}\CC\cong \U_q(\wh\gl_N),
\eeq
where the $\CC[z,z^{-1}]$-module $\CC$ is defined
via the evaluation of the Laurent polynomials at $z=q$.
The next proposition takes care of the weak part of
the PBW theorem. We use a particular total order on the generators
of the algebra for which the argument appears to be
the most straightforward. For the purposes of representation
theory a different order is more useful and we will take care
of that one in Corollary~\ref{cor:pbwaffrep} below.

We associate the triple $(i,a,r)$ to
each nonzero generator $t_{ia}^{(r)}$ or $\bar t_{ia}^{\ts(r)}$
of $\U^{\circ}_z(\wh\gl_N)$.
If $(i,a,r)<(j,b,s)$ in the lexicographical order
then we will say that each generator associated with $(i,a,r)$
precedes each generator associated with $(j,b,s)$.
Moreover, we will suppose that $t_{ia}^{(r)}$
precedes $\bar t_{ia}^{\ts(r)}$ for each triple $(i,a,r)$
such that both generators are nonzero.

\bpr\label{prop:span}
The ordered monomials in the generators span the algebra
$\U^{\circ}_z(\wh\gl_N)$ over $\CC[z,\zin]$.
\epr

\bpf
Let $r$ and $s$ be nonnegative integers.
Multiply both sides of the
relation \eqref{dfaff1} with $q$ replaced with $z$
by
\ben
\frac{1}{z^{-\delta_{ij}} u - z^{\tss\delta_{ij}}v}
=\sum_{k=1}^{\infty} z^{(2k-1)\tss\delta_{ij}}u^{-k}\ts v^{k-1}
\een
and equate the coefficients of $u^{-r}v^{-s}$. This provides
an expression for the product $t_{ia}^{(r)}t_{jb}^{(s)}$
with $i>j$ as a $\CC[z,\zin]$-linear
combination of the elements of the form
$t_{jc}^{(k)}t_{id}^{(l)}$. Furthermore,
taking $i=j$ in \eqref{dfaff1} with $q$ replaced with $z$ we obtain
\begin{multline}
\label{partde}
(z^{-ab}u - z^{ab} v)\ts t_{ib}(v)\ts t_{ia}(u)\\
=(z^{-1} u - z\tss v)\ts
t_{ia}(u)\ts t_{ib}(v)- (z^{-1} - z)\ts (u\ts \delta_{a<b}
+ v\ts \delta_{a>b})\ts t_{ia}(v)\ts t_{ib}(u).
\end{multline}
This allows us to express the product $t_{ib}^{(r)}t_{ia}^{(s)}$
with $b>a$ as a $\CC[z,\zin]$-linear
combination of the elements of the form
$t_{ia}^{(k)}t_{ib}^{(l)}$. Taking $a=b$ in \eqref{partde},
we find that the generators $t_{ia}^{(r)}$ and $t_{ia}^{(s)}$
commute for any $r$ and $s$.

Applying similar arguments to the relations
\eqref{dfaff2} and \eqref{dfaff3} with $q$ replaced with $z$
and using induction on the length of monomials we conclude that any
monomial in the generators of $\U^{\circ}_z(\wh\gl_N)$
can be written as a $\CC[z,\zin]$-linear
combination of the ordered monomials.
\epf

\medskip

Recall now that the algebra $\U^{\circ}_z(\wh\gl_N)$
possesses a Hopf algebra structure with the coproduct
\ben
\Delta:\U^{\circ}_z(\wh\gl_N)\to \U^{\circ}_z(\wh\gl_N)
\ot \U^{\circ}_z(\wh\gl_N),
\een
where the tensor product is taken over $\CC[z,\zin]$,
defined by
\beql{copraff}
\Delta\big(t_{ij}(u)\big)=\sum_{k=1}^N t_{ik}(u)\ot t_{kj}(u),
\qquad
\Delta\big(\bar t_{ij}(u)\big)=\sum_{k=1}^N
\bar t_{ik}(u)\ot\bar t_{kj}(u).
\eeq
The quantized enveloping algebra
$\U^{\circ}_z(\gl_N)$ is a Hopf subalgebra of
$\U^{\circ}_z(\wh\gl_N)$
defined by the embedding
\beql{emb}
t_{ij}\mapsto t_{ij}^{(0)},\qquad \bar t_{ij}
\mapsto\bar t_{ij}^{\ts(0)}.
\end{equation}
Moreover, the mapping
\beql{eval}
\pi:T(u)\mapsto T+\overline T\ts u^{-1},\qquad
\overline T(u)\mapsto \overline T+T\ts u
\end{equation}
defines a $\CC[z,\zin]$-algebra homomorphism
$\U^{\circ}_z(\wh\gl_N)\to \U^{\circ}_z(\gl_N)$
called the evaluation homomorphism.

In our proof of
the PBW theorem for the quantum affine algebra
we will follow the approach
used in \cite{bk:pp}
to prove the corresponding theorem for the Yangian for $\gl_N$;
see also \cite{g:gd} for the super-version of the same approach.

We will need a simple lemma which is easily verified
by induction.
Let $x_1,\dots,x_l$
be indeterminates. For $r=0,\dots,l-1$
and $k=1,\dots,l$ consider the elementary symmetric polynomials
in $l-1$ variables, where the variable $x_k$ is skipped:
\ben
e_{rk}=e_r(x_1,\dots,\wh{x}_k,\dots x_l)=\sum x_{i_1}\dots x_{i_r},
\een
summed over indices $i_a\ne k$
with $1\leqslant i_1<\dots<i_r\leqslant l$.

\ble\label{lem:elemsf}
We have
\beql{matresf}
\det\begin{bmatrix} e_{01}&e_{02}&\cdots&e_{0l}\\
                  e_{11}&e_{12}&\cdots&e_{1l}\\
                  \vdots&\vdots&\ddots&\vdots\\
                  e_{l-1,1}&e_{l-1,2}&\cdots&e_{l-1,l}
\end{bmatrix}
=\prod_{1\leqslant i<j\leqslant l}(x_i-x_j).
\eeq
In particular, the determinant
is nonzero under any specialization of variables
$x_i=a_i$, $i=1,\dots,l$, where the $a_i$ are
distinct complex numbers.
\qed
\ele

For each positive integer $l$ introduce the
$\CC[z,\zin]$-algebra homomorphism
\ben
\kappa_l:\U^{\circ}_z(\wh\gl_N)\to \U^{\circ}_z(\gl_N)^{\ot l}
\een
by setting
\beql{kappal}
\kappa_l=\pi^{\ot l}\circ \Delta^{(l-1)},
\eeq
where
\ben
\Delta^{(l-1)}: \U^{\circ}_z(\wh\gl_N)\to
\U^{\circ}_z(\wh\gl_N)^{\ot l}
\een
denotes the coproduct iterated $l-1$ times. The explicit
formulas for the images of the generators of
$\U^{\circ}_z(\wh\gl_N)$ under the homomorphism $\kappa_l$
have the following form:
\beql{kat}
\bal
\kappa_l:t_{ij}^{(r)}&\mapsto
\sum_{p_1<\dots<p_r}\sum_{i_1,\dots,i_l}
t_{ii_1}\ot t_{i_1i_2}\ot\dots\ot \bar t_{i_{p_1-1}i_{p_1}}\ot
\dots\ot \bar t_{i_{p_r-1}i_{p_r}}\ot\dots\ot t_{i_{l-1}j},\\
\kappa_l:\bar t_{ij}^{\ts(r)}&\mapsto
\sum_{p_1<\dots<p_r}\sum_{i_1,\dots,i_l}
\bar t_{ii_1}\ot \bar t_{i_1i_2}\ot\dots\ot t_{i_{p_1-1}i_{p_1}}\ot
\dots\ot t_{i_{p_r-1}i_{p_r}}\ot\dots\ot \bar t_{i_{l-1}j},
\eal
\eeq
where the indices $i_1,\dots,i_l$ in each formula
run over the set $\{1,\dots,N\}$ and the indices
$\{p_1,\dots,p_r\}\subset\{1,\dots,l\}$ indicate the
places taken by the barred generators $\bar t_{kl}$
(resp. unbarred generators $t_{kl}$)
in the first (resp. second) formula.
The images of $t_{ij}^{(r)}$ and $\bar t_{ij}^{\ts(r)}$
under the homomorphism $\kappa_l$ are zero unless $l\geqslant r$.

With the order on the generators of $\U^{\circ}_z(\wh\gl_N)$
introduced before Proposition~\ref{prop:span}
consider the corresponding ordered monomials.
The zero generators
$t_{ij}^{(0)}$ for $i<j$
and $\bar t_{ij}^{\ts(0)}$ for $i>j$ will be excluded.
Moreover, using the relation
$t_{ii}^{(0)}\tss\bar t_{ii}^{\ts(0)}=1$
we will suppose that for each $i=1,\dots,N$
each monomial contains either a nonnegative power of $t_{ii}^{(0)}$
or a positive power of $\bar t_{ii}^{\ts(0)}$.
With these conventions we have the following version of the
PBW theorem.

\bth\label{thm:pbwaff}
The ordered monomials in
the generators $t_{ij}^{(r)}$ and $\bar t_{ij}^{\ts(r)}$
form a basis of the algebra $\U^{\circ}_z(\wh\gl_N)$
over $\CC[z,\zin]$.
\eth

\bpf
Due to Proposition~\ref{prop:span}, we only need to verify
that the ordered monomials are linearly independent.
We will argue by contradiction. Suppose that
a nontrivial linear combination of the ordered monomials
is zero.
Let $m$ be the minimum nonnegative integer
such that for all generators $t_{ij}^{(r)}$ and $\bar t_{ij}^{\ts(r)}$
occurring in the combination we have $0\leqslant r\leqslant m$.
Consider the homomorphism $\kappa_l$ defined in \eqref{kappal}
with $l=2m+1$ and apply it to the linear combination.
We then get the respective nontrivial $\CC[z,\zin]$-linear
combination of elements of the algebra
$\U^{\circ}_z(\gl_N)^{\ot l}$ equal to zero.
By Proposition~\ref{prop:pbwfin}
we may suppose that at least one coefficient
of the combination does not vanish at $z=1$.

On the other hand, due to \eqref{isomqfone}
we have the isomorphism
\ben
\U^{\circ}_z(\gl_N)^{\ot l}\ot_{\CC[z,\zin]}\CC\cong \Pc_N^{\ot l}.
\een
Taking the image of the linear combination under this
isomorphism we get a nontrivial $\CC$-linear combination
of elements of the polynomial algebra $\Pc_N^{\ot l}$
equal to zero.

We will regard $\Pc_N^{\ot l}$ as the algebra of
polynomials in $l$ sets of variables
$\{x_{ij}^{[k]}, \bar x_{ij}^{[k]}\}$, where
the parameter $k\in\{1,\dots,l\}$ indicates the
$k$-th copy of $\Pc_N$ in the tensor product.
Thus, the proof of the theorem is now reduced to verifying
the following claim. Consider
the elements $y_{ij}^{(r)}$ and $\bar y_{ij}^{(r)}$
of the algebra $\Pc_N^{\ot l}$
defined by the relations
\ben
\bal
y_{ij}^{(r)}&=
\sum_{p_1<\dots<p_r}\sum_{i_1,\dots,i_l}
x^{[1]}_{ii_1}\ts x^{[2]}_{i_1i_2}\dots \bar
x^{[p_1]}_{i_{p_1-1}i_{p_1}}
\dots\bar x^{[p_r]}_{i_{p_r-1}i_{p_r}}\dots x^{[l]}_{i_{l-1}j},\\
\bar y_{ij}^{(r)}&=
\sum_{p_1<\dots<p_r}\sum_{i_1,\dots,i_l}
\bar x^{[1]}_{ii_1}\ts \bar x^{[2]}_{i_1i_2}\dots
x^{[p_1]}_{i_{p_1-1}i_{p_1}}
\dots  x^{[p_r]}_{i_{p_r-1}i_{p_r}}\dots \bar x^{[l]}_{i_{l-1}j},
\eal
\een
with the same conditions on the summation indices as in
\eqref{kat} together with the relations
$x^{[s]}_{ij}=\bar x^{[s]}_{ji}=0$ for $i<j$.
We need to verify that modulo
the relations
\ben
\bal
y_{ij}^{(0)}=\bar y_{ji}^{\ts(0)}&=0,
\qquad 1 \leqslant i<j\leqslant N,\\
y_{ii}^{(0)}\tss \bar y_{ii}^{\ts(0)}&=1,\qquad 1\leqslant i\leqslant N,
\eal
\een
the polynomials $y_{ij}^{(r)},\bar y_{ij}^{(r)}$
with $1\leqslant i,j\leqslant N$ and $0\leqslant r\leqslant m$
are algebraically
independent. It will be sufficient to show that
the corresponding
differentials $dy_{ij}^{(r)},d\bar y_{ij}^{(r)}$
are linearly independent. In order to do this,
we calculate the matrix of the map
\ben
\big(dx_{ij}^{[s]},d\bar x_{ij}^{[s]}\big)\to
\big(dy_{ij}^{(r)},d\bar y_{ij}^{(r)})
\een
and show that its determinant is nonzero even when
the variables are specialized to
\beql{speci}
x_{ij}^{[s]}=\de_{ij}c_s,\qquad \bar x_{ij}^{[s]}=\de_{ij}c^{-1}_s,
\eeq
where $c_1,\dots,c_l$ are distinct nonzero complex numbers.

If $i>j$ then under the specialization \eqref{speci} we have
\ben
dy_{ij}^{(r)}=c_1\dots c_l\ts\sum_{s=1}^l c_s^{-1}
e_r(c_1^{-2},\dots,\wh {c_s^{-2}},\dots,c_l^{-2})\ts
dx_{ij}^{[s]}
\een
for $r=0,1,\dots,m$, and
\ben
d\bar y_{ij}^{(r)}=c_1\dots c_l\ts\sum_{s=1}^l c_s^{-1}
e_{l-r}(c_1^{-2},\dots,\wh {c_s^{-2}},\dots,c_l^{-2})\ts
dx_{ij}^{[s]}
\een
for $r=1,\dots,m$. Similarly, for $i<j$ we have
\ben
d\bar y_{ij}^{(r)}=c^{-1}_1\dots c^{-1}_l\ts\sum_{s=1}^l c_s\ts
e_r(c^2_1,\dots,\wh {c^2_s},\dots,c^2_l)\ts
d\bar x_{ij}^{[s]}
\een
for $r=0,1,\dots,m$, and
\ben
dy_{ij}^{(r)}=c^{-1}_1\dots c^{-1}_l\ts\sum_{s=1}^l c_s\ts
e_{l-r}(c^2_1,\dots,\wh {c^2_s},\dots,c^2_l)\ts
d\bar x_{ij}^{[s]}
\een
for $r=1,\dots,m$. Note that since
$x_{ii}^{[s]}\bar x_{ii}^{[s]}=1$, we have
$d\bar x_{ii}^{[s]}=
-\big(x_{ii}^{[s]}\big)^{-2}\ts dx_{ii}^{[s]}$.
Therefore, setting $e_{-1}=0$,
for $i=j$ we obtain
\ben
dy_{ii}^{(r)}=c_1\dots c_l\ts\sum_{s=1}^l c_s^{-1}
\Big(e_r(c_1^{-2},\dots,\wh {c_s^{-2}},\dots,c_l^{-2})
-c^{-2}_s\ts e_{r-1}(c_1^{-2},\dots,\wh {c_s^{-2}},\dots,c_l^{-2})
\Big)
\ts
dx_{ii}^{[s]}
\een
for $r=0,1,\dots,m$, and
\ben
d\bar y_{ii}^{(r)}=c_1\dots c_l\ts\sum_{s=1}^l c_s^{-1}
\Big(e_{l-r}(c_1^{-2},\dots,\wh {c_s^{-2}},\dots,c_l^{-2})
-c^{-2}_s\ts e_{l-r-1}(c_1^{-2},\dots,\wh {c_s^{-2}},\dots,c_l^{-2})
\Big)
\ts
dx_{ii}^{[s]}
\een
for $r=1,\dots,m$.

It follows from Lemma~\ref{lem:elemsf} that
in each of the three cases, the determinant of the
$l\times l$ matrix is nonzero. This proves that the
differentials $dy_{ij}^{(r)}$ and $d\bar y_{ij}^{(r)}$
are linearly independent (excluding $dy_{ij}^{(0)}$ for $i<j$
and $d\bar y_{ij}^{(0)}$ for $i\geqslant j$),
thus completing the proof.
\epf

The following corollary is immediate from the
isomorphism \eqref{isomaff}.

\bco\label{cor:pbwaffspe}
Let $q$ be a nonzero complex number.
With the same order on the generators as
in Theorem~\ref{thm:pbwaff},
the ordered monomials in
the generators $t_{ij}^{(r)}$ and $\bar t_{ij}^{\ts(r)}$
form a basis of the algebra $\U_q(\wh\gl_N)$
over $\CC$.
\qed
\eco

Note that the proof of the linear independence
of the ordered monomials in
$\U^{\circ}_z(\wh\gl_N)$ over $\CC[z,\zin]$
does not rely on the ordering used.
Therefore, Theorem~\ref{thm:pbwaff} holds in the
same form for any other ordering, provided
that the corresponding weak form of the
PBW theorem holds;
cf. Proposition~\ref{prop:span}. We will prove
this weak form for another ordering which is useful for
the description of representations of $\U_q(\wh\gl_N)$.

To each nonzero generator $t_{ia}^{(r)}$ and $\bar t_{ia}^{\ts(r)}$
of $\U^{\circ}_z(\wh\gl_N)$
we now associate the triple of integers $(a-i,i,r)$.
The generators will now be ordered in accordance with
the lexicographical order on the corresponding triples and
we will also suppose that $t_{ia}^{(r)}$
precedes $\bar t_{ia}^{\ts(r)}$ for each triple $(a-i,i,r)$
such that both generators are nonzero.
We have the following version of
the PBW theorem.

\bco\label{cor:pbwaffrep}
Let $q$ be a nonzero complex number.
With the order on the generators defined above,
the ordered monomials in
the generators $t_{ij}^{(r)}$ and $\bar t_{ij}^{\ts(r)}$
form a basis of the algebra $\U_q(\wh\gl_N)$
over $\CC$.
\eco

\bpf
As we pointed out above, the linear independence of
the corresponding monomials in the $\CC[z,\zin]$-algebra
$\U^{\circ}_z(\wh\gl_N)$
will follow by the argument used
in the proof of Theorem~\ref{thm:pbwaff}.
We only need to show that the ordered monomials
span this algebra over $\CC[z,\zin]$.
The corollary will then follow from
the isomorphism \eqref{isomaff}.

Arguing as in the proof of Proposition~\ref{prop:span},
we derive from the
relation \eqref{dfaff1} that
\beql{relorg1}
t_{ia}^{(r)}t_{jb}^{(s)} =
\text{linear
combination of}\quad
t_{jb}^{(k)}t_{ia}^{(l)}\quad\text{and}\quad
t_{ja}^{(m)}t_{ib}^{(p)}
\eeq
for some $k,l,m,p$.
Swapping $i$ with $j$ and $a$ with $b$ in \eqref{dfaff1}
we also obtain
\beql{relorg2}
t_{ia}^{(r)}t_{jb}^{(s)} =
\text{linear
combination of}\quad
t_{jb}^{(k)}t_{ia}^{(l)}\quad\text{and}\quad
t_{ib}^{(m)}t_{ja}^{(p)}.
\eeq
Suppose now that $a-i>b-j$. If $a-j\ne b-i$, then
we use the formula \eqref{relorg1} or \eqref{relorg2}
depending on whether $a-j< b-i$ or $b-i< a-j$
to write $t_{ia}^{(r)}t_{jb}^{(s)}$ as a linear combination
of the ordered products of the generators.
If $a-j=b-i$, then either $j<i$ or $i<j$;
the equality $i=j$ is impossible due to the condition $a-i>b-j$.
Again, the product $t_{ia}^{(r)}t_{jb}^{(s)}$ is then
written as a linear combination
of the ordered products of the generators by
\eqref{relorg1} or \eqref{relorg2}, respectively.

Further, suppose that $a-i=b-j$ and $i>j$. Then $a>b$
and $b-i<a-j$ so that \eqref{relorg2} provides
an expression of $t_{ia}^{(r)}t_{jb}^{(s)}$
as a linear combination
of the ordered products of the generators.

The same arguments relying on \eqref{dfaff2}
instead of \eqref{dfaff1}
prove the corresponding statement for
the products of the generators
$\bar t_{ia}^{\ts(r)}$.

Finally, relation \eqref{dfaff3} implies the following
counterpart of \eqref{relorg1}:
\beql{relmixone}
\bar t_{ia}^{\ts(r)}t_{jb}^{(s)} =
\text{linear
combination of}\quad
t_{jb}^{(k)}\bar t_{ia}^{\ts(l)},\quad
\bar t_{ja}^{\ts(m)}t_{ib}^{(p)}\quad\text{and}\quad
t_{ja}^{(h)}\bar t_{ib}^{\ts(n)}.
\eeq
The corresponding counterpart of \eqref{relorg2}
is obtained from \eqref{dfaff4} and it has the form
\beql{relmixtwo}
\bar t_{ia}^{\ts(r)}t_{jb}^{(s)} =
\text{linear
combination of}\quad
t_{jb}^{(k)}\bar t_{ia}^{\ts(l)},\quad
t_{ib}^{(p)}\bar t_{ja}^{\ts(m)}\quad\text{and}\quad
\bar t_{ib}^{\ts(n)}t_{ja}^{(h)}.
\eeq
The above argument can now be applied to the products
$\bar t_{ia}^{\ts(r)}t_{jb}^{(s)}$ allowing one to write it
as a linear combination
of the ordered products of the generators.

Recalling that each
of the generators $t_{ia}^{(r)}$ and $\bar t_{ia}^{\ts(r)}$
commutes with each of $t_{ia}^{(s)}$ and $\bar t_{ia}^{\ts(s)}$
for all $r$ and $s$, we conclude by
an easy induction
that any monomial in
$\U^{\circ}_z(\wh\gl_N)$
can be written as a
$\CC[z,\zin]$-linear combination
of the ordered products of the generators.
\epf

As with the quantized enveloping algebra,
we need to introduce an extended quantum affine
algebra.
We denote by $\U^{\tss\text{\rm ext}}_q(\wh\gl_N)$
the algebra over $\CC$ with countably many
generators $t_{ij}^{(r)}$ and $\bar t_{ij}^{\ts(r)}$,
$1\leqslant i,j\leqslant N$ and $r\geqslant 0$,
together with $t^{(0)-1}_{ii}$ and $\bar t^{\ts(0)-1}_{ii}$
with $1\leqslant i\leqslant N$,
subject to the defining relations \eqref{defrelaff},
where the second set of relations is replaced with
\ben
t_{ii}^{(0)}\tss \bar t_{ii}^{\tss(0)}
=\bar t_{ii}^{\tss(0)}\tss t_{ii}^{(0)},\qquad
t^{(0)}_{ii}\tss t^{(0)-1}_{ii}
=t^{(0)-1}_{ii}\tss t^{(0)}_{ii}=1,\qquad
\bar t^{\tss(0)}_{ii}\tss \bar t_{ii}^{\tss(0)-1}=
\bar t^{\tss(0)-1}_{ii}\tss \bar t_{ii}^{\tss(0)}=1,
\een
for $i=1,\dots, N$.
We have the natural
epimorphism
\beql{uextepi}
\U^{\tss\text{\rm ext}}_q(\wh\gl_N)\to
\U_q(\wh\gl_N)
\eeq
whose kernel is the ideal
of $\U^{\tss\text{\rm ext}}_q(\wh\gl_N)$ generated
by the central elements $t_{ii}^{(0)}\ts \bar t_{ii}^{\ts(0)}-1$
for $i=1,\dots,N$.
We also define the algebra $\U^{\tss\text{\rm ext}}_z(\wh\gl_N)$
over $\CC[z,z^{-1}]$ with the same generators and
relations, where $q$ should be replaced with $z$.

It is straightforward to conclude that
the PBW theorem for the algebra $\U^{\tss\text{\rm ext}}_q(\wh\gl_N)$
holds in the same form as in Corollaries~\ref{cor:pbwaffspe}
and \ref{cor:pbwaffrep}, except for allowing
the generators $t_{ii}^{(0)}$
and $\bar t_{ii}^{\ts(0)}$ to
occur simultaneously in the monomials and their powers can now
run over the set of all integers.

Observe that given any tuple $(\phi_1,\dots,\phi_N)$ of nonzero
complex numbers, the mapping
\beql{aautoaffsignext}
t_{ij}(u)\mapsto \phi_i\ts t_{ij}(u),\qquad
\bar t_{ij}(u)\mapsto \phi_i\ts \bar t_{ij}(u),
\eeq
defines an automorphism of the algebra
$\U^{\tss\text{\rm ext}}_q(\wh\gl_N)$.

\subsection{Twisted $q$-Yangians
$\Y'_q(\oa_{N})$ and $\Y'_q(\spa_{2n})$}

\label{subsec:twqyang}

The {\it twisted $q$-Yangians\/} $\Y'_q(\oa_{N})$
and $\Y'_q(\spa_{2n})$ associated with
the orthogonal Lie algebra $\oa_N$ and symplectic
Lie algebra $\spa_{2n}$
were introduced in \cite{mrs:cs}.
By definition, $\Y'_q(\oa_{N})$
is the subalgebra of $\U_q(\wh\gl_N)$
generated by the coefficients $s_{ij}^{(r)}, r\geqslant 0,$ of
the series
\beql{scoeff}
s_{ij}(u)=\sum_{r=0}^{\infty} s_{ij}^{(r)}\ts u^{-r},
\qquad 1\leqslant i,j\leqslant N,
\eeq
where
\beql{sumatr}
s_{ij}(u)=\sum_{a=1}^N t_{ia}(u)\ts \bar t_{ja}(u^{-1}).
\end{equation}
In the symplectic case, we define $\Y'_q(\spa_{2n})$
as the subalgebra of $\U^{\tss\text{\rm ext}}_q(\wh\gl_{2n})$
generated by the coefficients
$s_{ij}^{(r)}, r\geqslant 0,$ of
the series
\beql{scoeffs}
s_{ij}(u)=\sum_{r=0}^{\infty} s_{ij}^{(r)}\ts u^{-r},
\qquad 1\leqslant i,j\leqslant 2n,
\end{equation}
where
\beql{sumatrs}
s_{ij}(u)=q\ts\sum_{a=1}^{n} t_{i,2a-1}(u)\ts \bar t_{j,2a}(u^{-1})
-\sum_{a=1}^{n} t_{i,2a}(u)\ts \bar t_{j,2a-1}(u^{-1}),
\end{equation}
and by the elements $s_{i,i+1}^{(0)-1}$ with $i=1,3,\dots,2n-1$.

\bre\label{rem:spextaff}
The twisted $q$-Yangian in the symplectic case was defined in
\cite{mrs:cs}
by the above formulas as a subalgebra of
the quantum affine algebra $\U_q(\wh\gl_{2n})$
without using its extension.
The generators of the corresponding algebra $\Y^{\rm tw}_q(\spa_{2n})$
satisfy some extra relations: for any odd $i$
\beql{sqcubeaff}
s^{(0)}_{i+1,i+1}\ts s^{(0)}_{ii}-q^2\ts s^{(0)}_{i+1,i}\ts
s^{(0)}_{i,i+1}=q^3.
\eeq
Moreover,
the elements $s^{(0)}_{i+1,i+1}\ts s^{(0)}_{ii}-q^2\ts s^{(0)}_{i+1,i}\ts
s^{(0)}_{i,i+1}$ are central in $\Y'_q(\spa_{2n})$
and its quotient by the relations \eqref{sqcubeaff}
is isomorphic to $\Y^{\rm tw}_q(\spa_{2n})$.
\qed
\ere

Both in the orthogonal and symplectic cases, the twisted
$q$-Yangians can be equivalently defined as abstract
algebras with quadratic defining relations.
Namely, consider the matrices $S(u)=T(u)\ts \overline T(u^{-1})^t$
and $S(u)=T(u)\ts G\ts\overline T(u^{-1})^t$,
where the matrix $G$ is defined in \eqref{g}.
Then the matrix elements of $S(u)$ are the formal series
$s_{ij}(u)$ given by \eqref{sumatr} and \eqref{sumatrs},
respectively. The coefficients $s_{ij}^{(r)}$ of these
series then satisfy the relations
\beql{rsrsaffs}
R(u,v)\ts S_1(u)\ts R^{\ts t_1}(u^{-1},v)\ts S_2(v)=S_2(v)
\ts R^{\ts t_1}(u^{-1},v)\ts S_1(u)\ts R(u,v),
\eeq
where the $R$-matrix $R(u,v)$ is defined in \eqref{trRm} and
\beql{rtaff}
\bal
R^{\ts t_1}(u,v)={}&(u-v)\sum_{i\ne j}E_{ii}\ot E_{jj}+(\qin u-q\tss v)
\sum_{i}E_{ii}\ot E_{ii} \\
{}+ {}&(\qin-q)\tss u\tss\sum_{i> j}E_{ji}\ot
E_{ji}+ (\qin-q)\tss v\tss\sum_{i< j}E_{ji}\ot E_{ji}.
\eal
\eeq
In terms of the generating series $s_{ij}(u)$ the relation
\eqref{rsrsaffs} takes the form
\beql{draffsu}
\bal
(q^{-\delta_{ij}}u-q^{\delta_{ij}}v)
\ts\al_{ijab}(u,v)
+ (\qin-q)(u\de_{j<i}+v\de_{i<j})&
\ts\al_{jiab}(u,v)\\
=
(q^{-\delta_{ab}}u-q^{\delta_{ab}}v)
\ts\al_{jiba}(v,u)
+ (\qin-q)(u\de_{a<b}+v\de_{b<a})&
\ts\al_{jiab}(v,u),
\eal
\end{equation}
where
\ben
\al_{ijab}(u,v)=(q^{-\delta_{aj}}-q^{\delta_{aj}}uv)
\ts s_{ia}(u)\ts s_{jb}(v)+
(\qin-q)(\de_{j<a}+uv\de_{a<j})\ts s_{ij}(u)\ts s_{ab}(v).
\een

All
coefficients $\bar s_{ij}^{(r)}$
of the matrix elements $\bar s_{ij}(u)$ of the matrices
\beql{matsbar}
\overline S(u)=\overline T(u)\ts T(u^{-1})^t\Fand
\overline S(u)=\overline T(u)\ts G\ts T(u^{-1})^t
\eeq
belong to the subalgebras $\Y'_q(\oa_{N})\subseteq\U_q(\wh\gl_N)$
and $\Y'_q(\spa_{2n})\subseteq\U^{\tss\text{\rm ext}}_q(\wh\gl_{2n})$,
respectively. Moreover, the relations between
the elements $s_{ij}^{(r)}$ and $\bar s_{ij}^{(r)}$
can be derived from those of the algebras
$\U_q(\wh\gl_N)$ and $\U^{\tss\text{\rm ext}}_q(\wh\gl_{2n})$.
They take the form
\beql{rsrsbar}
\bal
R(u,v)\ts \overline S_1(u)\ts R^{\ts t_1}(u^{-1},v)\ts S_2(v)&=S_2(v)
\ts R^{\ts t_1}(u^{-1},v)\ts \overline S_1(u)\ts R(u,v),\\
R(u,v)\ts S_1(u)\ts R^{\ts t_1}(u^{-1},v)\ts \overline S_2(v)&
=\overline S_2(v)
\ts R^{\ts t_1}(u^{-1},v)\ts S_1(u)\ts R(u,v),\\
R(u,v)\ts \overline S_1(u)\ts R^{\ts t_1}(u^{-1},v)\ts
\overline S_2(v)&=\overline S_2(v)
\ts R^{\ts t_1}(u^{-1},v)\ts \overline S_1(u)\ts R(u,v).
\eal
\eeq

The proof of the equivalence of the two definitions
of the twisted $q$-Yangians is based on analogues
of the PBW theorem whose proofs were outlined in \cite{mrs:cs}.
They use a specialization argument based on the fact
the twisted $q$-Yangians are deformations of universal
enveloping algebras. Here we give a different proof
relying on Theorem~\ref{thm:pbwaff}.

In the orthogonal case we consider the same
total order on the set of generators of $\Y'_q(\oa_N)$
as in \cite{mrs:cs}; we order the generators
$s_{ia}^{(r)}$ in accordance with the lexicographical order
on the corresponding triples $(i,a,r)$.

In the symplectic case use the function
$\vs: \{ 1, 2, \ldots, 2n\} \to \{\pm 1, \pm 3,\ldots, \pm(2n-1) \}$
defined in \eqref{nu} and order the generators
$s_{ia}^{(r)}$ of $\Y'_q(\spa_{2n})$
in accordance with the lexicographical order
on the corresponding triples $(\vs(i)+\vs(a),\vs(i),r)$.
Since for any odd $i$
the generators $s_{i,i+1}^{(0)}$ and $s_{i,i+1}^{(0)-1}$
commute, it is unambiguous to associate each of them
to the same triple $(0,i,0)$.

By the definition \eqref{sumatr}
we have
\beql{ortexc}
s_{ij}^{(0)}=0\quad\text{for}\quad i<j\qquad\text{and}\qquad
s_{ii}^{(0)}=1 \quad\text{for all}\quad i
\eeq
in the orthogonal case. Similarly, by \eqref{sumatrs}
in the symplectic case we have
\beql{symexc}
s_{ij}^{(0)}=0\quad\text{for}\quad i<j \quad\text{unless}\ \
j=i+1\quad\text{with}\quad i\ \  \text{odd}.
\eeq
Consequently,
the generators \eqref{ortexc} and \eqref{symexc}
will not occur in the ordered monomials.

\bpr\label{prop:pbwtw}
Let $q$ be a nonzero complex number.
With the orders on the generators chosen as above,
the ordered monomials in the generators
form a basis of the respective algebra $\Y'_q(\oa_{N})$
and $\Y'_q(\spa_{2n})$.
\epr

\bpf
The weak form of the PBW theorem was proved in
\cite[Lemma~3.2]{mrs:cs} for the order
used in the orthogonal case. The proof for the order we chose
in the symplectic case is obtained by obvious
modifications of the same arguments; cf. the proof of
Proposition~\ref{prop:modPBWforU}.
Thus, in both cases
the ordered monomials span the algebra
$\Y'_q(\oa_{N})$ or $\Y'_q(\spa_{2n})$,
respectively.

Furthermore, by Theorem~\ref{thm:pbwaff}, we have the isomorphism
\beql{isomatone}
\U^{\circ}_z(\wh\gl_N)\ot_{\CC[z,z^{-1}]}\CC\cong \wh\Pc_N,
\eeq
where the $\CC[z,z^{-1}]$-module $\CC$ is defined
via the evaluation of the Laurent polynomials at $z=1$
and $\wh\Pc_N$ is the algebra of polynomials
in the variables $x_{ij}^{(r)},\bar x_{ij}^{(r)}$ with
$1\leqslant i,j\leqslant N$ and $r\geqslant 0$
subject to the relations
$x_{ij}^{(0)}=\bar x_{ji}^{(0)}=0$ for $i<j$ and
$x_{ii}^{(0)}\bar x_{ii}^{(0)}=1$ for all $i$.

Define the algebra $\Y^{\circ}_z(\oa_{N})$ over $\CC[z,z^{-1}]$ as
the $\CC[z,z^{-1}]$-subalgebra of $\U^{\circ}_z(\wh\gl_N)$
generated by the elements $s_{ij}^{(r)}$ defined
in \eqref{scoeff} and \eqref{sumatr}.
Suppose that
a nontrivial $\CC[z,z^{-1}]$-linear combination
of the ordered monomials in the generators of $\Y^{\circ}_z(\oa_{N})$
is zero. By Theorem~\ref{thm:pbwaff}
we may suppose that at least one coefficient
of the combination does not vanish at $z=1$.
Using the isomorphism \eqref{isomatone} we then
get a nontrivial $\CC$-linear combination
of the ordered monomials in the images of the generators
in the polynomial algebra $\wh\Pc_N$. We will come
to a contradiction if we show that the images
of the generators of $\Y^{\circ}_z(\oa_{N})$
in $\wh\Pc_N$ are algebraically independent.

Let $\si_{ij}^{(r)}$ denote the image of $s_{ij}^{(r)}$
in $\wh\Pc_N$. Then
\ben
\si_{ij}^{(r)}=\sum_{a=1}^N
\sum_{k+l=r} x_{ia}^{(k)}\ts \bar x_{ja}^{(l)}.
\een
It suffices to verify that the differentials
$d\si_{ij}^{(r)}$ are linearly independent.
Calculate the differentials in terms of
$dx_{ia}^{(k)}$ and
$d\bar x_{ia}^{(k)}$ and specialize the coefficient matrix
by setting
\beql{sce}
x_{ij}^{(k)}=\bar x_{ij}^{(k)}=\de_{ij}\ts\de_{k0}.
\eeq
Then
\ben
d\si_{ij}^{(r)}=d\bar x_{ji}^{(r)}+dx_{ij}^{(r)}
\een
which implies that the differentials $d\si_{ij}^{(r)}$
are linearly independent even under the specialization \eqref{sce}.
This completes the proof in the orthogonal case.

In the symplectic case define
the algebra $\Y^{\circ}_z(\spa_{2n})$ over $\CC[z,z^{-1}]$ as
the $\CC[z,z^{-1}]$-subalgebra of
$\U^{\tss\text{\rm ext}}_z(\wh\gl_{2n})$
generated by the elements $s_{ij}^{(r)}$ defined
by \eqref{scoeffs} and \eqref{sumatrs}
with $q$ replaced
by $z$ and with the same additional generators
$s_{i,i+1}^{(0)-1}$. Suppose that
a nontrivial $\CC[z,z^{-1}]$-linear combination
of the ordered monomials in the generators of $\Y^{\circ}_z(\spa_{2n})$
is zero. We may ignore the generators $s_{i,i+1}^{(0)-1}$
because they
can be excluded from the linear
combination by multiplying it by appropriate powers of the elements
$s_{i,i+1}^{(0)}$ and using the following consequence
of \eqref{draffsu} (with $q$ replaced by $z$):
\beql{siiponeskl}
s^{(0)}_{i,i+1}\ts s^{}_{kl}(u)=
z^{\de_{ik}+\de_{il}-\de_{i+1,k}-\de_{i+1,l}}\ts
s^{}_{kl}(u)\ts s^{(0)}_{i,i+1},\qquad
i=1,3,\dots,2n-1.
\eeq
By the PBW theorem for the algebra
$\U^{\tss\text{\rm ext}}_z(\wh\gl_{2n})$ (see Sec.~\ref{subsec:quaff})
we have the isomorphism
\beql{isomatonesp}
\U^{\tss\text{\rm ext}}_z(\wh\gl_{2n})
\ot_{\CC[z,z^{-1}]}\CC\cong \wh\Pc^{\tss\text{\rm ext}}_{2n},
\eeq
where the $\CC[z,z^{-1}]$-module $\CC$ is defined
via the evaluation of the Laurent polynomials at $z=1$
and $\wh\Pc^{\tss\text{\rm ext}}_{2n}$ is the algebra of polynomials
in the variables $x_{ij}^{(r)},\bar x_{ij}^{(r)}$ with
$1\leqslant i,j\leqslant N$ and $r\geqslant 0$
together with $x_{ii}^{(0)-1}$ and $\bar x_{ii}^{(0)-1}$
subject to the relations
$x_{ij}^{(0)}=\bar x_{ji}^{(0)}=0$ for $i<j$ and
\ben
x_{ii}^{(0)} x_{ii}^{(0)-1}=1,\qquad
\bar x_{ii}^{(0)} \bar x_{ii}^{(0)-1}=1
\een
for all $i$. As in the orthogonal case, it suffices to verify
that the images $\si_{ij}^{(r)}$
of the generators $s_{ij}^{(r)}$
in $\wh\Pc^{\tss\text{\rm ext}}_{2n}$ are algebraically independent.
Calculating the differentials and specializing
the variables as in \eqref{sce}, we get
\ben
d\si_{ij}^{(r)}=\begin{cases}
       d\bar x_{j,i+1}^{(r)}+dx_{i,j-1}^{(r)}
       \qquad&\text{if}\quad i \ \ \text{is odd,}
       \quad j\ \ \text{is even,}\\[5pt]
       d\bar x_{j,i+1}^{(r)}-dx_{i,j+1}^{(r)}
       \qquad&\text{if}\quad i \ \ \text{is odd,}
       \quad j\ \ \text{is odd,}\\[5pt]
       d\bar x_{j,i-1}^{(r)}+dx_{i,j-1}^{(r)}
       \qquad&\text{if}\quad i \ \ \text{is even,}
       \quad j\ \ \text{is even,}\\[5pt]
       d\bar x_{j,i-1}^{(r)}-dx_{i,j+1}^{(r)}
       \qquad&\text{if}\quad i \ \ \text{is even,}
       \quad j\ \ \text{is odd,}
                      \end{cases}
\een
which shows that the differentials are
linearly independent in this case
as well.
\epf

\section{Representations of the quantum affine algebra}
\label{sec:qaff}
\setcounter{equation}{0}

As in the Lie algebra representation theory, the
representations of the quantum affine algebra
associated with $\sll_2$ plays a key role
in the description of the representations of the
quantum affine algebras $\U_q(\wh\agot)$; see
\cite{cp:qaa}, \cite{cp:gq}. Finite-dimensional
irreducible representations of $\U_q(\wh\sll_2)$
were classified in \cite{cp:qaa}. In our proofs
below the case of the twisted $q$-Yangian $\Y'_q(\spa_2)$
will be similarly important for the general
classification theorem. In order to make our arguments
clearer, we first reproduce a proof of the classification
theorem for the representations of $\U_q(\wh\gl_2)$
following an approach
used for the Yangian representations and which goes back
to pioneering work of Tarasov~\cite{t:sq, t:im}.
This approach is alternative
to \cite{cp:qaa} and it also allows one to obtain a
description of the finite-dimensional irreducible
representations of $\U_q(\wh\gl_2)$ as tensor products
of the evaluation modules
as an immediate corollary. The corresponding
arguments were outlined in \cite[Sec.~3.5]{m:yc}.

Suppose that the complex number $q$
is nonzero and not a root of unity.
A representation $L$ of $\U_q(\wh\gl_N)$ is called
a {\it highest weight representation\/}
if $L$ is generated
by a nonzero vector $\ze$ (the {\it highest vector\/})
such that
\begin{alignat}{3}
t_{ij}(u)\ts\ze&=0,\qquad &&\bar t_{ij}(u)\ts\ze=0 \qquad
&&\text{for} \quad
1\leqslant i<j\leqslant N,
\non\\
t_{ii}(u)\ts\ze&=\nu_i(u)\ts\ze,\qquad
&&\bar t_{ii}(u)\ts\ze=\bar \nu_i(u)\ts\ze
\qquad &&\text{for} \quad 1\leqslant i\leqslant N,
\non
\end{alignat}
where $\nu(u)=(\nu_1(u),\dots,\nu_N(u))$ and
$\bar\nu(u)=(\bar\nu_1(u),\dots,\bar\nu_N(u))$ are certain $N$-tuples
of formal power series in $u^{-1}$ and $u$, respectively:
\beql{nui}
\nu_i(u)=\sum_{r=0}^{\infty} \nu_i^{(r)}\tss u^{-r},
\qquad
\bar\nu_i(u)=\sum_{r=0}^{\infty} \bar\nu_i^{(r)}\tss u^{r}.
\eeq
We have $\nu_i^{(0)}\bar\nu_i^{(0)}=1$ for each $i$
due to the second set of relations in \eqref{defrelaff}.

Note that this definition corresponds to
{\it pseudo-highest weight representations\/}
of the quantum loop algebras
in the terminology of \cite[Def.~12.2.4]{cp:gq}.

A standard argument shows that any
finite-dimensional irreducible representation
of $\U_q(\wh\gl_N)$ is a highest weight representation;
cf. \cite[Prop.~12.2.3]{cp:gq}. Furthermore,
Corollary~\ref{cor:pbwaffrep} implies that
given any formal series of the form \eqref{nui}
with $\nu_i^{(0)}\bar\nu_i^{(0)}=1$ for all $i$,
there exists a nontrivial {\it Verma module\/}
$M(\nu(u);\bar\nu(u))$ which is defined as the quotient
of $\U_q(\wh\gl_N)$ by the left ideal generated by
all coefficients of the series $t_{ij}(u)$, $\bar t_{ij}(u)$
for $i<j$ and $t_{ii}(u)-\nu_i(u)$, $\bar t_{ii}(u)-\bar\nu_i(u)$
for all $i$. Moreover, $M(\nu(u);\bar\nu(u))$ has a unique
irreducible quotient $L(\nu(u);\bar\nu(u))$.
Therefore, in order to describe all finite-dimensional
irreducible representations of the algebra $\U_q(\wh\gl_N)$,
we need to determine for which highest weights
$(\nu(u);\bar\nu(u))$ the representation $L(\nu(u);\bar\nu(u))$
is finite-dimensional. By considering `simple root embeddings'
$\U_q(\wh\gl_2)\hookrightarrow\U_q(\wh\gl_N)$, the problem
is largely reduced to the particular case $N=2$.

\subsection{Representations of $\U_q(\wh\gl_2)$}
\label{subsec:qafftwo}

Consider an arbitrary irreducible highest weight representation
$L\big(\nu_1(u),\nu_2(u);\bar\nu_1(u),\bar\nu_2(u)\big)$
of the algebra $\U_q(\wh\gl_2)$.

\bpr\label{prop:polyn}
Suppose that
$\dim L\big(\nu_1(u),\nu_2(u);\bar\nu_1(u),\bar\nu_2(u)\big)<\infty$.
Then there exist polynomials $Q(u)$ and $R(u)$ in $u$
of the same degree such that
the product of the constant term and the leading coefficient
of each polynomial is equal to $1$, and
\beql{fracrel}
\frac{\nu_1(u)}{\nu_2(u)}=
\frac{Q(u)}{R(u)}
=\frac{\bar\nu_1(u)}{\bar\nu_2(u)},
\eeq
where the first equality
is understood in the sense
that the ratio of polynomials has to be
expanded as a power series in $u^{-1}$,
while for the second equality the same ratio has to be expanded
as a power series in $u$.
\epr

\bpf
By twisting the representation
$L\big(\nu_1(u),\nu_2(u);\bar\nu_1(u),\bar\nu_2(u)\big)$
with an appropriate
automorphism of the form \eqref{automsln}, we may assume
without loss of generality that $\nu_2(u)=\bar\nu_2(u)=1$.
Consider the vector subspace $L$ of
$L\big(\nu_1(u),1;\bar\nu_1(u),1\big)$
spanned by all vectors $t_{21}^{(i)}\ze$, $i\geqslant 0$
and $\bar t_{21}^{\ts(j)}\ze$, $j\geqslant 1$.
Since $\dim L<\infty$, the space $L$ is spanned
by the vectors $t_{21}^{(i)}\ze$ and $\bar t_{21}^{\ts(j)}\ze$,
where $i$ and $j$ run over some finite sets of values.
This implies that for sufficiently large $n$ and $m$
any vector $c_n\ts t_{21}^{(n)}\ze+d_m\bar t_{21}^{\ts(m)}\ze$
is a linear combination of the spanning vectors
$t_{21}^{(i)}\ze$ and $\bar t_{21}^{\ts(j)}\ze$.
Hence, there exist integers $n\geqslant 0$, $m\geqslant 1$
and complex numbers $c_i,d_j$ such that
\ben
\sum_{i=0}^n c_i\ts t_{21}^{(i)}\ze
+\sum_{j=1}^m d_j\ts \bar t_{21}^{\ts(j)}\ze=0,
\een
and $c_n,d_m\ne 0$. Denote the linear combination which occurs
on the left-hand side by $\xi$. Then $t_{12}^{(r)}\ts\xi=0$
for all $r\geqslant 1$. On the other hand, by the defining relations
\eqref{dfaff1} we have
\ben
(u - v)\ts
\big(t_{12}(u)\ts t_{21}(v)-t_{21}(v)\ts t_{12}(u)\big)=
(q - q^{-1})\ts v\ts \big(t_{22}(u)\ts t_{11}(v)-
t_{22}(v)\ts t_{11}(u)\big).
\een
Now multiply both sides by
\ben
\frac{1}{u-v}=\sum_{k=1}^{\infty}u^{-k}v^{k-1}
\een
and take the coefficients of $u^{-r}v^{-s}$
on both sides to get
\ben
t_{12}^{(r)}\ts t_{21}^{(s)}-t_{21}^{(s)}\ts t_{12}^{(r)}=
(q-\qin)\ts \sum_{p=1}^r
\big(t_{22}^{(r-p)}\ts t_{11}^{(s+p)}-
t_{22}^{(s+p)}\ts t_{11}^{(r-p)}\big).
\een
Similarly, using \eqref{dfaff4} we get
\ben
t_{12}^{(r)}\ts \bar t_{21}^{\ts(s)}-
\bar t_{21}^{\ts(s)}\ts t_{12}^{(r)}=
(q-\qin)\ts \sum_{p=1}^{\min\{r,s\}}
\big(t_{22}^{(r-p)}\ts \bar t_{11}^{\ts(s-p)}-
\bar t_{22}^{\ts(s-p)}\ts t_{11}^{(r-p)}\big).
\een
Since $t_{22}^{(r)}\ts\ze=\bar t_{22}^{\ts(r)}\ts\ze=0$
for $r\geqslant 1$, we find that
\ben
t_{12}^{(r)}\ts t_{21}^{(s)}\ts\ze=(q-\qin)\ts\nu_1^{(r+s)}\ts\ze,
\qquad
t_{12}^{(r)}\ts \bar t_{21}^{\ts(s)}\ts\ze=
(q-\qin)\ts(\bar\nu_1^{(s-r)}-\nu_1^{(r-s)})\ts\ze,
\een
where $s\geqslant 1$ in the second relation and
we assume
$\nu_1^{(s)}=\bar\nu_1^{(s)}=0$ for $s<0$.
Hence, taking the coefficient of $\ze$ in $t_{12}^{(r)}\ts\xi=0$
we get
\ben
\sum_{i=0}^n c_i\tss\nu_1^{(r+i)}+
\sum_{j=1}^m d_j\big(\bar\nu_1^{(j-r)}-\nu_1^{(r-j)}\big)=0
\een
for all $r\geqslant 1$.
This is equivalent to the relation
\ben
\bal
\nu_1(u)\Bigg(\sum_{i=0}^n c_i u^i-\sum_{j=1}^m d_j u^{-j}\Bigg)
{}&=\sum_{i=0}^n c_i u^i
\ts\big(\nu_1^{(0)}+\dots+\nu_1^{(i)}u^{-i}\big)\\
{}&-\sum_{j=1}^m d_j u^{-j}
\ts\big(\bar\nu_1^{(0)}+\dots+\bar\nu_1^{(j-1)}u^{j-1}\big).
\eal
\een
Now use the relations $\bar t_{12}^{\ts(r)}\ts\xi=0$,
$r\geqslant 0$. The defining relations \eqref{dfaff2}
give
\ben
(u - v)\ts
\big(\bar t_{12}(u)\ts \bar t_{21}(v)-
\bar t_{21}(v)\ts \bar t_{12}(u)\big)=
(q - q^{-1})\ts v\ts \big(\bar t_{22}(u)\ts \bar t_{11}(v)-
\bar t_{22}(v)\ts \bar t_{11}(u)\big).
\een
Divide both sides by $u-v$ and use
the expansion
\ben
\frac{v}{u-v}=-\sum_{k=0}^{\infty}u^{k}v^{-k}.
\een
Comparing the coefficients of $u^{r}v^{s}$
on both sides we get
\ben
\bar t_{12}^{\ts(r)}\ts \bar t_{21}^{\ts(s)}-
\bar t_{21}^{\ts(s)}\ts \bar t_{12}^{\ts(r)}=
(q-\qin)\ts \sum_{p=0}^r
\big(\bar t_{22}^{\ts(s+p)}\ts \bar t_{11}^{\ts(r-p)}-
\bar t_{22}^{\ts(r-p)}\ts \bar t_{11}^{\ts(s+p)}\big).
\een
Similarly, using \eqref{dfaff3} we get
\ben
\bar t_{12}^{\ts(r)}\ts t_{21}^{(s)}-
t_{21}^{(s)}\ts \bar t_{12}^{\ts(r)}=
(q-\qin)\ts \sum_{p=0}^{\min\{r,s\}}
\big(t_{22}^{(s-p)}\ts \bar t_{11}^{\ts(r-p)}-
\bar t_{22}^{\ts(r-p)}\ts t_{11}^{(s-p)}\big).
\een
Hence,
\ben
\bar t_{12}^{\ts(r)}\ts \bar t_{21}^{\ts(s)}\ts\ze=
(\qin-q)\ts\bar\nu_1^{\ts(r+s)}\ts\ze,
\qquad
\bar t_{12}^{\ts(r)}\ts t_{21}^{(s)}\ts\ze=
(q-\qin)\ts(\bar\nu_1^{(r-s)}-\nu_1^{(s-r)})\ts\ze,
\een
where $s\geqslant 1$ in the first relation.
Taking the coefficient of $\ze$ in $\bar t_{12}^{\ts(r)}\ts\xi=0$
we get
\ben
\sum_{i=0}^n c_i\tss\big(\bar\nu_1^{(r-i)}-\nu_1^{(i-r)}\big)-
\sum_{j=1}^m d_j\ts\bar\nu_1^{(r+j)}=0
\een
for all $r\geqslant 0$. This is equivalent to the relation
\ben
\bal
\bar\nu_1(u)\Bigg(\sum_{i=0}^n c_i u^i-\sum_{j=1}^m d_j u^{-j}\Bigg)
{}&=\sum_{i=0}^n c_i u^i
\ts\big(\nu_1^{(0)}+\dots+\nu_1^{(i)}u^{-i}\big)\\
{}&-\sum_{j=1}^m d_j u^{-j}
\ts\big(\bar\nu_1^{(0)}+\dots+\bar\nu_1^{(j-1)}u^{j-1}\big).
\eal
\een
Thus, both series $\nu_1(u)$ and $\bar\nu_1(u)$
are expansions of the same rational function in $u$,
\ben
\nu_1(u)=
\frac{Q(u)}{R(u)}
=\bar\nu_1(u),
\een
where the polynomials $Q(u)$ and $R(u)$ have the required
properties.
\epf

Write decompositions
\ben
\bal
Q(u)&=(\al^{}_1u+\al_1^{-1})\dots (\al^{}_k u+\al_k^{-1}),\\
R(u)&=(\be^{}_1u+\be_1^{-1})\dots (\be^{}_k u+\be_k^{-1}),
\eal
\een
where $\al_i$ and $\be_i$ are nonzero complex numbers.
By twisting the finite-dimensional
representation $L\big(\nu_1(u),\nu_2(u);\bar\nu_1(u),\bar\nu_2(u)\big)$
by an appropriate automorphism
of the form \eqref{automsln},
we can get another finite-dimensional
representation such that the components of the highest weight
have the form:
\beql{hwalbe}
\bal
\nu_1(u)&=
(\al^{}_1+\al^{-1}_1u^{-1})\dots (\al^{}_k+\al^{-1}_ku^{-1}),
\\
\nu_2(u)&=
(\be^{}_1+\be^{-1}_1u^{-1})\dots (\be^{}_k+\be^{-1}_ku^{-1}),
\\
\bar\nu_1(u)&=
(\al^{}_1u+\al^{-1}_1)\dots (\al^{}_ku+\al^{-1}_k),
\\
\bar\nu_2(u)&=
(\be^{}_1u+\be^{-1}_1)\dots (\be^{}_ku+\be^{-1}_k).
\eal
\eeq
For any pair of nonzero
complex numbers $\al$ and $\be$ consider the corresponding
irreducible highest weight representation $L(\al,\be)$
of $\U_q(\gl_2)$.
That is, $L(\al,\be)$ is generated by a nonzero
vector $\ze$ such that
\ben
\bar t_{12}\ts\ze=0,\qquad
t_{11}\ts\ze=\al\tss\ze, \qquad
t_{22}\ts\ze=\be\tss\ze.
\een
The representation $L(\al,\be)$ is finite-dimensional
if and only if $\al/\be=\pm\tss q^m$ for some nonnegative
integer $m$.
We make $L(\al,\be)$ into a module over
the quantum affine algebra
$\U_q(\wh\gl_2)$ via the evaluation
homomorphism $\U_q(\wh\gl_2)\to\U_q(\gl_2)$
given by the formulas \eqref{eval}.
This evaluation module is a highest weight representation
of $\U_q(\wh\gl_2)$ with the highest weight
\ben
\big(\al+\al^{-1}u^{-1},\be+\be^{-1}u^{-1};\
\al\ts u+\al^{-1},\be\ts u+\be^{-1}\big).
\een

The comultiplication map \eqref{copraff} allows us
to regard the tensor product
\beql{tenprqaff}
L(\al_1,\be_1)\ot L(\al_2,\be_2)\ot\dots\ot L(\al_k,\be_k)
\eeq
as a representation of $\U_q(\wh\gl_2)$. Moreover,
it follows easily from \eqref{copraff} that
the cyclic span $\U_q(\wh\gl_2)(\ze_1\ot\dots\ot\ze_k)$
is a highest weight representation of $\U_q(\wh\gl_2)$
with the highest weight given by the formulas
\eqref{hwalbe}; here $\ze_i$ denotes the highest vector
of $L(\al_i,\be_i)$. Our next goal is to show that under
some additional conditions on the parameters $\al_i$ and $\be_i$
the tensor product module \eqref{tenprqaff}
is irreducible and, hence, isomorphic to
$L(\nu_1(u),\nu_2(u);\bar\nu_1(u),\bar\nu_2(u))$.
Namely, we will suppose that for every $i=1,\dots, k-1$
the following condition holds:
if the multiset $\{\al_r/\be_s\ |\ i\leqslant r,s \leqslant k\}$
contains numbers of the form $\pm\tss q^m$ with
nonnegative integers $m$,
then $\al_i/\be_i=\pm\tss q^{m_0}$ and $m_0$
is minimal amongst these nonnegative integers.
The following proposition goes back to \cite{t:im}.

\bpr\label{prop:tpr}
If the above condition on the parameters
$\al_i$ and $\be_i$ holds, then the
representation \eqref{tenprqaff} of
$\U_q(\wh\gl_2)$ is irreducible.
\epr

\bpf
We follow the corresponding argument used in the Yangian case;
see e.g. \cite[Prop.~3.3.2]{m:yc}.
Denote the representation \eqref{tenprqaff} by $L$.
We start by proving the following claim: any vector $\xi\in L$
satisfying $t_{12}(u)\tss\xi=0$ is proportional to
$\ze=\ze_1\ot\dots\ot\ze_k$.
We use the induction on $k$. The claim
is obvious for $k=1$ so suppose that $k\geqslant 2$.
Write any such vector $\xi$, which is assumed
to be nonzero, in the form
\ben
\xi=\sum_{r=0}^p(t_{21})^r \ze_1\ot \xi_r,
\qquad
\text{where}
\quad
\xi_r\in L(\al_2,\be_2)\ot\dots\ot L(\al_k,\be_k)
\een
and $p$ is some nonnegative integer. Moreover,
if $\al_1/\be_1=\pm\tss q^m$ for some nonnegative
integer $m$, then we will assume that
$p\leqslant m$.
We will also assume that $\xi_p\ne 0$. Using the coproduct
formulas \eqref{copraff}, we get
\beql{srpte}
\sum_{r=0}^p\Big(t_{11}(u)(t_{21})^r \ze_1\ot t_{12}(u)\tss \xi_r
+t_{12}(u)(t_{21})^r \ze_1\ot t_{22}(u)\tss \xi_r\Big)=0.
\eeq
By \eqref{eval}, we obtain
\beql{toneone}
t_{11}(u)(t_{21})^r \ze_1=(t_{11}+\bar t_{11}u^{-1})(t_{21})^r \ze_1
=\big(q^{-r}\al^{}_1+q^{r}\al^{-1}_1\tss u^{-1}\big)(t_{21})^r \ze_1,
\eeq
and
\beql{tonetwo}
t_{12}(u)(t_{21})^r \ze_1=u^{-1}\bar t_{12}\tss(t_{21})^r \ze_1=
u^{-1}\ts (q^r-q^{-r})\big(q^{r-1}\be_1/\al_1-q^{-r+1}\al_1/\be_1\big)
(t_{21})^{r-1} \ze_1.
\eeq
Taking the coefficient of $(t_{21})^p \ze_1$ in
\eqref{srpte} we get
$t_{12}(u)\tss\xi_p=0$. By the induction hypothesis,
applied to the representation
$
L(\al_2,\be_2)\ot\dots\ot L(\al_k,\be_k),
$
the vector
$\xi_p$ must be proportional to $\ze_2\ot\dots\ot\ze_k\tss$.
Therefore, using again \eqref{copraff} and \eqref{eval},
we obtain
\beql{tuetap}
t_{22}(u)\tss\xi_p=(\be^{}_2+\be^{-1}_2u^{-1})\dots
(\be^{}_k+\be^{-1}_k u^{-1}) \ts\xi_p.
\eeq
The proof of the claim will be completed if show that $p=0$.
Suppose on the contrary that $p\geqslant 1$. Then
taking the coefficient of
$(t_{21})^{p-1} \ze_1$ in \eqref{srpte} we derive
\begin{multline}
\big(q^{-p+1}\al^{}_1+q^{p-1}\al^{-1}_1 u^{-1}\big)
\ts t_{12}(u)\tss\xi_{p-1}\\
{}+u^{-1}\ts
(q^p-q^{-p})\big(q^{p-1}\be_1/\al_1-q^{-p+1}\al_1/\be_1\big)
\ts t_{22}(u)\ts\xi_p=0.
\non
\end{multline}
Note that
$t_{12}(u)\tss\xi_{p-1}$ is a polynomial
in $u^{-1}$.
Taking $u=-q^{2p-2}\al^{-2}_1$ and using \eqref{tuetap}
we obtain the relation
\ben
\big(q^{2p-2}-(\al_1/\be_1)^2\big)\big(q^{2p-2}-(\al_1/\be_2)^2\big)
\dots \big(q^{2p-2}-(\al_1/\be_k)^2\big)=0,
\een
where we also used the assumption that
$q$ is not a root of unity.
However, this is impossible due to the conditions on the parameters
$\al_i$ and $\be_i$. Thus, $p$ must be zero and the claim follows.

Now suppose that $M$ is a nonzero submodule of $L$.
Then $M$ must contain a nonzero vector $\xi$ such that
$t_{12}(u)\tss\xi=0$. This can be seen by
considering $\U_q(\gl_2)$-weights of $M$. If $\eta\in M$
is a vector of weight $(\mu_1,\mu_2)$, i.e.,
$t_{11}\eta=\mu_1\eta$ and $t_{22}\eta=\mu_2\eta$,
then $t_{12}(u)\eta$ has weight $(q\mu_1,\qin\mu_2)$.
So it suffices to observe that
the set of $\U_q(\gl_2)$-weights of $L$ has
a maximal element with respect to the natural
ordering on the set of weights.

Due to the claim proved above, the highest vector $\ze$
belongs to $M$. It remains to show that
the vector $\ze$ is cyclic in $L$, that is,
the submodule $K=\U_q(\wh\gl_2)\tss\ze$ coincides with $L$.

Note that all $\U_q(\gl_2)$-weight spaces of $L$
are finite-dimensional.
Denote by $L^*$ the restricted dual vector space to $L$
which is the direct sum of the dual vector spaces
to the $\U_q(\gl_2)$-weight spaces of $L$.
We equip $L^*$ with
the $\U_q(\wh\gl_2)$-module structure
defined by
\beql{dualstarm}
(y\ts\om)(\eta)=\om(\varkappa(y)\ts\eta)
\quad
\text{for}
\quad
y\in\U_q(\wh\gl_2)
\Fand
\om\in L^*,\ \eta\in L,
\eeq
where $\varkappa$ is the involutive anti-automorphism of the algebra
$\U_q(\wh\gl_2)$,
defined by
\beql{varsigma}
\varkappa:\ t_{ij}(u)\mapsto \bar t_{3-i,3-j}(u^{-1}),
\qquad \bar t_{ij}(u)\mapsto t_{3-i,3-j}(u^{-1}).
\eeq
The latter is the composition of the automorphism \eqref{autotransp}
and the anti-automorphism \eqref{antiautopr}.
The anti-automorphism $\varkappa$ commutes with the
comultiplication $\Delta$ in the sense that
\ben
\Delta\circ\varkappa=
(\varkappa\ot\varkappa)\circ\Delta.
\een
This implies the isomorphism of $\U_q(\wh\gl_2)$-modules:
\beql{isomlstar}
L^*\cong L(\be^{-1}_1,\al^{-1}_1)\ot\dots\ot L(\be^{-1}_k,\al^{-1}_k).
\eeq
Moreover, the highest vector $\ze_i^*$ of the module
$L(\be^{-1}_i,\al^{-1}_i)\cong L(\al_i,\be_i)^*$ can be identified with
the element of $L(\al_i,\be_i)^*$ such that
$\ze_i^*(\ze_i)=1$
and $\ze_i^*(\eta_i)=0$ for all
$\U_q(\gl_2)$-weight vectors $\eta_i\in L(\al_i,\be_i)$
whose weights are different from the weight of $\ze_i$.

Now suppose on the contrary that the
submodule $K=\U_q(\wh\gl_2)\tss\ze$ of $L$ is proper.
The annihilator of $K$ defined by
\beql{annih}
\text{Ann\ts} K=\{\om\in L^*\ |\ \om(\eta)=0
\quad\text{for all}\quad\eta\in K\}
\eeq
is a submodule of the $\U_q(\wh\gl_2)$-module $L^*$, which does not
contain the vector $\ze_1^*\ot\dots\ot\ze_k^*$.
However, this contradicts the claim verified
in the first part of the proof,
because the condition on the parameters $\al_i$ and $\be_i$
remain satisfied after we replace each $\al_i$ by $\be^{-1}_i$ and each
$\be_i$ by $\al^{-1}_i$.
\epf

We can now describe the finite-dimensional irreducible
representations of the algebra $\U_q(\wh\gl_2)$.

\bth\label{thm:classgltwo}
The irreducible highest weight representation
$L\big(\nu_1(u),\nu_2(u);\bar\nu_1(u),\bar\nu_2(u)\big)$
of $\U_q(\wh\gl_2)$ is
finite-dimensional if and only if there exists a
polynomial $P(u)$ in $u$ with constant term $1$ such that
\beql{dpone}
\frac{\ve_1\ts\nu_1(u)}{\ve_2\ts\nu_2(u)}=q^{-\deg P}\cdot
\frac{P(u\tss q^2)}{P(u)}
=\frac{\ve_1\ts\bar\nu_1(u)}{\ve_2\ts\bar\nu_2(u)}
\eeq
for some $\ve_1,\ve_2\in\{-1,1\}$.
In this case $P(u)$ is unique.
\eth

\bpf
Suppose that the representation
$L\big(\nu_1(u),\nu_2(u);\bar\nu_1(u),\bar\nu_2(u)\big)$
is finite-dimensional. As was shown above,
we may assume without loss of generality that
the components of the highest weight have the form
\eqref{hwalbe}. Moreover, we may re-enumerate
the parameters $\al_i$ and $\be_i$ to satisfy
the conditions of Proposition~\ref{prop:tpr}.
By that proposition, the representation
$L\big(\nu_1(u),\nu_2(u);\bar\nu_1(u),\bar\nu_2(u)\big)$
is isomorphic to the tensor product \eqref{tenprqaff}.
Therefore, all ratios $\al_i/\be_i$ must have the form
$\pm q^{m_i}$, where each $m_i$ is a nonnegative integer.
Then the polynomial
\beql{dptpepr}
P(u)=\prod_{i=1}^k (1+\be_i^2\tss u)(1+\be_i^2\tss q^2\tss u)\dots
(1+\be_i^2\tss q^{2m_i-2}\tss u)
\eeq
satisfies \eqref{dpone} with an appropriate
choice of the signs $\ve_1,\ve_2\in\{-1,1\}$.

Conversely, suppose \eqref{dpone} holds for
a polynomial $P(u)=(1+\ga_1u)\dots(1+\ga_pu)$ and
some $\ve_1,\ve_2\in\{-1,1\}$.
Choose square roots $\be_i$ so that $\be_i^2=\ga_i$
for $i=1,\dots,p$ and
set
\ben
\bal
\mu_1(u)&=
(\be^{}_1q+\be^{-1}_1\qin u^{-1})\dots
(\be^{}_pq+\be^{-1}_p\qin u^{-1}),
\\
\mu_2(u)&=
(\be^{}_1+\be^{-1}_1u^{-1})\dots (\be^{}_p+\be^{-1}_pu^{-1}),
\\
\bar\mu_1(u)&=
(\be^{}_1q\ts u+\be^{-1}_1\qin)\dots (\be^{}_pq\ts u+\be^{-1}_p\qin),
\\
\bar\mu_2(u)&=
(\be^{}_1u+\be^{-1}_1)\dots (\be^{}_pu+\be^{-1}_p).
\eal
\een
Consider the tensor product module
\ben
L(\be^{}_1\tss q,\be^{}_1)\ot L(\be^{}_2\tss q,\be^{}_2)\ot\dots\ot
L(\be^{}_p\tss q,\be^{}_p)
\een
of $\U_q(\wh\gl_2)$.
This module is finite-dimensional and
the cyclic $\U_q(\wh\gl_2)$-span
of the tensor product of the highest vectors of
$L(\be^{}_i\tss q,\be^{}_i)$
is a highest weight module with the highest weight
$(\mu_1(u),\mu_2(u);\bar\mu_1(u),\bar\mu_2(u))$.
Hence, the
irreducible highest weight module
$L(\mu_1(u),\mu_2(u);\bar\mu_1(u),\bar\mu_2(u))$
is finite-dimensional.
Since
\ben
\frac{\mu_1(u)}{\mu_2(u)}=\frac{\ve_1\ts\nu_1(u)}{\ve_2\ts\nu_2(u)},
\qquad
\frac{\bar\mu_1(u)}{\bar\mu_2(u)}
=\frac{\ve_1\ts\bar\nu_1(u)}{\ve_2\ts\bar\nu_2(u)}
\een
there exist automorphisms of $\U_q(\wh\gl_2)$ of the form
\eqref{automsln} and \eqref{aautoaffsign}
such that their composition with the representation
$L(\mu_1(u),\mu_2(u);\bar\mu_1(u),\bar\mu_2(u))$ is isomorphic to
the irreducible highest weight representation
$L\big(\nu_1(u),\nu_2(u);\bar\nu_1(u),\bar\nu_2(u)\big)$.
Thus, the latter
is also finite-dimensional.

The uniqueness of $P(u)$ is easily verified.
\epf

The above arguments imply that,
up to twisting with
an automorphism of the form \eqref{automsln},
every finite-dimensional irreducible
representation of $\U_q(\wh\gl_2)$
is isomorphic to a tensor product representation
of the form \eqref{tenprqaff}.
We will now establish a criterion of irreducibility of
such representations which we will use
in Sec.~\ref{subsec:reptw} below. It is essentially
a version
of the well-known results;
see \cite[Ch.~12]{cp:gq}, \cite{t:im}.

We will define a $q$-{\it string\/} to be
any subset of $\CC$ of the form $\{\be,\be\tss q,\dots,\be\tss q^p\}$,
where $\be\ne 0$ and $p$ is a nonnegative integer.
Since $q$ is not a root of unity, the $q$-strings
$\{\be,\be\tss q,\dots,\be\tss q^p\}$ and
$\{-\be,-\be\tss q,\dots,-\be\tss q^p\}$ have no common
elements. Their union will be called a $q$-{\it spiral\/}.
Two $q$-spirals $S_1$ and $S_2$ are
{\it in general position\/}
if either
\begin{itemize}
\item[(i)] $S_1\cup S_2$ is not a $q$-spiral; or
\item[(ii)] $S_1\subset S_2$, or $S_2\subset S_1$.
\end{itemize}

Given a pair of nonzero complex numbers
$(\al,\be)$ with $\al/\be=\pm q^{m}$ and $m\in\ZZ_+$
the corresponding $q$-spiral is defined as
\ben
S_q(\al,\be)=
\{\be,\be\tss q,\dots,\be\ts q^{m-1}\}\cup
\{-\be,-\be\tss q,\dots,-\be\ts q^{m-1}\}.
\een
If $\al=\be$, then the set $S_q(\al,\be)$ is regarded to be
empty. Note that changing sign of $\al$ or $\be$ does not
affect the $q$-spiral.

Denote by $L$ the tensor product \eqref{tenprqaff},
where for all ratios we have $\al_i/\be_i=\pm q^{m_i}$
for some nonnegative integers $m_i$.

\bco\label{cor:irrcrit}
The representation $L$
of $\U_q(\wh\gl_2)$ is irreducible
if and only if the $q$-spirals
$S_q(\al_1,\be_1),\dots,S_q(\al_k,\be_k)$ are pairwise in
general position.
\eco

\bpf
Suppose that the $q$-spirals are pairwise in general position and
assume first that the nonnegative integers $m_i$
satisfy the inequalities
$m_1\leqslant\dots\leqslant m_k$. This implies that
the condition of Proposition~\ref{prop:tpr}
on the parameters $\al_i$ and $\be_i$
holds. Indeed, if this is not the case, then $\al_r/\be_s=\pm q^p$
for some $i\leqslant r,s\leqslant k$ and a nonnegative
integer $p$ with $p<m_i$. By our assumption, $r\ne s$.
Suppose that $r>s$. Then $m_i\leqslant m_s$ so we may
assume that $s=i$. The condition $\al_r=\pm \be_i\tss q^p$
means that $\al_r$ belongs to the $q$-spiral $S_q(\al_i,\be_i)$.
Hence, the $q$-spirals $S_q(\al_i,\be_i)$ and $S_q(\al_r,\be_r)$
are not in general position, a contradiction.
The opposite inequality $r<s$ leads to a similar contradiction.

Thus,
$L$ is irreducible by Proposition~\ref{prop:tpr}.
It is easy to verify (cf.
\cite[Prop.~3.2.10]{m:yc}) that
any permutation of the tensor factors yields
an isomorphic irreducible
representation.

Conversely, let $k=2$ and
let $L(\al_1,\be_1)\ot L(\al_2,\be_2)$
be irreducible. Suppose that the $q$-spirals
$S_q(\al_1,\be_1)$ and $S_q(\al_2,\be_2)$
are not in general position. Then
the $q$-spirals
$S_q(\al_1,\be_2)$ and $S_q(\al_2,\be_1)$ are in general position.
Hence, the representation $L(\al_1,\be_2)\ot L(\al_2,\be_1)$ of
$\U_q(\wh\gl_2)$
is irreducible due to the first part of the proof.
Comparing the dimension of this representation
with the dimension of $L(\al_1,\be_1)\ot L(\al_2,\be_2)$
we come to a contradiction.

The case of general $k\geqslant 3$ is reduced to $k=2$
by permuting the tensor factors in
\eqref{tenprqaff}, if necessary. Indeed,
if $L$ is irreducible, but
a pair of $q$-spirals $S_q(\al_i,\be_i)$ and $S_q(\al_j,\be_j)$
is not in general position, then
we may assume that $i$ and $j$ are adjacent. However,
the representation $L(\al_i,\be_i)\ot L(\al_j,\be_j)$ of
$\U_q(\wh\gl_2)$
is reducible as shown above. This implies that
$L$ is reducible, a contradiction.
\epf

\bre\label{rem:gauss}
An isomorphism between the $RTT$-presentation
of the algebra $\U_q(\wh\gl_N)$ and its new realization
is provided in \cite{df:it} by using the Gauss decomposition
of the matrices $T(u)$ and $\overline T(u)$; cf.
\cite{bk:pp}, \cite{fm:ha}.
Thus, Theorem~\ref{thm:classgltwo} provides a
description of finite-dimensional irreducible
representations of the algebras $\U_q(\wh\gl_2)$
and $\U_q(\wh\sll_2)$
in terms of the new realization
via this isomorphism.
This argument is alternative to \cite{cp:qaa} and it is
straightforward to apply this description to prove
the classification theorem for
finite-dimensional irreducible
representations of an arbitrary
quantum affine algebra $\U_q(\wh\agot)$; cf.
the $A$ type case considered below.
\qed
\ere

\subsection{Representations of $\U_q(\wh\gl_N)$}

The evaluation homomorphism $\pi:\U_q(\wh\gl_N)\to\U_q(\gl_N)$
defined in \eqref{eval} allows us to regard
any $\U_q(\gl_N)$-module as a $\U_q(\wh\gl_N)$-module.
In particular, we thus obtain the {\it evaluation modules\/}
$L(\mu)$ over $\U_q(\wh\gl_N)$; see Sec.~\ref{subsec:que}.

As we pointed out in the beginning of Sec.~\ref{sec:qaff},
in order to describe all finite-dimensional
irreducible representations of the algebra $\U_q(\wh\gl_N)$,
we need to determine for which highest weights
$(\nu(u);\bar\nu(u))$ the representation $L(\nu(u);\bar\nu(u))$
is finite-dimensional. These conditions are provided
by the following theorem which is essentially equivalent
to \cite[Theorem~12.2.6]{cp:gq}.

\bth\label{thm:clgen}
The irreducible highest weight representation
$L\big(\nu(u);\bar\nu(u)\big)$
of $\U_q(\wh\gl_N)$ is
finite-dimensional if and only if there exist polynomials
$P_1(u),\dots,P_{N-1}(u)$ in $u$,
all with constant term $1$,
such that
\beql{nunupone}
\frac{\ve_i\ts\nu_i(u)}{\ve_{i+1}\ts\nu_{i+1}(u)}=q^{-\deg P_i}\cdot
\frac{P_i(u\tss q^2)}{P_i(u)}
=\frac{\ve_i\ts\bar\nu_i(u)}{\ve_{i+1}\ts\bar\nu_{i+1}(u)}
\eeq
for $i=1,\dots,N-1$ and some
$\ve_j\in\{-1,1\}$.
The polynomials $P_1(u),\dots,P_{N-1}(u)$
are determined uniquely while the tuple $(\ve_1,\dots,\ve_N)$
is determined uniquely, up to the simultaneous
change of sign of the $\ve_i$.
\eth

\bpf
Suppose that the representation $L\big(\nu(u);\bar\nu(u)\big)$
is finite-dimensional.
Let us fix $k\in\{0,\dots,N-2\}$ and let $\U_q(\wh\gl_2)$
act on $L\big(\nu(u);\bar\nu(u)\big)$ via
the homomorphism $\U_q(\wh\gl_2)\to\U_q(\wh\gl_N)$
which sends $t_{ij}(u)$ and $\bar t_{ij}(u)$
to $t_{k+i,k+j}(u)$ and $\bar t_{k+i,k+j}(u)$,
respectively, for any $i,j\in\{1,2\}$.
The cyclic $\U_q(\wh\gl_2)$-span of the highest vector
of $L\big(\nu(u);\bar\nu(u)\big)$
is a highest weight representation of the algebra
$\U_q(\wh\gl_2)$ with the highest weight
\ben
\big(\nu_{k+1}(u),\nu_{k+2}(u);\bar\nu_{k+1}(u),\bar\nu_{k+2}(u)\big).
\een
Its irreducible quotient
is finite-dimensional, and so
the required conditions follow from Theorem~\ref{thm:classgltwo}.

Conversely, taking into account the automorphisms \eqref{automsln}
and \eqref{aautoaffsign},
it is enough to show that given any
set of polynomials $P_1(u),\dots,P_{N-1}(u)$ in $u$
with constant terms equal to $1$,
there exists an irreducible finite-dimensional representation
whose highest weight satisfies
\eqref{nunupone}. Such a representation can be constructed
by using the following inductive procedure.
Consider the irreducible highest weight representation
$L(\la)$ of $\U_q(\gl_N)$ with the highest weight
\ben
\la=(d\tss q^{m_1},\dots,d\tss q^{m_N}),
\een
where $d$ is a nonzero complex number
and the integers $m_i$ satisfy
$m_1\geqslant\dots\geqslant m_N$. This representation
is finite-dimensional and we regard it as the evaluation
module over $\U_q(\wh\gl_N)$ by using
the homomorphism \eqref{eval}.
By the first part of the proof
we can associate a family of polynomials
$P_1(u),\dots,P_{N-1}(u)$ to any finite-dimensional
representation $L(\nu(u);\bar\nu(u))$.
Let $\ze$ and $\xi$ be the highest vectors
of the representations $L(\la)$ and $L(\nu(u);\bar\nu(u))$,
respectively, and equip $L(\la)\ot L(\nu(u);\bar\nu(u))$
with the $\U_q(\wh\gl_N)$-module structure by
using the coproduct \eqref{copraff}. It is easily verified that
the cyclic span $\U_q(\wh\gl_N)(\ze\ot\xi)$
is a highest weight representation of
$\U_q(\wh\gl_N)$ such that
\ben
\bal
t_{ii}(u)(\ze\ot\xi)&=\big(d\ts q^{m_i}
+d^{-1}\tss q^{-m_i}u^{-1}\big)\ts \nu_i(u)\ts(\ze\ot\xi),\\
\bar t_{ii}(u)(\ze\ot\xi)&=\big(d^{-1}\tss q^{-m_i}
+d\ts q^{m_i}u\big)\ts \bar\nu_i(u)\ts(\ze\ot\xi).
\eal
\een
Hence,
the irreducible quotient of this representation
corresponds to
the family of polynomials $Q_1(u)P_1(u),\dots,Q_{N-1}(u)P_{N-1}(u)$,
where
\ben
Q_i(u)=(1+d^{2}q^{2m_{i+1}}u)(1+d^{2}q^{2m_{i+1}+2}u)
\dots (1+d^{2}q^{2m_{i}-2}u),\qquad i=1,\dots,N-1.
\een
Starting from the trivial representation
of $\U_q(\wh\gl_N)$ and choosing appropriate
parameters $d$ and $m_i$ we will be able to produce
a finite-dimensional highest weight representation
associated with an arbitrary family of
polynomials $P_1(u),\dots,P_{N-1}(u)$ by iterating
this construction.
The last statement of the theorem is easily verified.
\epf

The polynomials $P_1(u),\dots,P_{N-1}(u)$ introduced
in Theorem~\ref{thm:clgen} are called the {\it Drinfeld
polynomials\/} of the representation $L\big(\nu(u);\bar\nu(u)\big)$.
Moreover,
any finite-dimensional irreducible representation of
$\U_q(\wh\gl_N)$ is obtained from a
representation $L\big(\nu(u);\bar\nu(u)\big)$
associated with the tuple $(\ve_1,\dots,\ve_N)=(1,\dots,1)$
by twisting with
an automorphism of the form \eqref{aautoaffsign}.
Note also that the evaluation
module $L(\mu)$ with
\ben
\mu_i=q^{m_i},\qquad i=1,\dots,N,
\een
where $m_1\geqslant m_2\geqslant\dots\geqslant m_N$
are arbitrary integers, is a representation
associated with the tuple $(1,\dots,1)$.
Its Drinfeld polynomials are given by
\ben
P_i(u)=(1+q^{2m_{i+1}}\tss u)(1+q^{2m_{i+1}+2}\tss u)
\dots(1+q^{2m_{i}-2}\tss u),
\een
for $i=1,\dots,N-1$.

A description of finite-dimensional irreducible
representations of the extended algebra
$\U^{\tss\text{\rm ext}}_q(\wh\gl_N)$ can be easily obtained
from that of the quantum affine algebra
$\U_q(\wh\gl_N)$. Namely, every finite-dimensional
irreducible representation of $\U^{\tss\text{\rm ext}}_q(\wh\gl_N)$
is isomorphic to the highest weight representation
$L(\nu(u);\bar\nu(u))$. The latter is defined in the same way as for
the algebra $\U_q(\wh\gl_N)$ except that
the relations $\nu_i^{(0)}\bar\nu_i^{(0)}=1$
for the series \eqref{nui} are replaced by
the condition that all constants $\nu_i^{(0)}$
and $\bar\nu_i^{(0)}$ are nonzero.
We have the following corollary of Theorem~\ref{thm:clgen}.

\bco\label{cor:clgenext}
The irreducible highest weight representation
$L\big(\nu(u);\bar\nu(u)\big)$
of the algebra $\U^{\tss\text{\rm ext}}_q(\wh\gl_N)$ is
finite-dimensional if and only if there exist polynomials
$P_1(u),\dots,P_{N-1}(u)$ in $u$,
all with constant term $1$, and nonzero
constants $\phi_1,\dots,\phi_N$ such that
\beql{nunuponeext}
\frac{\phi_i\ts\nu_i(u)}{\phi_{i+1}\ts\nu_{i+1}(u)}=q^{-\deg P_i}\cdot
\frac{P_i(u\tss q^2)}{P_i(u)}
=\frac{\phi_i\ts\bar\nu_i(u)}{\phi_{i+1}\ts\bar\nu_{i+1}(u)}
\eeq
for $i=1,\dots,N-1$.
The polynomials $P_1(u),\dots,P_{N-1}(u)$
are determined uniquely while the tuple $(\phi_1,\dots,\phi_N)$
is determined uniquely, up to a common factor.
\eco

\bpf
By twisting the representation $L\big(\nu(u);\bar\nu(u)\big)$
by an appropriate automorphism of the algebra
$\U^{\tss\text{\rm ext}}_q(\wh\gl_N)$ of the form
\eqref{aautoaffsignext}, we can get
the representation where
all central elements $t_{ii}^{(0)}\ts \bar t_{ii}^{\ts(0)}$,
$i=1,\dots,N$, act as the identity operators.
Therefore, due to \eqref{uextepi}, we get
the irreducible highest weight representation
of the algebra $\U_q(\wh\gl_N)$, such that the components
of the highest weight have the form
$\phi_i\tss\nu_i(u)$ and $\phi_i\tss\bar\nu_i(u)$.
Now all statements follow from
Theorem~\ref{thm:clgen}.
\epf

\section{Representations of the twisted $q$-Yangians}
\label{sec:twqyang}
\setcounter{equation}{0}

We will combine the approaches of Sec.~\ref{sec:qaff}
and \cite[Ch.~4]{m:yc} to classify
the finite-dimensional irreducible representations
of the twisted $q$-Yangians $\Y'_q(\spa_{2n})$.
As with the quantum affine algebras, the case $n=1$
will play a key role.
We start by proving some general results about
highest weight representation of the twisted $q$-Yangians.

\subsection{Highest weight representations}

As we recalled in Sec.~\ref{subsec:twqyang},
the twisted $q$-Yangian $\Y'_q(\spa_{2n})$
can be defined as the algebra generated
by the elements $s_{ij}^{(r)}$ with $r\geqslant 0$
and $1\leqslant i,j\leqslant 2n$
and by the elements $s_{i,i+1}^{(0)-1}$ with $i=1,3,\dots,2n-1$.
The defining relations are written in terms
of the generating series \eqref{scoeffs}
and they take the form \eqref{draffsu}
together with the relations
\ben
s_{ij}^{(0)}=0\quad\text{for}\quad i<j \quad\text{unless}\ \
j=i+1\quad\text{with}\quad i\ \  \text{odd}
\een
and
\beql{inverti}
s_{i,i+1}^{(0)}\tss s_{i,i+1}^{(0)-1}
=s_{i,i+1}^{(0)-1}\tss s_{i,i+1}^{(0)}=1,\qquad
i=1,3,\dots,2n-1.
\eeq
Observe that given any formal series $g(u)$ in $u^{-1}$
of the form
\ben
g(u)=g_0+g_1u^{-1}+g_2u^{-2}+\dots,\qquad g_0\ne 0,
\een
the mapping
\beql{multautotw}
s_{ij}(u)\mapsto g(u)\ts s_{ij}(u)
\eeq
defines an automorphism of the algebra $\Y'_q(\spa_{2n})$.

Given any tuple $(\psi_1,\dots,\psi_{2n})$ of nonzero
complex numbers, the mapping
\beql{aautoaffsigntw}
s_{ij}(u)\mapsto \psi_i\ts\psi_j\ts s_{ij}(u)
\eeq
defines another automorphism of the algebra
$\Y'_q(\spa_{2n})$.

Furthermore, the mapping
\beql{antiautw}
\varkappa:s_{ij}(u)\mapsto s_{2n-j+1,2n-i+1}(u)
\eeq
defines an involutive anti-automorphism of the algebra
$\Y'_q(\spa_{2n})$. This can be verified
directly from the defining relations. Alternatively,
one can show that the mappings \eqref{autotransp}
and \eqref{antiautopr} respectively define an automorphism
and anti-automorphism of the extended algebra
$\U^{\tss\text{\rm ext}}_q(\wh\gl_{2n})$, and their
composition $\varkappa$ preserves the subalgebra
$\Y'_q(\spa_{2n})$.

We will use the elements $\bar s_{ij}^{(r)}$
of the algebra $\Y'_q(\spa_{2n})$ which are defined
as the coefficients of the power series $\bar s_{ij}(u)$
in $u$; see \eqref{matsbar}. The relationship
between the elements is given by the
formulas\footnote{The corresponding relation (3.62) in \cite{mrs:cs}
should be corrected by swapping $\delta_{i<j}$ and $\delta_{j<i}$.}
\begin{multline}\label{relsbarssympl}
(u^{-1}q-uq^{-1})\ts\bar s_{ij}(u) = \\
(uq^{\delta_{ij}}-u^{-1}q^{-\delta_{ij}})\, s_{ji}(u^{-1})
+(q-q^{-1})(u^{-1}\delta_{i<j}+u\tss\delta_{j<i})\ts s_{ij}(u^{-1}).
\end{multline}

For the rest of this section
we will suppose that the complex number $q$
is nonzero and not a root of unity.
Recall the function
$\vs: \{ 1, 2, \ldots, 2n\} \to \{\pm 1, \pm 3,\ldots, \pm(2n-1) \}$
defined in \eqref{nu}.

\bde\label{def:hwtw}
A representation $V$ of the algebra $\Y'_q(\spa_{2n})$ is called
a {\it highest weight representation\/}
if $V$ is generated
by a nonzero vector $\xi$ (the {\it highest vector\/})
such that
\begin{alignat}{2}
s_{kl}(u)\ts\xi&=0,\qquad
&&\text{for} \quad
\vs(k)+\vs(l)>0,
\non\\
s_{2i,2i-1}(u)\ts\xi&=\mu_i(u)\ts\xi,\qquad
\qquad &&\text{for} \quad 1\leqslant i\leqslant n,
\label{hwcond}\\
\bar s_{2i,2i-1}(u)\ts\xi&=\bar \mu_i(u)\ts\xi,\qquad
\qquad &&\text{for} \quad 1\leqslant i\leqslant n,
\non
\end{alignat}
where $\mu(u)=(\mu_1(u),\dots,\mu_n(u))$ and
$\bar\mu(u)=(\bar\mu_1(u),\dots,\bar\mu_n(u))$ are certain $n$-tuples
of formal power series in $u^{-1}$ and $u$, respectively:
\beql{muitw}
\mu_i(u)=\sum_{r=0}^{\infty} \mu_i^{(r)}\tss u^{-r},
\qquad
\bar\mu_i(u)=\sum_{r=0}^{\infty} \bar\mu_i^{(r)}\tss u^{r}.
\eeq
Due to \eqref{relsbarssympl},
the first relation in \eqref{hwcond}
is equivalent to
\beql{hwannih}
s_{ij}(u)\ts\xi=\bar s_{ij}(u)\ts\xi=0
\eeq
for $j=1,3,\dots,2n-1$ and $i=1,2,\dots,j$.
\qed
\ede

The definition of the highest weight representation
is consistent with a particular choice
of the positive root system of type $C_n$. Namely,
the root system $\Phi$ is the subset
of vectors in $\RR^n$ of the form
\ben
{}\pm 2\ts\ve_i\quad\text{with}\quad 1\leqslant i\leqslant n\Fand
{}\pm\ve_i\pm\ve_j\quad\text{with}\quad 1\leqslant i<j\leqslant n,
\een
where $\ve_i$ denotes the $n$-tuple which has
$1$ on the $i$-th position
and zeros elsewhere. Partition this set into positive
and negative roots $\Phi=\Phi^+\cup(-\Phi^+)$, where the set
of positive roots $\Phi^+$ consists of the vectors
\beql{posroots}
2\ts\ve_i\quad\text{with}\quad 1\leqslant i\leqslant n\Fand
\ve_i+\ve_j,\quad {}-\ve_i+\ve_j\quad\text{with}
\quad 1\leqslant i<j\leqslant n.
\eeq
We will regard $\U'_q(\spa_{2n})$ as a subalgebra of
$\Y'_q(\spa_{2n})$ defined via the embedding
\beql{embedd}
s_{ij}\mapsto s^{(0)}_{ij}.
\eeq
By \eqref{siiponeskl} we have
\beql{siiponesklnew}
s^{}_{i,i+1}\ts s^{}_{kl}(u)=
q^{\de_{ik}+\de_{il}-\de_{i+1,k}-\de_{i+1,l}}\ts
s^{}_{kl}(u)\ts s^{}_{i,i+1},
\eeq
for any $i=1,3,\dots,2n-1$. Hence, the generating series
$s_{kl}(u)$ with $\vs(k)+\vs(l)>0$ can be associated
with the elements of $\Phi^+$ listed in \eqref{posroots}
as follows:
\beql{assosroots}
\bal
2\ts\ve_i\longleftrightarrow s_{2i-1,2i-1}(u),\qquad
\ve_i+\ve_j &\longleftrightarrow
\{s_{2i-1,2j-1}(u),s_{2j-1,2i-1}(u)\},\\
\qquad
-\ve_i+\ve_j &\longleftrightarrow \{s_{2i,2j-1}(u),s_{2j-1,2i}(u)\}.
\eal
\eeq
Here the commutative subalgebra of $\U'_q(\spa_{2n})$
generated by the elements $s^{}_{i,i+1}$ with
$i=1,3,\dots,2n-1$ plays the role of a Cartan subalgebra; cf.
\cite{m:rtq}.

\bth\label{thm:fdhwtw}
Any finite-dimensional irreducible representation $V$
of the algebra $\Y'_q(\spa_{2n})$ is
a highest weight representation.
Moreover, $V$ contains a unique,
up to a constant factor, highest vector.
\eth

\bpf
We use a standard argument with some necessary modifications;
cf.~\cite[Sec.~12.2]{cp:gq} and \cite[Sec.~4.2]{m:yc}.
Set
\beql{vo}
\bal
V^{\tss0}=\{\eta\in V\ |\ s_{ij}(u)\ts\eta
&=\bar s_{ij}(u)\ts\eta=0,\\
j&=1,3,\dots,2n-1\Fand i=1,2,\dots,j\tss\}.
\eal
\eeq
Equivalently, $V^{\tss0}$ is spanned by the vectors
annihilated by all operators $s_{kl}(u)$
such that $\vs(k)+\vs(l)>0$.
Let us show that $V^{\tss0}$ is nonzero.
The operators $s^{}_{2i-1,2i}=s_{2i-1,2i}^{(0)}$ with $i=1,\dots,n$
pairwise commute on $V$ and so $V$ contains a common
eigenvector $\theta$ for these operators:
\ben
s_{2i-1,2i}\ts\theta=\rho_i\ts\theta,\qquad i=1,\dots,n.
\een
By \eqref{siiponesklnew}, every coefficient of the series
$s_{kl}(u)\ts\theta$ with $\vs(k)+\vs(l)>0$ is again
a common eigenvector for the operators $s^{}_{2i-1,2i}$,
whose eigenvalues have the form
$\rho_i\ts q^{\al_i}$, $i=1,\dots,n$, where
$\al=\al_1\ve_1+\dots+\al_n\ve_n$ is the element of $\Phi^+$
associated with $s_{kl}(u)$ by \eqref{assosroots}.
Since the sets of eigenvalues
obtained in this way are distinct and
$\dim V<\infty$, we can conclude that
there exists a nonzero element of $V$ annihilated by
all operators $s_{kl}(u)$ with $\vs(k)+\vs(l)>0$.
Thus, $V^{\tss0}\ne\{0\}$.

Next, we show that the subspace $V^{\tss0}$ is invariant
with respect to the action of all operators
$s_{b+1,b}(v)$ and $\bar s_{b+1,b}(v)$ with odd $b$.
Let us show first that if $\eta\in V^{\tss0}$, then for any odd
$a$ and $i\leqslant a$ we have
\beql{siasbbpone}
s_{ia}(u)\tss s_{b+1,b}(v)\tss \eta =0.
\eeq
If $i\leqslant b$, then this follows by
\eqref{draffsu} with $j=b+1$.
If $i>b$, then $b+1\leqslant i\leqslant a$ and
\eqref{siasbbpone} follows by
the application of \eqref{draffsu},
where $i,j,a,b$ are respectively replaced
with $b+1,i,b,a$, and $u$ is swapped with $v$.

Furthermore, three more relations of the form
\eqref{siasbbpone} where $s_{ia}(u)$
is replaced with $\bar s_{ia}(u)$ or
$s_{b+1,b}(v)$ is replaced with $\bar s_{b+1,b}(v)$,
are verified by exactly the same argument with the use
of the corresponding relation in
\eqref{rsrsbar} instead of \eqref{rsrsaffs}.

A similar argument shows that all operators
$s_{a+1,a}^{(r)}$ and $\bar s_{a+1,a}^{\ts(r)}$
on the space $V^{\tss0}$
with odd $a$ and $r\geqslant 0$ pairwise commute.
Indeed, suppose that both $a$ and $b$ are odd
and $a\leqslant b$.
The application of \eqref{draffsu}
with $i=a+1$ and $j=b+1$ proves that
all operators $s_{a+1,a}^{(r)}$ pairwise commute.
The remaining commutativity relations are verified
in the same way with the use of \eqref{rsrsbar}.

Thus, the operators $s_{a+1,a}^{(r)}$ and $\bar s_{a+1,a}^{\ts(r)}$
on the space $V^{\tss0}$
with odd $a$ and $r\geqslant 0$ are
simultaneously diagonalizable.  We let
$\xi$ be a simultaneous eigenvector for all
these operators.  Then $V = \Y'_q(\spa_{2n})\tss \xi$
by the irreducibility of $V$, and
$\xi$ is a highest weight vector.
By considering the $\U'_q(\spa_{2n})$-weights
of $V$ and using \eqref{siiponesklnew} we may also
conclude that $\xi$ is determined uniquely,
up to a constant factor.
\epf

Consider now the tuples $\mu(u)=(\mu_1(u),\dots,\mu_n(u))$ and
$\bar\mu(u)=(\bar\mu_1(u),\dots,\bar\mu_n(u))$
of arbitrary formal power series
of the form \eqref{muitw}.

The {\it Verma module\/}
$M(\mu(u);\bar\mu(u))$
over the twisted $q$-Yangian $\Y'_q(\spa_{2n})$
is the quotient of $\Y'_q(\spa_{2n})$
by the left ideal generated by all coefficients of the series
$s_{kl}(u)$ for $\vs(k)+\vs(l)>0$,
$s_{2i,2i-1}(u)-\mu_i(u)$ and $\bar s_{2i,2i-1}(u)-\bar\mu_i(u)$ for
$i=1,\dots,n$.

Clearly, the Verma module $M(\mu(u);\bar\mu(u))$
is a highest weight representation
of $\Y'_q(\spa_{2n})$ with the highest weight $(\mu(u);\bar\mu(u))$.
Moreover, any highest weight representation
with the same highest weight
is isomorphic to a quotient
of $M(\mu(u);\bar\mu(u))$.

The Poincar\'e--Birkhoff--Witt theorem for the algebra
$\Y'_q(\spa_{2n})$ (Proposition~\ref{prop:pbwtw})
implies that the ordered monomials of the form
\beql{monompbwtw}
s_{i_1j_1}^{(r_1)}\dots
s_{i_mj_m}^{(r_m)}\ts 1,\qquad m\geqslant 0,\qquad
\vs(i_a)+\vs(j_a)<0,
\eeq
form a basis of $M(\mu(u);\bar\mu(u))$. Moreover, using
\eqref{siiponesklnew} and considering the
weights of $M(\mu(u);\bar\mu(u))$
with respect to the operators $s_{i,i+1}$ with odd $i$,
we derive that the Verma module $M(\mu(u);\bar\mu(u))$
possesses a unique maximal proper submodule $K$.

The {\it irreducible highest weight representation\/}
$V(\mu(u);\bar\mu(u))$ of
$\Y'_q(\spa_{2n})$ with the highest weight $(\mu(u);\bar\mu(u))$
is defined as the
quotient of the Verma module $M(\mu(u);\bar\mu(u))$ by
the submodule $K$.

Due to Theorem~\ref{thm:fdhwtw}, all finite-dimensional
irreducible representations of the twisted $q$-Yangian $\Y'_q(\spa_{2n})$
have the form $V(\mu(u);\bar\mu(u))$ for a
certain highest weight $(\mu(u);\bar\mu(u))$.
Hence, in order to classify such representations
it remains to describe the set of highest weights
$(\mu(u);\bar\mu(u))$ such that
$V(\mu(u);\bar\mu(u))$ is finite-dimensional.
As with the quantum affine algebra $\U_q(\wh\gl_N)$
(see Sec.~\ref{sec:qaff}), a key role will be played by
the particular case $\Y'_q(\spa_{2})$ which we consider
in the next section.

\subsection{Representations of $\Y'_q(\spa_{2})$}
\label{subsec:reptw}

The irreducible highest weight representations
$V(\mu(u);\bar\mu(u))$ of $\Y'_q(\spa_{2})$
are parameterized by formal series of the form
\beql{muusp}
\bal
\mu(u)&=\mu^{(0)}+\mu^{(1)}u^{-1}+\mu^{(2)}u^{-2}+\dots,\\
\bar\mu(u)&=\bar\mu^{\tss(0)}+\bar\mu^{\tss(1)}u
+\bar\mu^{\tss(2)}u^{2}+\dots,\qquad
\mu^{(r)},\bar\mu^{\tss(r)}\in\CC.
\eal
\eeq
The highest vector $\xi$ of $V(\mu(u);\bar\mu(u))$
satisfies the conditions
\ben
s_{11}(u)\ts\xi=0,\qquad
s_{21}(u)\ts\xi=\mu(u)\ts\xi,\qquad
\bar s_{21}(u)\ts\xi=\bar \mu(u)\ts\xi.
\een
Note that due to the relations \eqref{relsbarssympl},
the vector $\xi$ is also an eigenvector for the
operators $s_{12}(u)$ and $\bar s_{12}(u)$,
\beql{eigen}
\bal
s_{12}(u)\ts\xi&=\mu'(u)\ts\xi,\qquad
\mu'(u)
=\frac{(q^2-1)\tss\mu(u)+(1-u^2q^2)\ts\bar\mu(u^{-1})}
{q(u^2-1)},\\
\bar s_{12}(u)\ts\xi&=\bar\mu^{\ts\prime}(u)\ts\xi,\qquad
\bar\mu^{\ts\prime}(u)
=\frac{(q^2-1)\tss\bar\mu(u)+(1-u^2q^2)\ts\mu(u^{-1})}
{q(u^2-1)},
\eal
\eeq
where $\mu'(u)$ and $\bar\mu^{\ts\prime}(u)$ are regarded
as formal series in $u^{-1}$ and $u$, respectively.

\bpr\label{prop:nonzero}
If the representation
$V(\mu(u);\bar\mu(u))$ of $\Y'_q(\spa_{2})$
is finite-dimensional, then both coefficients
$\mu^{(0)}$ and $\bar\mu^{\ts(0)}$ in \eqref{muusp}
are nonzero.
\epr

\bpf
The constant term $\mu^{\prime\ts(0)}$ of the series $\mu'(u)$
is nonzero due to the relation \eqref{inverti}.
By \eqref{eigen} we have $\mu^{\prime\ts(0)}=-q\tss\bar\mu^{(0)}$
and so, $\bar\mu^{(0)}\ne 0$. Furthermore,
consider the restriction of $V(\mu(u);\bar\mu(u))$
to the subalgebra $\U'_q(\spa_{2})$ defined by the
embedding \eqref{embedd}. The cyclic span $\U'_q(\spa_{2})\ts\xi$
of the highest vector
is a finite-dimensional representation
of the subalgebra with the highest weight
$(\mu^{(0)};\mu^{\prime\ts(0)})$. However, as was pointed out
in the proof of Proposition~\ref{prop:fdco}, this
implies that $\mu^{(0)}\ne0$.
\epf

The following is an analogue
of Proposition~\ref{prop:polyn} and its proof follows
a similar approach.

\bpr\label{prop:polyntw}
Suppose that
$\dim V(\mu(u);\bar\mu(u))<\infty$.
Then there exists
a polynomial $Q(u)$ in $u$ of even degree
with the constant term equal to $1$
such that
\beql{ratfun}
\frac{\bar\mu(u^{-1})}{\mu(u)}=
\frac{u^{\deg Q}\ts Q(u^{-1})}{Q(u)}.
\eeq
\epr

\bpf
By twisting the representation $V(\mu(u);\bar\mu(u))$
with an
automorphism of $\Y'_q(\spa_{2})$ of the form
\eqref{multautotw} we get a representation
isomorphic to $V(g(u)\mu(u);g(u^{-1})\bar\mu(u))$.
Since $\mu^{\prime\ts(0)}\ne 0$, we may consider
such an automorphism with
\ben
g(u)=\big(\mu'(u)+\qin u^{-2}\mu(u)\big)^{-1}.
\een
Hence, we may assume
without loss of generality that
the highest weight of the representation $V(\mu(u);\bar\mu(u))$
satisfies the condition $\mu'(u)+\qin u^{-2}\mu(u)=1$.

As in the proof of Proposition~\ref{prop:polyn},
the assumption $\dim V(\mu(u);\bar\mu(u))<\infty$ implies that
\beql{singvect}
\sum_{l=0}^k c_l\ts s^{(l)}_{22}\ts \xi=0
\eeq
for some $k\geqslant 0$ and
some $c_l\in\CC$ with $c_k\ne 0$.

On the other hand,
the defining relations \eqref{draffsu} imply
\ben
\bal
(u^{-1}-v^{-1})&(1-u^{-1}v^{-1})\ts s_{11}(u)\ts s_{22}(v)\ts\xi\\
=(\qin-q)\Big({}&{}u^{-1}\big(\mu'(u)+q\ts\mu(u)\big)
\big(\mu'(v)+\qin v^{-2}\mu(v)\big)\\
{}-{}&{}v^{-1}\big(\mu'(v)+q\ts\mu(v)\big)
\big(\mu'(u)+\qin u^{-2}\mu(u)\big)\Big)\xi.
\eal
\een
Taking into account the assumption $\mu'(u)+\qin u^{-2}\mu(u)=1$,
we can write these relations as
\beql{relsf}
(1-u^{-1}v^{-1})\ts s_{11}(u)\ts s_{22}(v)\ts\xi
=(\qin-q)\ts
\frac{u^{-1}\rho(u)-v^{-1}\rho(v)}{u^{-1}-v^{-1}}\ts\xi,
\eeq
where
\ben
\rho(u)=\mu'(u)+q\ts\mu(u).
\een
Write
\beql{muu}
\rho(u)=\sum_{r=0}^{\infty} \rho^{(r)}\ts u^{-r}.
\eeq
Divide both sides of \eqref{relsf} by $1-u^{-1}v^{-1}$
and compare
the coefficients of $u^{-m}v^{-l}$. This gives
\ben
s^{(m)}_{11}\tss s^{(l)}_{22}\ts \xi=
(\qin-q)\ts \wt\rho^{\ts(m,\tss l)}\ts \xi,
\een
where
\ben
\wt\rho^{\ts(m,\tss l)}=
\sum_{i=0}^{\min\{m,\tss l\}}\rho^{(m+l-2i)}.
\een
Hence, applying the operator
$s^{(m)}_{11}$ to the vector \eqref{singvect}
and taking the coefficient of $\xi$
we get
\beql{rhoti}
\sum_{l=0}^k c_l\ts \wt\rho^{\ts(m,\tss l)}=0
\eeq
for all $m\geqslant 0$. Our next step is to demonstrate
that this set of relations for the coefficients
of the series \eqref{muu} implies that $\rho(u)$ is
the expansion of a rational function in $u$
with the property
\beql{rhopro}
\rho(u)=u^2\rho(u^{-1}).
\eeq
To this end, introduce the coefficients $d_{-k},d_{-k+1},\dots,d_k$
by the formulas
\beql{drcr}
d_r=d_{-r}=
\begin{cases}
c_r+c_{r+2}+\cdots+c_k \qquad &\text{if $k-r$ is even},\\
c_r+c_{r+2}+\cdots+c_{k-1} \qquad &\text{if $k-r$ is odd},
\end{cases}
\eeq
where $r=0,1,\dots,k$. For any $m\geqslant k$ the relation
\eqref{rhoti} takes the form
\ben
\sum_{l=0}^k c_l\ts \big(\rho^{(m+l)}+\rho^{(m+l-2)}+\dots+
\rho^{(m-l)}\big)=0,
\een
which can written as
\ben
\sum_{r=-k}^k d_r\ts \rho^{(m+r)}=0.
\een
Therefore, we have
\begin{multline}\label{rhopolin}
\big(d_k u^k+d_{k-1}u^{k-1}+\dots+d_{-k}u^{-k}\big)\ts
\rho(u)\\
=\sum_{r=-k+1}^k \big(d_k\rho^{(k-r)}+d_{k-1}\rho^{(k-r-1)}
+\dots+d_r\rho^{(0)}\big)\ts u^r.
\end{multline}
This shows that
$\rho(u)$
is a rational function in $u$. The property
\eqref{rhopro} is equivalent to the relations
\ben
d_k\rho^{(k-r)}+d_{k-1}\rho^{(k-r-1)}
+\dots+d_r\rho^{(0)}
=d_k\rho^{(k+r-2)}+d_{k-1}\rho^{(k+r-3)}
+\dots+d_{-r+2}\rho^{(0)}
\een
for $r=2,3,\dots,k+1$, where we assume that $\rho^{(i)}$
and $d_i$ with out-of-range indices are zero.
However, the relations are easily verified by using
\eqref{drcr}: after writing them in terms of the coefficients
$c_l$ they take the form of \eqref{rhoti} with $m=r-2$.

The argument is now completed by
noting that the series $\mu(u)$ and $\mu'(u)$
are expressed in terms of $\rho(u)$ as
\ben
\mu(u)=\frac{q\ts (1-\rho(u))}{u^{-2}-q^2},\qquad
\mu'(u)=\frac{\rho(u)u^{-2}-q^2}{u^{-2}-q^2}.
\een
Hence, using \eqref{eigen}, we obtain
\beql{rhoze}
\frac{\bar\mu(u^{-1})}{\mu(u)}=
\frac{1-u^{-2}\rho(u)}{1-\rho(u)}.
\eeq
Now write \eqref{rhopolin} as $D(u)\rho(u)=F(u)$ so that
$F(u)$ and $D(u)$ are Laurent polynomials in $u$
with $D(u^{-1})=D(u)$. Moreover, $u^2 F(u^{-1})=F(u)$
due to \eqref{rhopro}.
Hence, \eqref{rhoze} implies
\ben
\frac{\bar\mu(u^{-1})}{\mu(u)}=
\frac{D(u)-u^{-2}F(u)}
{D(u)-F(u)}.
\een
Recalling that $d_k=d_{-k}=c_k\ne 0$ set
\ben
Q(u)=d_k^{-1} u^k \big(D(u)-F(u)\big)=(1-\rho^{(0)})u^{2k}+
\dots+1.
\een
The coefficient $1-\rho^{(0)}$ is nonzero by
\eqref{rhoze}
and Proposition~\ref{prop:nonzero}.
Hence $Q(u)$ is a polynomial in $u$ of degree $2k$
with the constant term equal to $1$ and
$Q(u)$ satisfies \eqref{ratfun}.
\epf

Proposition~\ref{prop:polyntw} implies that
if $\dim V(\mu(u);\bar\mu(u))<\infty$ then
there exist nonzero constants $\ga_1,\dots,\ga_{2k}$
such that
\beql{gacond}
\frac{\bar\mu(u^{-1})}{\mu(u)}=
\frac{(u+\ga_1)\dots(u+\ga_{2k})}{(1+\ga_1u)\dots(1+\ga_{2k}u)}.
\eeq
Therefore, in order to determine which
representations $V(\mu(u);\bar\mu(u))$
are finite-dimensional, we may restrict our attention
to those whose highest weights satisfy \eqref{gacond}.
We now aim to prove a tensor product decomposition
for such representations analogous to Proposition~\ref{prop:tpr};
cf. \cite[Prop.~4.3.2]{m:yc}.

We will re-enumerate
the numbers $\ga_i$, if necessary, so that
for each $i=1,\dots,k$
the following condition holds:
if the multiset $\{\ga_r\ga_s\ |\ 2i-1\leqslant r<s \leqslant 2k\}$
contains numbers of the form $q^{-2m}$ with
nonnegative integers $m$,
then $\ga_{2i-1}\ga_{2i}=q^{-2m^{}_0}$ and $m^{}_0$
is minimal amongst these nonnegative integers.

Assuming that the $\ga_i$ satisfy these conditions,
let us choose square roots $\al_i$ and $\be_i$ so that
\beql{albega}
\al_i^2=\ga_{2i-1}^{-1},\qquad \be_i^2=\ga^{}_{2i},
\qquad i=1,2,\dots,k.
\eeq
Recall the evaluation modules $L(\al,\be)$
over the algebra $\U_q(\wh\gl_2)$
defined in Sec.~\ref{subsec:qafftwo}. Each of them
may also be regarded as a module over the extended
algebra $\U^{\tss\text{\rm ext}}_q(\wh\gl_2)$
via the epimorphism \eqref{uextepi}
so that
the elements $t_{11}^{(0)}\ts \bar t_{11}^{\ts(0)}$
and $t_{22}^{(0)}\ts \bar t_{22}^{\ts(0)}$
act as the identity operators. More generally, the
tensor product
\beql{tenprtwe}
L(\al_1,\be_1)\ot\dots\ot L(\al_k,\be_k)
\eeq
can be regarded as a representation of the algebra
$\U^{\tss\text{\rm ext}}_q(\wh\gl_2)$ and hence
as a representation over its subalgebra $\Y'_q(\spa_2)$.
In other words, as far as the action of
$\Y'_q(\spa_2)$ on the space \eqref{tenprtwe} is concerned,
the operators $s_{ij}(u)$ are related with
the action of the generators of
the algebra $\U_q(\wh\gl_2)$ by the formulas \eqref{sumatrs}.

\bpr\label{prop:tprtw}
If the above condition on the parameters
$\ga_i$ holds, then there exists an automorphism
of the algebra $\Y'_q(\spa_2)$ of the form
\eqref{multautotw} such that its composition
with the representation $V(\mu(u);\bar\mu(u))$
is isomorphic to the representation
\eqref{tenprtwe}
of $\Y'_q(\spa_2)$.
\epr

\bpf
Due to \eqref{sumatrs}
the generators of $\Y'_q(\spa_2)$ act on
the tensor product module \eqref{tenprtwe}
by the formulas
\beql{svt}
\bal
s_{11}(u)&=q\ts t_{11}(u)\ts \bar t_{12}(u^{-1})
-t_{12}(u)\ts \bar t_{11}(u^{-1}),\\
s_{21}(u)&=q\ts t_{21}(u)\ts \bar t_{12}(u^{-1})
-t_{22}(u)\ts \bar t_{11}(u^{-1}),\\
\bar s_{21}(u)&=q\ts \bar t_{21}(u)\ts  t_{12}(u^{-1})
-\bar t_{22}(u)\ts t_{11}(u^{-1}).
\eal
\eeq
Consider
the tensor product $\ze=\ze_1\ot\dots\ot\ze_k$
of the highest vectors
of the representations
$L(\al_i,\be_i)$. As we pointed out in Sec.~\ref{subsec:qafftwo},
this vector generates a highest weight submodule
of the tensor product module \eqref{tenprtwe}
over $\U_q(\wh\gl_2)$.
Therefore, the formulas \eqref{hwalbe} and \eqref{svt} imply that
\beql{hwcsp}
\bal
s_{11}(u)\ts\ze&=0,\\
s_{21}(u)\ts\ze&=-\prod_{i=1}^k(\al_i^{-1}+\al^{}_iu^{-1})
(\be^{}_i+\be_i^{-1}u^{-1})\ts\ze,\\
\bar s_{21}(u^{-1})\ts\ze&=-\prod_{i=1}^k(\al^{}_i+\al_i^{-1}u^{-1})
(\be_i^{-1}+\be^{}_iu^{-1})\ts\ze.
\eal
\eeq
Hence the ratio of the eigenvalues of $\bar s_{21}(u^{-1})$
and $s_{21}(u)$ equals
\ben
\prod_{i=1}^k \frac{(\al^{}_i+\al_i^{-1}u^{-1})
(\be_i^{-1}+\be^{}_iu^{-1})}{(\al_i^{-1}+\al^{}_iu^{-1})
(\be^{}_i+\be_i^{-1}u^{-1})}
=\prod_{i=1}^k \frac{(u+\al_i^{-2})
(u+\be^2_i)}{(1+\al_i^{-2}u)
(1+\be_i^{2}u)}=\prod_{i=1}^{2k} \frac{u+\ga_i}{1+\ga_iu},
\een
which coincides with $\bar\mu(u^{-1})/\mu(u)$ by \eqref{gacond}.
We may conclude that there exists an automorphism
of the algebra $\Y'_q(\spa_2)$ of the form
\eqref{multautotw} such that its composition
with the representation $V(\mu(u);\bar\mu(u))$
is isomorphic to the irreducible quotient
of the cyclic span $\Y'_q(\spa_2)\ts\ze$.
In order to complete the argument,
we will now be proving that
$\Y'_q(\spa_2)$-module \eqref{tenprtwe}
is irreducible and so it coincides with the cyclic span of $\ze$.

Denote the representation \eqref{tenprtwe}
of $\Y'_q(\spa_2)$ by $L$.
We first prove the following claim: any vector $\xi\in L$
satisfying $s_{11}(u)\tss\xi=0$ is proportional to
$\ze$.
We use the induction on $k$ and suppose that $k\geqslant 1$.
Write
\ben
\xi=\sum_{r=0}^p(t_{21})^r \ze_1\ot \xi_r,
\qquad
\text{where}
\quad
\xi_r\in L(\al_2,\be_2)\ot\dots\ot L(\al_k,\be_k)
\een
and $p$ is some nonnegative integer.
Moreover,
if $\al_1/\be_1=\pm\tss q^m$ for some nonnegative
integer $m$, then we will assume that
$p\leqslant m$.
We will also assume that $\xi_p\ne 0$. Using
\eqref{sumatrs} and the coproduct
formulas \eqref{copraff}, we get
\ben
\bal
s_{11}(u)\big((t_{21})^r \ze_1\ot \xi_r\big)&=
t_{11}(u)\ts \bar t_{11}(u^{-1})(t_{21})^r
\ze_1\ot s_{11}(u)\ts\xi_r\\
{}&+t_{11}(u)\ts \bar t_{12}(u^{-1})(t_{21})^r
\ze_1\ot s_{12}(u)\ts\xi_r\\
{}&+t_{12}(u)\ts \bar t_{11}(u^{-1})(t_{21})^r
\ze_1\ot s_{21}(u)\ts\xi_r\\
{}&+t_{12}(u)\ts \bar t_{12}(u^{-1})(t_{21})^r
\ze_1\ot s_{22}(u)\ts\xi_r.
\eal
\een
Now use relations \eqref{toneone} and \eqref{tonetwo}
together with the following formulas which are implied by
\eqref{eval}:
\ben
\bar t_{11}(u^{-1})(t_{21})^r \ze_1
=(\bar t_{11}+t_{11}u^{-1})(t_{21})^r \ze_1
=\big(q^{r}\al^{-1}_1+q^{-r}\al^{}_1\tss u^{-1}\big)(t_{21})^r \ze_1,
\een
and
\ben
\bar t_{12}(u^{-1})(t_{21})^r \ze_1
=\bar t_{12}\tss(t_{21})^r \ze_1=
(q^r-q^{-r})\big(q^{r-1}\be_1/\al_1-q^{-r+1}\al_1/\be_1\big)
(t_{21})^{r-1} \ze_1.
\een
Taking the coefficient of $(t_{21})^p \ze_1$ in
the expansion of $s_{11}(u)\tss\xi$ we get
$s_{11}(u)\tss\xi_p=0$. By the induction hypothesis,
applied to the representation
$
L(\al_2,\be_2)\ot\dots\ot L(\al_k,\be_k),
$
the vector
$\xi_p$ must be proportional to $\ze_2\ot\dots\ot\ze_k\tss$.
As we observed in Sec.~\ref{subsec:qafftwo},
the cyclic $\U_q(\wh\gl_2)$-span of the vector $\xi_p$
is a highest weight representation of $\U_q(\wh\gl_2)$
whose highest weight is found by formulas \eqref{hwalbe}.
Hence, using \eqref{sumatrs}, we find that
\ben
s_{21}(u)\ts\xi_p=\nu(u)\ts\xi_p,\qquad
\bar s_{21}(u)\ts\xi_p=\bar\nu(u)\ts\xi_p,
\qquad s_{12}(u)\ts\xi_p=\nu'(u)\ts\xi_p,
\een
where
\ben
\bal
\nu(u)&=-\prod_{i=2}^k(\al_i^{-1}+\al^{}_iu^{-1})
(\be^{}_i+\be_i^{-1}u^{-1}),\\
\bar\nu(u^{-1})&=-\prod_{i=2}^k(\al^{}_i+\al_i^{-1}u^{-1})
(\be_i^{-1}+\be^{}_iu^{-1}),
\eal
\een
and
\ben
\nu'(u)
=\frac{(q^2-1)\tss\nu(u)+(1-u^2q^2)\ts\bar\nu(u^{-1})}
{q(u^2-1)}.
\een
Note that these are polynomials in $u^{-1}$.
To complete the proof of the claim, we need to show that $p=0$.
Suppose on the contrary that $p\geqslant 1$. Then
taking the coefficient of
$(t_{21})^{p-1} \ze_1$ in the expansion
of $s_{11}(u)\tss\xi$ we get
\beql{relxip}
\bal
&\big(q^{-p+1}\al^{}_1+q^{p-1}\al^{-1}_1 u^{-1}\big)
\big(q^{p-1}\al^{-1}_1+q^{-p+1}\al^{}_1 u^{-1}\big)
\ts s_{11}(u)\tss\xi_{p-1}\\
{}+{}&(q^p-q^{-p})\big(q^{p-1}\be_1/\al_1-q^{-p+1}\al_1/\be_1\big)\\
{}&\times
\Big(\big(q^{-p+1}\al^{}_1+q^{p-1}\al^{-1}_1 u^{-1}\big)\ts\nu'(u)
+u^{-1}\big(q^{p}\al^{-1}_1+q^{-p}\al^{}_1 u^{-1}\big)\ts\nu(u)\Big)
\ts \xi_p=0.
\eal
\eeq
By the definition of the action of the algebra $\Y'_q(\spa_2)$
on the vector space \eqref{tenprtwe},
the expression
$s_{11}(u)\tss\xi_{p-1}$ is a polynomial
in $u^{-1}$. Now we consider two cases. Suppose first
that the expression $q^{p}\al^{-1}_1+q^{-p}\al^{}_1 u^{-1}$
does not vanish at $u=-q^{2p-2}\al^{-2}_1$. Then
putting this value of $u$ into \eqref{relxip}
and recalling the notation \eqref{albega} we get
the relation
\ben
\big(\ga_1\ga_2-q^{-2p+2}\big)\big(\ga_1\ga_3-q^{-2p+2}\big)
\dots \big(\ga_1\ga_{2k}-q^{-2p+2}\big)=0.
\een
However, this is impossible due to the conditions on the parameters
$\ga_i$. Thus, in the case under consideration, $p$ must be zero.

Now suppose that the expression $q^{p}\al^{-1}_1+q^{-p}\al^{}_1 u^{-1}$
vanishes at $u=-q^{2p-2}\al^{-2}_1$
so that $\ga_1=\al_1^{-2}=\ve\ts q^{-2p+1}$ for some $\ve\in\{-1,1\}$.
In this case we may simplify \eqref{relxip} by canceling
the common factor $q^{-p+1}\al^{}_1+q^{p-1}\al^{-1}_1 u^{-1}$.
Then setting $u=-\ve\ts q$ we obtain
\ben
\nu'(-\ve\ts q)-\qin \nu(-\ve\ts q)=-(q+\qin)\ts\bar\nu(-\ve\ts\qin)=0
\een
which gives the relation
\ben
\big(\ga_3-\ve\ts q\big)\big(\ga_4-\ve\ts q\big)
\dots \big(\ga_{2k}-\ve\ts q\big)=0.
\een
Hence, $\ga_1\ga_j=q^{-2p+2}$ for some $j\in\{3,\dots,2k\}$
which contradicts the condition of the $\ga_i$.
Thus, $p$ must be zero is this case as well, and the claim
is proved.

Suppose that $M$ is a nonzero submodule of $L$.
Then $M$ must contain a nonzero vector $\xi$ such that
$s_{11}(u)\tss\xi=0$. By the claim proved above,
$\xi$ is proportional to the highest vector $\ze$, and so $\ze$
belongs to $M$. It remains to show that
the vector $\ze$ is cyclic in $L$, that is,
the submodule $K=\Y'_q(\spa_2)\tss\ze$ coincides with $L$.
We will do this by employing the dual space $L^*$
introduced in the proof of Proposition~\ref{prop:tpr}.
We equip $L^*$ with
a $\Y'_q(\spa_2)$-module structure
by using the anti-automorphism \eqref{antiautw}.
Namely, we set
\beql{dualstarmtw}
(y\ts\om)(\eta)=\om(\varkappa(y)\ts\eta)
\quad
\text{for}
\quad
y\in\Y'_q(\spa_2)
\Fand
\om\in L^*,\ \eta\in L.
\eeq
Since $\varkappa$
is obtained as the restriction of the
anti-automorphism \eqref{varsigma},
we conclude that \eqref{isomlstar}
is a $\Y'_q(\spa_2)$-module isomorphism.
Arguing as in the proof of Proposition~\ref{prop:tpr},
suppose now that the
submodule $K=\Y'_q(\spa_2)\tss\ze$ of $L$ is proper.
The annihilator
\ben
\text{Ann\ts} K=\{\om\in L^*\ |\ \om(\eta)=0
\quad\text{for all}\quad\eta\in K\}
\een
is a submodule of the $\Y'_q(\spa_2)$-module $L^*$, which does not
contain the vector $\ze_1^*\ot\dots\ot\ze_k^*$.
However, this contradicts the claim verified
in the first part of the proof,
because the tensor product in \eqref{isomlstar}
is associated with the set of parameters
obtained by swapping $\ga_{2i-1}$
and $\ga_{2i}$ for each $i=1,\dots,k$ so that
the condition on the parameters
remain satisfied after this swap.
\epf

Proposition~\ref{prop:tpr} allows us to
describe the finite-dimensional irreducible
representations of the algebra $\Y'_q(\spa_2)$.

\bth\label{thm:classspwo}
The irreducible highest weight representation
$V\big(\mu(u);\bar\mu(u)\big)$
of $\Y'_q(\spa_2)$ is
finite-dimensional if and only if there exists a
polynomial $P(u)$ in $u$ of even degree
with constant term $1$ such that
$u^{\deg P}\ts P(u^{-1})=q^{-\deg P}\ts P(u\tss q^2)$ and
\beql{dponetw}
\frac{\bar\mu(u^{-1})}{\mu(u)}=q^{-\deg P}\cdot
\frac{P(u\tss q^2)}{P(u)}.
\eeq
In this case $P(u)$ is unique.
\eth

\bpf
Suppose that the representation
$V\big(\mu(u);\bar\mu(u)\big)$
is finite-dimensional.
By Proposition~\ref{prop:polyntw}
there exist constants $\ga_i$ such that \eqref{gacond}
holds. Re-enumerate these constants to satisfy
the assumptions of Proposition~\ref{prop:tprtw}.
This proposition implies that
each representation $L(\al_i,\be_i)$
occurring in \eqref{tenprtwe} is finite-dimensional.
Therefore, all ratios $\al_i/\be_i$ must have the form
$\pm q^{m_i}$, where each $m_i$ is a nonnegative integer.
Then the polynomial
\ben
\bal
P(u)&=\prod_{i=1}^k (1+\be_i^2\tss u)(1+\be_i^2\tss q^2\tss u)\dots
(1+\be_i^2\tss q^{2m_i-2}\tss u)\\
{}&\times
\prod_{i=1}^k (1+\al_i^{-2}\tss u)(1+\al_i^{-2}\tss q^2\tss u)\dots
(1+\al_i^{-2}\tss q^{2m_i-2}\tss u)
\eal
\een
has the property
$u^{\deg P}\ts P(u^{-1})=q^{-\deg P}\ts P(u\tss q^2)$
and satisfies \eqref{dponetw}.

Conversely, suppose \eqref{dponetw} holds for
a polynomial $P(u)=(1+\ga_1u)\dots(1+\ga_{2k}u)$
with the property
$u^{\deg P}\ts P(u^{-1})=q^{-\deg P}\ts P(u\tss q^2)$.
This property implies that
the multiset of parameters $\ga_i$
can be written in the form
\ben
\{\ga_1,\dots,\ga_{2k}\}
=\{\al_1,\dots,\al_k,\al^{-1}_1\tss q^{-2},\dots,
\al^{-1}_k\tss q^{-2}\}.
\een
Consider
the irreducible highest weight representation
$L(\nu_1(u),\nu_2(u);\bar\nu_1(u),\bar\nu_2(u))$
of the algebra $\U_q(\wh\gl_2)$, where
the components of the highest weight are given by
\ben
\bal
\nu_1(u)&=
(\al^{}_1\tss q+\qin u^{-1})\dots
(\al^{}_k\tss q+\qin u^{-1}),
\\
\nu_2(u)&=
(\al^{}_1+u^{-1})\dots
(\al^{}_k+u^{-1}),
\\
\bar\nu_1(u)&=
(\qin+\al^{}_1q\ts u)\dots (\qin+\al^{}_kq\ts u),
\\
\bar\nu_2(u)&=
(1+\al^{}_1u)\dots (1+\al^{}_ku).
\eal
\een
By Theorem~\ref{thm:classgltwo}, the representation
$L(\nu_1(u),\nu_2(u);\bar\nu_1(u),\bar\nu_2(u))$
is finite-dimensional as
\ben
\frac{\nu_1(u)}{\nu_2(u)}=q^{-\deg Q}\cdot
\frac{Q(u\tss q^2)}{Q(u)}
=\frac{\bar\nu_1(u)}{\bar\nu_2(u)}
\een
with $Q(u)=(1+\al_1u)\dots(1+\al_{k}u)$.
We will regard $L(\nu_1(u),\nu_2(u);\bar\nu_1(u),\bar\nu_2(u))$
as a representation of the extended algebra
$\U^{\tss\text{\rm ext}}_q(\wh\gl_2)$ via the epimorphism
\eqref{uextepi}. The formulas \eqref{svt} imply that
the cyclic span $\Y'_q(\spa_2)\ts\ze$ of the highest
vector $\ze$ of $L(\nu_1(u),\nu_2(u);\bar\nu_1(u),\bar\nu_2(u))$
is a highest weight representation
of $\Y'_q(\spa_2)$ with the highest weight $(\la(u);\bar\la(u))$,
where
\ben
\la(u)=-\nu_2(u)\tss\bar\nu_1(u^{-1}),
\qquad
\bar\la(u)=-\bar\nu_2(u)\tss\nu_1(u^{-1}).
\een
Therefore, the irreducible highest weight representation
$V(\la(u);\bar\la(u))$ is finite-dimensional and
\ben
\frac{\bar\la(u^{-1})}{\la(u)}=\frac{\nu_1(u)}{\nu_2(u)}
\cdot \frac{\bar\nu_2(u^{-1})}{\bar\nu_1(u^{-1})}
=\frac{Q(u\tss q^2)\ts Q(u^{-1})}{Q(u)\ts Q(u^{-1}q^2)}
=q^{-\deg P}\cdot
\frac{P(u\tss q^2)}{P(u)}.
\een
Due to \eqref{dponetw} there exists an automorphism
of the algebra $\Y'_q(\spa_2)$ of the form \eqref{multautotw}
such that the composition of the representation
$V(\la(u);\bar\la(u))$ with this automorphism
is isomorphic to $V(\mu(u);\bar\mu(u))$. Hence
the latter is also finite-dimensional.

The uniqueness of $P(u)$ is easily verified.
\epf

In the following corollary we use the $q$-spirals
introduced in Sec.~\ref{subsec:qafftwo}.
Denote by $L$ the tensor product \eqref{tenprtwe},
where for all $i=1,\dots,k$ we have $\al_i/\be_i=\pm q^{m_i}$
for some nonnegative integers $m_i$.

\bco\label{cor:irrcritw}
The representation $L$
of $\Y'_q(\spa_2)$ is irreducible
if and only if each pair of
the $q$-spirals
\ben
S_q(\al_i,\be_i),\quad S_q(\al_j,\be_j)\Fand
S_q(\be^{-1}_i,\al^{-1}_i),\quad S_q(\al_j,\be_j)
\een
is in general position for all $1\leqslant i<j\leqslant k$.
\eco

\bpf
Suppose that the condition on the
$q$-spirals is satisfied. Then Corollary~\ref{cor:irrcrit}
implies that $L$ is irreducible as
a representation of the algebra $\U_q(\wh\gl_2)$.
Moreover, as we pointed out in the proof of that corollary,
any permutation of the tensor factors in \eqref{tenprtwe}
yields an isomorphic representation.
Hence we may assume that the nonnegative integers $m_i$
satisfy the inequalities
$m_1\leqslant\dots\leqslant m_k$.
Let us verify that in this case
the condition
of Proposition~\ref{prop:tprtw} is satisfied.
Indeed, if this is not the case, then
$\ga_r\ga_s=q^{-2p}$ for some $2i-1\leqslant r<s \leqslant 2k$
and a nonnegative integer $p$ such that $p<m_i$.
Suppose first that $r$ and $s$ are both odd.
Then we may assume that $r=2i-1$ and $s=2j-1$ for some $j>i$.
By \eqref{albega} we have $\ga_{2i-1}=\al_i^{-2}$ and
$\ga_{2j-1}=\al_j^{-2}$. Hence, $\al_j=\pm \al_i^{-1}q^p$
which means that $\al_j$ belong to the $q$-spiral
$S_q(\be^{-1}_i,\al^{-1}_i)$. However, the condition
$m_i\leqslant m_j$ then implies that the $q$-spirals
$S_q(\be^{-1}_i,\al^{-1}_i)$ and $S_q(\al_j,\be_j)$
are not in general position. This contradicts the assumptions
of the proposition. The remaining cases, where $r$ or $s$ is
even lead to similar contradictions. Thus,
Proposition~\ref{prop:tprtw} allows us to
conclude that $L$ is irreducible as a representation
of $\Y'_q(\spa_2)$.

Conversely, suppose that the representation $L$ of $\Y'_q(\spa_2)$
is irreducible. Then $L$ is irreducible as
a representation of $\U_q(\wh\gl_2)$. By Corollary~\ref{cor:irrcrit}
the $q$-spirals $S_q(\al_i,\be_i)$ and $S_q(\al_j,\be_j)$
are in general position for all $i<j$.
Now fix an index $i\in\{1,\dots,k\}$ and consider
the $\Y'_q(\spa_2)$-module $L'$ obtained by replacement
of the tensor factor
$L(\al_i,\be_i)$ by
$L(\be^{-1}_i,\al^{-1}_i)$. We claim that
$L'$ is isomorphic to $L$.
Indeed, the formulas \eqref{hwcsp} show that
the highest
weight of the cyclic $\Y'_q(\spa_2)$-span of the tensor
product of the highest vectors of the tensor factors
occurring in $L'$ is unchanged under the
replacement $\al_i\mapsto \be_i^{-1}$, $\be_i\mapsto \al_i^{-1}$.
This implies that the module $L$ is isomorphic
to the irreducible quotient of this span.
Since $\dim L=\dim L'$, the claim follows.

Thus, $L'$ is irreducible as a $\Y'_q(\spa_2)$-module and,
hence, as a $\U_q(\wh\gl_2)$-module.
By Corollary~\ref{cor:irrcrit},
the $q$-spiral $S_q(\be^{-1}_i,\al^{-1}_i)$ is
in general position with any $q$-spiral $S_q(\al_j,\be_j)$ for $i\ne j$.
This gives the required condition on the $q$-spirals.
\epf

\subsection{Classification theorem}

We can now prove the classification theorem for
finite-dimensional irreducible representations
of the twisted $q$-Yangian $\Y'_q(\spa_{2n})$
for arbitrary $n\geqslant 1$.
By Theorem~\ref{thm:fdhwtw}, all finite-dimensional
irreducible representations of the twisted
$q$-Yangian $\Y'_q(\spa_{2n})$
have the form $V(\mu(u);\bar\mu(u))$ for a
certain highest weight $(\mu(u);\bar\mu(u))$.

\bth\label{thm:clgentw}
The irreducible highest weight representation
$V(\mu(u);\bar\mu(u))$
of the algebra $\Y'_q(\spa_{2n})$ is
finite-dimensional if and only if there exist polynomials
$P_1(u),\dots,P_{n}(u)$ in $u$,
all with constant term $1$,
where $P_1(u)$ is of even degree and satisfies
$u^{\deg P_1}\ts P_1(u^{-1})=q^{-\deg P_1}\ts P_1(u\tss q^2)$, and
nonzero constants $\phi_1,\dots,\phi_n$ such that
\beql{nunuponetw}
\frac{\phi_{i-1}\mu_{i-1}(u)}
{\phi_{i}\ts\mu_{i}(u)}=q^{-\deg P_i}\cdot
\frac{P_i(u\tss q^2)}{P_i(u)}
=\frac{\phi_{i-1}\bar\mu_{i-1}(u)}{\phi_{i}\ts\bar\mu_{i}(u)}
\eeq
for $i=2,\dots,n$ and
\beql{dponetwgen}
\frac{\bar\mu_1(u^{-1})}{\mu_1(u)}=q^{-\deg P_1}\cdot
\frac{P_1(u\tss q^2)}{P_1(u)}.
\eeq
The polynomials $P_1(u),\dots,P_{n}(u)$
are determined uniquely, while the tuple $(\phi_1,\dots,\phi_n)$
is determined uniquely, up to a common factor.
\eth

\bpf
Suppose first that $\dim V(\mu(u);\bar\mu(u))<\infty$.
Let $J$ be the left ideal of $\Y'_q(\spa_{2n})$ generated by all
coefficients of the series $s_{ij}(u)$ with
$i,j=1,3,\dots,2n-1$. Due to \eqref{relsbarssympl},
all coefficients of the series $\bar s_{ij}(u)$ with
$i,j=1,3,\dots,2n-1$ also belong to $J$.
Consider the subspace $V^{J}$ of $V(\mu(u);\bar\mu(u))$ defined by
\ben
V^{J}=\{\eta\in V(\mu(u);\bar\mu(u))\ |\ s_{ij}(u)\ts\eta=0
\quad\text{for all}\quad i,j=1,3,\dots,2n-1\}.
\een
Note that the highest
vector $\xi$ of $V(\mu(u);\bar\mu(u))$ belongs to $V^{J}$.
The defining relations \eqref{draffsu}
together with \eqref{rsrsbar} imply
that if the indices $i,a,b$ are odd and $j$ is even, then
\ben
s_{ia}(u)\tss s_{jb}(v),\quad
s_{ia}(u)\tss \bar s_{jb}(v)\in J.
\een
Therefore the subspace $V^{J}$ is stable under the action
of the operators $s_{2i,2a-1}(u)$ and $\bar s_{2i,2a-1}(u)$.
Moreover, regarding the relations
\eqref{draffsu} and \eqref{rsrsbar}
modulo the left ideal $J$, we find that
the mapping
\beql{mappact}
t_{ia}(u)\mapsto s_{2i,2a-1}(u),\qquad
\bar t_{ia}(u)\mapsto \bar s_{2i,2a-1}(u),\qquad i,a=1,\dots,n,
\eeq
defines an action of the algebra
$\U^{\tss\text{\rm ext}}_q(\wh\gl_n)$ on
the space $V^{J}$.
The cyclic span $\U^{\tss\text{\rm ext}}_q(\wh\gl_n)\ts\xi$
is a finite-dimensional
highest weight representation of
$\U^{\tss\text{\rm ext}}_q(\wh\gl_n)$
with the highest weight
$(\mu(u);\bar\mu(u))$. It follows from
Corollary~\ref{cor:clgenext} that the highest weight
satisfies the conditions \eqref{nunuponetw}
for appropriate nonzero constants $\phi_i$.

Furthermore, the twisted $q$-Yangian $\Y'_q(\spa_{2})$
act on $V(\mu(u);\bar\mu(u))$ via
the homomorphism
$\Y'_q(\spa_{2})\to\Y'_q(\spa_{2n})$ which sends $s_{ij}(u)$
to the series with the same name in $\Y'_q(\spa_{2n})$.
The cyclic span $\Y'_q(\spa_{2})\tss\xi$
is a highest weight representation of $\Y'_q(\spa_{2})$ with
the highest weight $(\mu_1(u);\bar\mu_1(u))$.
Its irreducible quotient
is finite-dimensional, and so \eqref{dponetwgen}
follows from Theorem~\ref{thm:classspwo}.

In order to prove the converse statement, note that
given two irreducible highest weight representations
$V(\mu(u);\bar\mu(u))$ and $V(\la(u);\bar\la(u))$
such that the components of the highest weights satisfy
the conditions \eqref{nunuponetw} and \eqref{dponetwgen}
with the same set of polynomials $P_1(u),\dots,P_n(u)$,
there exist automorphisms of the form \eqref{multautotw}
and \eqref{aautoaffsigntw} such that
the composition of the representation $V(\mu(u);\bar\mu(u))$
with these automorphisms is isomorphic to $V(\la(u);\bar\la(u))$.
Hence, it suffices to show that given any
set of polynomials $P_1(u),\dots,P_n(u)$
of the form described in the formulation of the theorem,
there exists a finite-dimensional representation
$V(\mu(u);\bar\mu(u))$ whose highest weight satisfies
\eqref{nunuponetw} and \eqref{dponetwgen} with
$\phi_i=1$ for all $i$. We will use a result from \cite[Sec.~6]{m:rtq}
concerning a particular irreducible
highest weight representation
$L(\nu)$ of the quantized enveloping algebra $\U_q(\gl_{2n})$;
see Sec.~\ref{subsec:que} above for the definition.
The highest weight $\nu$ has the form
\ben
\nu=(q^{r_n},\dots,q^{r_1},1,\dots,1),\qquad
r_n\geqslant\dots\geqslant r_1\geqslant 0,
\een
where the parameters $r_i$ are integers. The representation
$L(\nu)$ is finite-dimensional and it admits
a basis parameterized by the Gelfand--Tsetlin patterns
associated with $\nu$; see \cite{j:qr}.
As in \cite{m:rtq}
consider the pattern $\Omega^{\tss 0}$ such that
for each $k=1,2,\dots,n$ its row $2k-1$ counted from the bottom
is $(r_k,r_{k-1},\dots,r_1,0,\dots,0)$ with $k-1$ zeros, while
the row $2k$ from the bottom
is $(r_k,r_{k-1},\dots,r_1,0,\dots,0)$ with $k$ zeros.
Then the corresponding basis vector $\ze^{}_{\tss\Om^{\tss 0}}$
has the properties
\beql{annihom}
\bal
\bar t_{ij}\ts \ze^{}_{\tss\Om^{\tss 0}}&=0
\qquad\text{if $j$ is even and $i<j$,}\\
t_{ij}\ts \ze^{}_{\tss\Om^{\tss 0}}&=0
\qquad\text{if $i$ is odd and $i>j$.}
\eal
\eeq
We denote by $L(d\ts\nu)$ the composition
of the representation $L(\nu)$ with the
automorphism of $\U_q(\gl_{2n})$ given in \eqref{automd}.
We will consider $L(d\ts\nu)$ as an
evaluation module over $\U_q(\wh\gl_{2n})$ by using
the homomorphism \eqref{eval}.

Suppose now that $V(\mu(u);\bar\mu(u))$ is
a finite-dimensional highest weight representation
of $\Y'_q(\spa_{2n})$ with the highest vector $\xi$.
By the first part of the proof
we can associate a family of polynomials
$P_1(u),\dots,P_n(u)$ to $V(\mu(u);\bar\mu(u))$.
The coproduct structure on $\U_q(\wh\gl_{2n})$
given by \eqref{copraff} allows us to equip
the vector space $L(d\ts\nu)\ot V(\mu(u);\bar\mu(u))$
with a structure of a $\Y'_q(\spa_{2n})$-module so that
for the action of the generators we have
\beql{copract}
s_{ij}(u)(\eta\ot\theta)=\sum_{k,l=1}^{2n}
t_{ik}(u)\tss\bar t_{jl}(u^{-1})\ts\eta\ot
s_{kl}(u)\ts \theta,\quad \eta\in L(d\ts\nu),
\quad\theta\in V(\mu(u);\bar\mu(u)).
\eeq
Let us verify that
$\ze^{}_{\tss\Om^{\tss 0}}\ot\xi$ is the highest vector
of the $\Y'_q(\spa_{2n})$-module
$\Y'_q(\spa_{2n})(\ze^{}_{\tss\Om^{\tss 0}}\ot\xi)$.
Take $\eta=\ze^{}_{\tss\Om^{\tss 0}}$
and $\ze=\xi$ in \eqref{copract} and suppose that $j$ is odd
and $i\leqslant j$.
Using
\eqref{eval} and \eqref{annihom}, we find that
\ben
\bar t_{jl}(u^{-1})\ts \ze^{}_{\tss\Om^{\tss 0}}
=u^{-1} t_{jl}\ts \ze^{}_{\tss\Om^{\tss 0}}=0
\een
for $j>l$.
If $l$ is even and $j<l$, then by \eqref{annihom}
\ben
\bar t_{jl}(u^{-1})\ts \ze^{}_{\tss\Om^{\tss 0}}
=\bar t_{jl}\ts \ze^{}_{\tss\Om^{\tss 0}}=0.
\een
Hence, we may assume that $l=j+2p$ for a nonnegative
integer $p$; in particular, $l$ is odd.
Then the index $k$ in \eqref{copract}
may be assumed to be even as otherwise
$s_{kl}(u)\ts \xi=0$. Furthermore, if $k\leqslant i$ then
$k\leqslant j+2p=l$ so that $s_{kl}(u)\ts \xi=0$
in this case too. Therefore, we may assume that $k>i$.
In this case we have
\beql{ttbarze}
t_{ik}(u)\tss\bar t_{jl}(u^{-1})\ts\ze^{}_{\tss\Om^{\tss 0}}
=u^{-1}\bar t_{ik}\ts (\bar t_{jl}+u^{-1}\de_{jl}\ts t_{jl})
\ts\ze^{}_{\tss\Om^{\tss 0}}.
\eeq
By the defining relations \eqref{defrelgbar} we have
\ben
\bar t_{ik}\ts \bar t_{jl}=
q^{\delta_{ij}}\ts \bar t_{jl}\ts \bar t_{ik}
-(q-\qin)\ts (\de_{k<l} -\de_{j<i})
\ts \bar t_{il}\ts \bar t_{jk}.
\een
Now, if $k<l$ then $s_{kl}(u)\ts \xi=0$. Otherwise,
$\de_{k<l} -\de_{j<i}\ne 0$ only if $j<i$. But in this case
$j<k$ and $\bar t_{jk}\ts\ze^{}_{\tss\Om^{\tss 0}}=0$.
Thus, in all cases \eqref{ttbarze} is zero due to
\eqref{annihom}. A similar calculation shows that
$\bar s_{ij}(u)(\ze^{}_{\tss\Om^{\tss 0}}\ot\xi)=0$
for odd $j$ and $i\leqslant j$. By \eqref{relsbarssympl}
this proves $s_{ij}(u)(\ze^{}_{\tss\Om^{\tss 0}}\ot\xi)=0$
if $\vs(i)+\vs(j)>0$; see Definition~\ref{def:hwtw}.

Let us now calculate the eigenvalues of
$\ze^{}_{\tss\Om^{\tss 0}}\ot\xi$ with respect to
the operators $s_{2i,2i-1}(u)$ and $\bar s_{2i,2i-1}(u)$.
The above arguments show that \eqref{copract}
with $\eta=\ze^{}_{\tss\Om^{\tss 0}}$
and $\ze=\xi$ simplifies to
\ben
\bal
s_{2i,2i-1}(u)(\ze^{}_{\tss\Om^{\tss 0}}\ot\xi)
{}&=t_{2i,2i}(u)\tss\bar t_{2i-1,2i-1}(u^{-1})\ts
\ze^{}_{\tss\Om^{\tss 0}}\ot
s_{2i,2i-1}(u)\ts \xi\\
{}&=\big(d+d^{-1}u^{-1}\big)\big(d^{-1}q^{-r_i}
+d\tss q^{r_i}u^{-1}\big)\ts \mu_i(u)
(\ze^{}_{\tss\Om^{\tss 0}}\ot\xi).
\eal
\een
Similarly,
\ben
\bal
\bar s_{2i,2i-1}(u)(\ze^{}_{\tss\Om^{\tss 0}}\ot\xi)
{}&=\bar t_{2i,2i}(u)\tss t_{2i-1,2i-1}(u^{-1})\ts
\ze^{}_{\tss\Om^{\tss 0}}\ot
\bar s_{2i,2i-1}(u)\ts \xi\\
{}&=\big(d^{-1}+d\tss u\big)\big(d\tss q^{r_i}
+d^{-1}\tss q^{-r_i}\big)\ts \bar\mu_i(u)
(\ze^{}_{\tss\Om^{\tss 0}}\ot\xi).
\eal
\een
The cyclic span $\Y'_q(\spa_{2n})(\ze^{}_{\tss\Om^{\tss 0}}\ot\xi)$
is finite-dimensional. By the above formulas,
the irreducible quotient of this representation
of $\Y'_q(\spa_{2n})$ corresponds to
the family of polynomials $Q_1(u)P_1(u),\dots,Q_n(u)P_n(u)$,
where
\ben
Q_i(u)=(1+d^{-2}q^{-2r_i}u)(1+d^{-2}q^{-2r_i+2}u)
\dots (1+d^{-2}q^{-2r_{i-1}-2}u),\qquad i=2,\dots,n,
\een
and
\ben
\bal
Q_1(u)&=(1+d^2 u)(1+d^2q^2 u)\dots (1+d^2q^{2r_1-2} u)\\
{}&\times(1+d^{-2}q^{-2r_1}u)(1+d^{-2}q^{-2r_1+2}u)
\dots (1+d^{-2}q^{-2}u).
\eal
\een
Thus,
starting from the trivial representation
$V(\mu(u);\bar\mu(u))$ and choosing appropriate
parameters $d$ and $r_i$ we will be able to produce
a finite-dimensional highest weight representation
of $\Y'_q(\spa_{2n})$ associated with an arbitrary family of
polynomials $P_1(u),\dots,P_n(u)$ by iterating
this construction; cf. the proof of Theorem~\ref{thm:clgen}.
The last statement of the theorem is easily verified.
\epf

We will call $P_1(u),\dots,P_n(u)$ the {\it Drinfeld polynomials}
of the finite-dimensional representation $V(\mu(u);\bar\mu(u))$.

We will now use Theorem~\ref{thm:clgentw} to describe
finite-dimensional irreducible representations of
the quotient algebra $\Y^{\rm tw}_q(\spa_{2n})$
of $\Y'_q(\spa_{2n})$
by the relations \eqref{sqcubeaff}; see Remark~\ref{rem:spextaff}.
Every finite-dimensional
irreducible representation of the algebra $\Y^{\rm tw}_q(\spa_{2n})$
is isomorphic to the highest weight representation
$V(\mu(u);\bar\mu(u))$ which is defined in the same way as for
the algebra $\Y'_q(\spa_{2n})$; see Definition~\ref{def:hwtw}.
This time the constant terms of the series \eqref{muitw}
should satisfy the conditions $\mu_i^{(0)}\bar\mu_i^{(0)}=1$
for $i=1,\dots,n$.

\bco\label{cor:clgetw}
The irreducible highest weight representation
$V(\mu(u);\bar\mu(u))$
of the algebra $\Y^{\rm tw}_q(\spa_{2n})$ is
finite-dimensional if and only if there exist polynomials
$P_1(u),\dots,P_{n}(u)$ in $u$,
all with constant term $1$,
where $P_1(u)$ is of even degree and satisfies
$u^{\deg P_1}\ts P_1(u^{-1})=q^{-\deg P_1}\ts P_1(u\tss q^2)$ such that
\beql{nunuponetwtw}
\frac{\ve_{i-1}\mu_{i-1}(u)}
{\ve_{i}\ts\mu_{i}(u)}=q^{-\deg P_i}\cdot
\frac{P_i(u\tss q^2)}{P_i(u)}
=\frac{\ve_{i-1}\bar\mu_{i-1}(u)}{\ve_{i}\ts\bar\mu_{i}(u)}
\eeq
for $i=2,\dots,n$ and
\beql{dponetwgentw}
\frac{\bar\mu_1(u^{-1})}{\mu_1(u)}=q^{-\deg P_1}\cdot
\frac{P_1(u\tss q^2)}{P_1(u)}
\eeq
for some $\ve_i\in\{-1,1\}$.
The polynomials $P_1(u),\dots,P_{n}(u)$
are determined uniquely, while the tuple $(\ve_1,\dots,\ve_n)$
is determined uniquely, up to a simultaneous change of sign.
\eco

\bpf
Suppose that $\dim V(\mu(u);\bar\mu(u))<\infty$.
We argue as in the proof of Theorem~\ref{thm:clgentw}.
The first part of that proof
is now modified so that the mapping \eqref{mappact}
defines an action of the algebra
$\U_q(\wh\gl_n)$ on
the corresponding space $V^{J}$. The necessary conditions
on the components of the highest weight
come from the application of Theorem~\ref{thm:clgen}.

Conversely, suppose that conditions \eqref{nunuponetwtw}
and \eqref{dponetwgentw} hold. Using the natural
epimorphism $\Y'_q(\spa_{2n})\to\Y^{\rm tw}_q(\spa_{2n})$
we may regard $V(\mu(u);\bar\mu(u))$ as a
$\Y'_q(\spa_{2n})$-module. This module is
finite-dimensional by Theorem~\ref{thm:clgentw}.
\epf

We conclude with a discussion of a particular class
of representations of the twisted $q$-Yangians
associated with the evaluation homomorphisms.
By \cite[Theorem~3.15]{mrs:cs} there exists a
homomorphism
$\Y^{\rm tw}_q(\spa_{2n})\to\U^{\rm tw}_q(\spa_{2n})$
which is identical on the subalgebra $\U^{\rm tw}_q(\spa_{2n})$.
The arguments used for the proof of that theorem
apply to the algebra $\Y'_q(\spa_{2n})$ without
any changes
so that we have the homomorphism
$\Y'_q(\spa_{2n})\to\U'_q(\spa_{2n})$ given by
\beql{evaltesp}
S(u)\mapsto S+q\tss u^{-1}\tss \overline S.
\eeq
It allows one to extend any representation of
$\U'_q(\spa_{2n})$ to the twisted $q$-Yangian
$\Y'_q(\spa_{2n})$. Consider the
highest weight representations $V(\mu;\mu')$ defined
in Sec.~\ref{subsec:twrep}. The $\Y'_q(\spa_{2n})$-modules
$V(\mu;\mu')$ will be called the {\it evaluation modules}.

Suppose that this representation is finite-dimensional
with the parameters $p_i$ as defined in
Proposition~\ref{prop:fdco}.

\bpr\label{prop:eval}
The Drinfeld polynomials
of the evaluation module $V(\mu;\mu')$ over
the algebra $\Y'_q(\spa_{2n})$
are given by
\ben
P_1(u)=(1+q\tss u)(1+q^3\tss u)\dots (1+q^{2p_1-1}\tss u)
(1+q^{-2p_1-1}\tss u)(1+q^{-2p_1+1}\tss u)\dots (1+q^{-3}\tss u)
\een
and
\ben
P_i(u)=(1+q^{-2p_i-1}\tss u)(1+q^{-2p_i+1}\tss u)\dots
(1+q^{-2p_{i-1}-3}\tss u)
\een
for $i=2,\dots,n$. The parameters $\phi_i$ are found by
\ben
\phi_i=\mu_i^{-1} q^{-p_i},\qquad i=1,\dots,n.
\een
\epr

\bpf
The highest vector of the representation
$V(\mu;\mu')$ of $\U'_q(\spa_{2n})$ is also
the highest vector of the evaluation module over $\Y'_q(\spa_{2n})$.
The claims are now verified by calculating the highest
weight of the $\Y'_q(\spa_{2n})$-module $V(\mu;\mu')$
with the use of \eqref{evaltesp} and the formulas
relating the matrix elements of the matrices $S$ and $\overline S$;
cf. \cite[(2.52)]{mrs:cs}. The components
of the highest weight are found by
\ben
\mu_i(u)=\mu^{}_i-u^{-1}\mu'_i,\qquad
\bar\mu_i(u)=\frac{1+q\tss u}{u+q}\ts(u\mu^{}_i-\mu'_i),\qquad
i=1,\dots,n.
\een
Together with \eqref{nunuponetw}
and \eqref{dponetwgen}
this implies all the statements.
\epf

As we pointed out in the proof of Proposition~\ref{prop:fdco},
if the highest weight $(\mu;\mu')$ satisfies
the additional conditions $\mu_i\ts\mu'_i=-q$ for $i=1,\dots,n$, then
$V(\mu;\mu')$ can be regarded as a representation
of the quotient algebra $\Y^{\rm tw}_q(\spa_{2n})$.
In this case we have $\mu_i=\ve_i\ts q^{-p_i}$
for all $i$ and some $\ve_i\in\{-1,1\}$.
The corresponding evaluation module over
$\Y^{\rm tw}_q(\spa_{2n})$ has the same Drinfeld polynomials
as given in Proposition~\ref{prop:eval}, while
$\phi_i=\ve_i$ for all $i$.

\bigskip

{\small\sc
\noindent
Lucy Gow:
Max-Planck-Institut f\"{u}r Gravitationsphysik\newline
Albert-Einstein-Institut\newline
Am M\"{u}hlenberg 1, 14476 Potsdam,
Germany

\noindent
E-mail address: {\rm lucy.gow@aei.mpg.de}

\bigskip

\noindent
Alexander Molev: School of Mathematics and Statistics\newline
University of Sydney, NSW 2006, Australia

\noindent
E-mail address: {\rm alexm@maths.usyd.edu.au}
}

\end{document}